\documentclass[12pt]{amsart}
\usepackage{amsmath}%
\usepackage{amsfonts}%
\usepackage{amssymb}%
\usepackage{graphicx}
\usepackage[draft]{todonotes}
\usepackage[margin=1.00in]{geometry}
\usepackage{enumerate}
\usepackage{csquotes}
\usepackage{url}
\newtheorem{theorem}{Theorem}
\theoremstyle{plain}

\newtheorem{corollary}{Corollary}

\newtheorem{definition}{Definition}

\newtheorem{lemma}{Lemma}

\newtheorem{proposition}{Proposition}

\numberwithin{equation}{section}

\DeclareMathOperator{\GL}{GL}
\DeclareMathOperator{\PGL}{PGL}
\DeclareMathOperator{\SL}{SL}

\DeclareMathOperator{\Mat}{Mat}
\DeclareMathOperator{\trace}{trace}
\DeclareMathOperator{\Aut}{Aut}

\begin{document}

\title{Local Conjugacy in $\GL_2(\mathbb{Z}/p^2\mathbb{Z})$}
\author{SPUR Final Paper, Summer 2016, Revised August 2017\\~\\
Hyun Jong Kim\\
Mentor: Atticus Christensen\\
Project Suggested by Andrew Sutherland}
\pagestyle{plain}
\date{\today}
\begin{abstract}
Subgroups $H_1$ and $H_2$ of a group $G$ are said to be locally conjugate if there is a bijection $f: H_1 \rightarrow H_2$ such that $h$ and $f(h)$ are conjugate in $G$ for every $h \in H_1$. This paper studies local conjugacy among subgroups of $\GL_2(\mathbb{Z}/p^2\mathbb{Z})$, where $p$ is an odd prime, building on Sutherland's categorizations of subgroups of $\GL_2(\mathbb{Z}/p\mathbb{Z})$ and local conjugacy among them. There are two conditions that locally conjugate subgroups $H_1$ and $H_2$ of $\GL_2(\mathbb{Z}/p^2\mathbb{Z})$ must satisfy: letting $\varphi: \GL_2(\mathbb{Z}/p^2\mathbb{Z}) \rightarrow \GL_2(\mathbb{Z}/p\mathbb{Z})$ be the natural homomorphism, $H_1 \cap \ker \varphi$ and $H_2 \cap \ker \varphi$ must be locally conjugate in $\GL_2(\mathbb{Z}/p^2\mathbb{Z})$ and $\varphi(H_1)$ and $\varphi(H_2)$ must be locally conjugate in $\GL_2(\mathbb{Z}/p\mathbb{Z})$. To identify $H_1$ and $H_2$ up to conjugation, we choose $\varphi(H_1)$ and $\varphi(H_2)$ to be similar to each other, then understand the possibilities for $H_1 \cap \ker \varphi$ and $H_2 \cap \ker \varphi$. This study fully categorizes local conjugacy in $\GL_2(\mathbb{Z}/p^2\mathbb{Z})$ through such casework.
\end{abstract}

\maketitle
\tableofcontents
\newpage

\section{Introduction}
	Given an elliptic curve $E$ over a number field $K$ and an integer $n$, the action of the absolute Galois group of $K$ on the $n$-torsion subgroup $E[n] \simeq (\mathbb{Z}/n\mathbb{Z})^2$ determines a subgroup of $\Aut(E[n]) \simeq \GL_2(\mathbb{Z}/n\mathbb{Z})$. Sutherland \cite{Sutherland} gives an efficient algorithm for computing the images of the Galois representation associated to an elliptic curve in $\GL_2(\mathbb{Z}/p\mathbb{Z})$ that determines the image up to local conjugacy. One needs to understand locally conjugate subgroups to determine the image up to conjugation. \cite{Sutherland} categorizes local conjugacy of subgroups in $\GL_2(\mathbb{Z}/p\mathbb{Z})$, see Theorem \ref{theoremLocConjGL2p}. A categorization of local conjugacy of subgroups in $\GL_2(\mathbb{Z}/p^2\mathbb{Z})$ is a first step in extending the categorization in $\GL_2(\mathbb{Z}/p\mathbb{Z})$ to the full $p$-adic image in $\GL_2(\mathbb{Z}_p)$. \par
	Locally conjugate subgroups $H_1$ and $H_2$ of a group $G$ form what is known as a Gassman triple $(G,H_1,H_2)$, which arise in the study of arithmetically equivalent number fields.  Such number fields have the 
same Dedekind zeta function but need not be isomorphic. Gassmann triples of the form $(\GL_2(\mathbb{Z}/n\mathbb{Z}), H_1, H_2 )$ can be used to explicitly construct arithmetically equivalent number fields $K_1$ and $K_2$ as subfields of the $n$-torsion field $\mathbb{Q}(E[n])$ of an elliptic curve $E/\mathbb{Q}$, see \cite{deSmit1998}. Non-conjugate subgroups $H_1$ and $H_2$ give rise to non-isomorphic number fields $K_1$ and $K_2$. \par
	Let $\varphi: \GL_2(\mathbb{Z}/p^2\mathbb{Z}) \rightarrow \GL_2(\mathbb{Z}/p\mathbb{Z})$ be the natural homomorphism and let $H_1$ and $H_2$ be locally conjugate subgroups of $\GL_2(\mathbb{Z}/p^2\mathbb{Z})$. The groups $H_1 \cap \ker \varphi$ and $H_2 \cap \ker \varphi$ must be locally conjugate in $\GL_2(\mathbb{Z}/p^2\mathbb{Z})$ and $\varphi(H_1)$ and $\varphi(H_2)$ must be locally conjugate in $\GL_2(\mathbb{Z}/p\mathbb{Z})$. However, the converse is not generally true. Section \ref{sectionKernel} of this paper identifies local conjugacy in $\GL_2(\mathbb{Z}/p^2\mathbb{Z})$ among the subgroups of $\ker \varphi$. Using the categorization of local conjugacy in $\GL_2(\mathbb{Z}/p\mathbb{Z})$, we replace $H_1$ and $H_2$ with conjugates so that $\varphi(H_1)$ and $\varphi(H_2)$ are similar to each other. Moreover, the elements of $\varphi(H_i)$ restrict $H_i \cap \ker \varphi$ for $i = 1,2$. We further replace $H_1$ and $H_2$ with conjugates so that they are in \enquote{nice} forms. \par
	Theorem \ref{theoremConjClass} completely classifies the conjugacy classes of $\GL_2(\mathbb{Z}/p^2\mathbb{Z})$. Furthermore, the elements of subgroups of $\GL_2(\mathbb{Z}/p^2\mathbb{Z})$ are often not difficult to identify. We use this information to determine when $H_1$ and $H_2$ are not locally conjugate given that $\varphi(H_1)$ and $\varphi(H_2)$ are locally conjugate and $H_1 \cap \ker \varphi$ and $H_2 \cap \ker \varphi$ are locally conjugate or when $H_1$ and $H_2$ are locally conjugate but not conjugate. \par
	Propositions \ref{propCartan} and \ref{propBorel} yield all of the non-conjugate locally conjugate subgroups, which we refer to as nontrivially locally conjugate subgroups, of $\GL_2(\mathbb{Z}/p^2\mathbb{Z})$ up to conjugation, see Theorem \ref{theoremConclusion}. One can choose generators of such subgroups to resemble each other. In fact, they are expressible in forms resembling the nontrivially locally conjugate subgroups of $\GL_2(\mathbb{Z}/p\mathbb{Z})$.

\section{Subgroups of $\GL_2(\mathbb{Z}/p^k\mathbb{Z})$}

Throughout this paper, $p$ is an odd prime and $\epsilon$ is taken to be some nonsquare in $\mathbb{Z}/p \mathbb{Z}$. Let $\GL_2(R)$, $\SL_2(R)$ and $\PGL_2(R)$ denote the general, special and projective linear groups of $2 \times 2$ matrices over a ring $R$. 

Define the following subgroups of $\GL_2(\mathbb{Z}/p^k\mathbb{Z})$ for $k \geq 1$:
\begin{align*}
	Z(p^k) &= \left \{ \begin{pmatrix} w & 0 \\ 0 & w \end{pmatrix} \in \GL_2(\mathbb{Z}/p^k\mathbb{Z}) \right \} \\
	C_s(p^k) &= \left \{ \begin{pmatrix} w & 0 \\ 0 & z \end{pmatrix} \in \GL_2(\mathbb{Z}/p^k\mathbb{Z}) \right \} \\
	C_{ns}(p^k) &= \left \{ \begin{pmatrix} w & \epsilon y \\ y & w \end{pmatrix} \in \GL_2(\mathbb{Z}/p^k\mathbb{Z}) \right \} \\
	B(p^k) &= \left \{ \begin{pmatrix} w & x \\ 0 & z \end{pmatrix} \in \GL_2(\mathbb{Z}/p^k\mathbb{Z}) \right \}.
\end{align*}
They are respectively called the \textbf{center}, \textbf{Cartan-split subgroup}, \textbf{Cartan-nonsplit subgroup}, and \textbf{Borel subgroup} of $\GL_2(\mathbb{Z}/p^k\mathbb{Z})$. For $H \leq \GL_2(\mathbb{Z}/p^k\mathbb{Z})$, let $N(H)$ denote the normalizer of $H$ in $\GL_2(\mathbb{Z}/p^k\mathbb{Z})$. In particular,
\begin{align*}
	N(C_s(p)) &= C_s(p) \cup \left \{ \begin{pmatrix} 0 & x \\ y & 0 \end{pmatrix} \in \GL_2(\mathbb{Z}/p \mathbb{Z}) \right \} \\
	&= C_s(p) \cup \begin{pmatrix} 0 & 1 \\ 1 & 0 \end{pmatrix} C_s(p) \\
	N(C_{ns}(p)) &= C_{ns}(p) \cup \left \{ \begin{pmatrix} w & \epsilon y \\ -y & -w \end{pmatrix} \in \GL_2(\mathbb{Z}/p\mathbb{Z}) \right \} \\ 
	&= C_{ns}(p) \cup \begin{pmatrix} 1 & 0 \\ 0 & -1 \end{pmatrix} C_{ns}(p).
\end{align*}

\section{Properties of Locally Conjugate Subgroups}
This section defines locally conjugate subgroups and discusses some properties of local conjugacy. For a group $G$ and an element $g \in G$, let $g^G$ denote the conjugacy class of $g$ in $G$. 
\begin{definition}
	Let $G$ be a group and $H_1, H_2 \leq G$. If there is a bijection $f: H_1 \rightarrow H_2$ such that $h$ and $f(h)$ are conjugate in $G$ for all $h \in H_1$, then $H_1$ and $H_2$ are \textbf{locally conjugate in $G$}. If $H_1$ and $H_2$ are locally conjugate in $G$ but not conjugate in $G$, then $H_1$ and $H_2$ are \textbf{nontrivially locally conjugate in $G$}.
\end{definition}

Conjugate subgroups of $G$ are always locally conjugate in $G$. Hence, conjugate subgroups are considered to be \enquote{trivially} locally conjugate. Moreover, it is possible for two subgroups $H_1$ and $H_2$ of a group $G'$, which is in turn a subgroup of $G$, to be locally conjugate in $G$ but not in $G'$. It is therefore important to emphasize the parent group in which two subgroups are locally conjugate in. Nevertheless, the parent group will be clear in context even if it is not explicitly stated.\par

Local conjugacy of subgroups in a fixed group $G$ is an equivalence relation. Furthermore, given that subgroups $H_2$ and $H_3$ of $G$ are conjugate, a subgroup $H_1$ of $G$ which is locally conjugate to $H_2$ is conjugate to $H_2$ if and only if $H_1$ is conjugate to $H_3$. With this in mind, we will categorize local conjugacy between subgroups of $\GL_2(\mathbb{Z}/p^2\mathbb{Z})$ up to conjugation of the subgroups. 

The following gives an alternate definition to local conjugacy:
\begin{proposition}\label{propEquivDefLocCon}
	Let $G$ be a group and $H_1, H_2 \leq G$. Then, $H_1$ and $H_2$ are locally conjugate in $G$ if and only if $|H_1 \cap C| = |H_2 \cap C|$ for all conjugacy classes of $G$.
	\begin{proof}
		Suppose that $H_1$ and $H_2$ are locally conjugate in $G$ via $f: H_1 \rightarrow H_2$. For every conjugacy class $C$ of $G$, $f \mid_{H_1 \cap C}$ maps into $H_2 \cap C$. Likewise, $f^{-1} \mid_{H_2 \cap C}$ maps into $H_1 \cap C$ and is the inverse of $f \mid_{H_1 \cap C}$. Thus, $|H_1 \cap C| = |H_2 \cap C|$. \par
		Conversely, suppose that $|H_1 \cap C| = |H_2 \cap C|$ for every conjugacy class $C$ of $G$. Choose some bijections $f_C: H_1 \cap C \rightarrow H_2 \cap C$ and define $f: H_1 \rightarrow H_2$ as $f(h) = f_{h^G} (h)$. Since the conjugacy classes of $G$ partition $G$, $f$ is a well defined bijection. Moreover, $h$ and $f(h)$ are in the same conjugacy class for every $h \in H_1$, and so $H_1$ and $H_2$ are locally conjugate. 
	\end{proof}	
\end{proposition}

The following two propositions give necessary conditions for local conjugacy on a finite group $G$ in terms of local conjugacy in a normal subgroup of $G$ and quotient groups of $G$.
\begin{proposition}\label{propNecessaryOne}
	Let $G$ be a group, $H_1, H_2 \leq G$ and $N \vartriangleleft G$. If $H_1$ and $H_2$ are locally conjugate in $G$, then $H_1 \cap N$ and $H_2 \cap N$ are locally conjugate in $G$. 
		\begin{proof}
			$N$ is the disjoint union of some conjugacy classes of $G$. Let $C$ be a conjugacy class of $G$. If $C \subseteq N$, then $|(H_i \cap N) \cap C| = |H_i \cap C|$ for $i = 1,2$. Otherwise, $|(H_i \cap N) \cap C| = 0$. $H_1 \cap N$ and $H_2 \cap N$ are therefore locally conjugate in $G$ by Proposition \ref{propEquivDefLocCon}.
		\end{proof}
\end{proposition}

\begin{proposition}\label{propNecessaryTwo}
	Let $G, G'$ be finite groups, $H_1, H_2 \leq G$ subgroups of $G$ and $\varphi: G \rightarrow G'$ a surjective homomorphism. If $H_1$ and $H_2$ are locally conjugate in $G$, then $\varphi(H_1)$ and $\varphi(H_2)$ are locally conjugate in $G'$.
	\begin{proof}
		Let $C'$ be any conjugacy class of $G'$ and let $U = \bigcup_{x \in C'} (\varphi^{-1}(x))^G$. We claim that $\varphi^{-1}(C') = U$. If $d \in \varphi^{-1}(C')$, then $\varphi(d) \in C'$, in which case $\varphi(d) \in (\varphi^{-1}(\varphi(d))^G \subseteq U$. Therefore, $\varphi^{-1}(C') \subseteq U$. Conversely, if $d \in (\varphi^{-1}(x))^G$ for some $x \in C'$, then $d = gyg^{-1}$ for some $g \in G$ and $y \in \varphi^{-1}(x)$. It follows that $\varphi(d) = \varphi(g) \varphi(y) \varphi(g)^{-1} = \varphi(g)x \varphi(g)^{-1}$, and so $d \in \varphi^{-1}(C')$. Hence, $U \subseteq \varphi^{-1}(C')$ as desired. In particular, $\varphi^{-1}(C')$ is the union of conjugacy classes of $G$. \par
	$\ker \varphi$ is the union of conjugacy classes of $G$ because it is normal in $G$. Moreover, since $H_1$ and $H_2$ are locally conjugate, $|H_1 \cap \ker \varphi| = |H_2 \cap \ker \varphi|$. Similarly, $|H_1 \cap \varphi^{-1}(C')| = |H_2 \cap \varphi^{-1}(C')|$. Note that $\varphi(H_i) \cap C' = \varphi(H_i \cap \varphi^{-1}(C'))$, and so $H_i \cap \ker \varphi$ has index $|\varphi(H_i \cap C')|$ in $H_i \cap \varphi^{-1}(C')$. Thus, $|\varphi(H_1) \cap C'| = |\varphi(H_2) \cap C'|$ and so $\varphi(H_1)$ and $\varphi(H_2)$ are locally conjugate by Proposition \ref{propEquivDefLocCon}. 
	\end{proof}
\end{proposition}

\section{The Kernel of the Natural Homomorphism $\varphi:\GL_2(\mathbb{Z}/p^2\mathbb{Z}) \rightarrow \GL_2(\mathbb{Z}/p\mathbb{Z})$} \label{sectionKernel}
For the rest of this paper, let $\varphi$ denote the natural homomorphism $\GL_2(\mathbb{Z}/p^2\mathbb{Z}) \rightarrow \GL_2(\mathbb{Z}/p\mathbb{Z})$. An element $\kappa$ of $\ker \varphi$ is of the form $\kappa = I + Ap$, where $A$ is identifiable as an element of $\Mat_2(\mathbb{Z}/p\mathbb{Z})$ and $A$ uniquely determines $\kappa$. We will refer to $A$ as the \textbf{$p$-part of $\kappa$} and define $p(\kappa)$ to be $A$. \par
	Note that $\ker \varphi$ is isomorphic to the $4$ dimensional $\mathbb{Z}/p\mathbb{Z}$ vector space because $(I + A_1p)(I + A_2p) = I + (A_1+A_2)p$ for all $A_1,A_2 \in \Mat_2(\mathbb{Z}/p\mathbb{Z})$. \par
	Suppose that $H_1$ and $H_2$ are subgroups of $\GL_2(\mathbb{Z}/p^2\mathbb{Z})$ that are locally conjugate in $\GL_2(\mathbb{Z}/p^2\mathbb{Z})$. By Proposition \ref{propNecessaryOne}, $H_1 \cap \ker \varphi$ and $H_2 \cap \ker \varphi$ are also locally conjugate in $\GL_2(\mathbb{Z}/p^2\mathbb{Z})$. This section determines the subgroups of $\ker \varphi$ which are locally conjugate in $\GL_2(\mathbb{Z}/p^2\mathbb{Z})$. \par

	Lemma \ref{lemmaKerEltConj} below determines when two elements of $\ker \varphi$ are conjugate in $\GL_2(\mathbb{Z}/p^2\mathbb{Z})$. 
	\begin{lemma} \label{lemmaKerEltConj}
		Let $\kappa_1, \kappa_2 \in \ker \varphi$, $A_1 = p(\kappa_1)$ and $A_2 = p(\kappa_2)$. Then, $\kappa_1$ and $\kappa_2$ are conjugate in $\GL_2(\mathbb{Z}/p^2\mathbb{Z})$ if and only if $A_1$ and $A_2$ are conjugate by an element of $\GL_2(\mathbb{Z}/p\mathbb{Z})$. 
	\begin{proof}
		If $\kappa_1 = I + A_1 p$ is conjugate to $\kappa_2 = I+A_2 p$ via $g \in \GL_2(\mathbb{Z}/p^2\mathbb{Z})$, i.e. $g(I+A_1 p) g^{-1} = I + A_2p$, then $I+gA_1g^{-1} p = I + A_2 p$ and so $A_1$ is conjugate to $A_2$ via $\varphi(g)$. Conversely, if $A_1$ is conjugate to $A_2$ are conjugate via some $g' \in \GL_2(\mathbb{Z}/p\mathbb{Z})$, then $I + A_1 p$ is conjugate to $I+A_2 p$ via any $g \in \varphi^{-1}(g')$. 
	\end{proof}
\end{lemma}
	
	Furthermore, Lemma \ref{lemmaConjDefined} below yields a well defined group action of $\GL_2(\mathbb{Z}/p\mathbb{Z})$ on $\ker \varphi$ in which $g \in \GL_2(\mathbb{Z}/p\mathbb{Z})$ sends $\kappa \in \ker \varphi$ to $\hat{g} \kappa \hat{g}^{-1}$, where $\hat{g}$ is any element of $\varphi^{-1}(g)$.
\begin{lemma} \label{lemmaConjDefined} 
Let $g,g' \in \GL_2(\mathbb{Z}/p^2\mathbb{Z})$ with $\varphi(g) = \varphi(g')$, i.e. $g$ and $g'$ are congruent modulo $p$. For any $\kappa \in \ker \varphi$, $g\kappa g^{-1} = g'\kappa g'^{-1}$.   
\begin{proof}
	Let $A = p(\kappa)$ and let $g' = g+Bp$ for some $B \in \Mat_2(\mathbb{Z}/p\mathbb{Z})$. It is not difficult to see that $g'^{-1} = g^{-1} - g^{-1} B g^{-1} p$. Therefore,
	\begin{align*}
		g' \kappa g'^{-1} &= g'(I+Ap)g'^{-1} \\
										&= (g+Bp)(I+Ap) (g^{-1}-g^{-1}Bg^{-1} p) \\
										&= I + (Bg^{-1} + gAg^{-1} - Bg^{-1}) p \\
										&= I + gAg^{-1} \\
										&= g(I+Ap)g^{-1} \\
										&= g \kappa g^{-1}.
	\end{align*}
\end{proof}
\end{lemma}

Lemma \ref{lemmaUsedMuch} restricts the possible combinations of $\varphi(H)$ and $H \cap \ker \varphi$ for subgroups $H$ of $\GL_2(\mathbb{Z}/p^2\mathbb{Z})$. 

\begin{lemma} \label{lemmaUsedMuch}
	If $H \leq \GL_2(\mathbb{Z}/p^2\mathbb{Z})$, $h \in \varphi(H)$ and $\kappa \in H \cap \ker \varphi$, then $h\kappa h^{-1} \in \ker \varphi$. Furthermore, $H \cap \ker \varphi$ is fixed under conjugation by $h$. 
	\begin{proof}
		This is because $H \cap \ker \varphi$ is a normal subgroup of $H$.
	\end{proof}
\end{lemma}

\subsection{Conjugacy Classes of $\GL_2(\mathbb{Z}/p^2\mathbb{Z})$}

For $k \geq 1$, the similarity classes of $\Mat(\mathbb{Z}/p^k\mathbb{Z})$ are defined as the orbits of $\Mat(\mathbb{Z}/p^k\mathbb{Z})$ under conjugation by elements of $\GL_2(\mathbb{Z}/p^k\mathbb{Z})$. The similarity classes of $\Mat(\mathbb{Z}/p^k\mathbb{Z})$ extend the conjugacy classes of $\GL_2(\mathbb{Z}/p^k\mathbb{Z})$ in that the conjugacy classes are themselves similarity classes. \par
\cite[Theorem 2.2]{0708.1608} yields a way to categorize the similarity classes of $\Mat_2(\mathbb{Z}/p^k\mathbb{Z})$ for $k \geq 1$. Just as in \cite[Section 1.2]{0708.1608}, fix a section $\mathbb{Z}/p\mathbb{Z} \hookrightarrow \mathbb{Z}/p^k\mathbb{Z}$ with image $K_1 \subset \mathbb{Z}/p^k\mathbb{Z}$. Further fix compatible sections $\mathbb{Z}/p^l \mathbb{Z} \rightarrow \mathbb{Z}/p^k\mathbb{Z}$ for $1 \leq l < k$. \cite[Lemma 2.1]{0708.1608} asserts that $\alpha \in \Mat_2(\mathbb{Z}/p^k\mathbb{Z})$ can be written in the form
	\begin{align*}
		\alpha = Id + \beta p^l
	\end{align*}
with $l \in \{0,\ldots,k\}$ maximal such that $\alpha$ is congruent to a scalar matrix modulo $p^l$, with unique $d \in K_l$ and unique nonscalar $\beta \in \Mat_2(\mathbb{Z}/p^{k-l}\mathbb{Z})$. \cite[Theorem 2.2]{0708.1608} concludes the following:
\begin{theorem} \label{theoremConjClass}
	With $\alpha \in \Mat_2(\mathbb{Z}/p^k\mathbb{Z})$ expressed in the form $\alpha = Id + \beta p^l$ as above, $l \in \{ 0,\ldots,k \}$, $d \in K_l$, and $\trace(\beta), \det(\beta) \in \mathbb{Z}/p^{k-l}\mathbb{Z}$ completely determine the conjugacy class of $g$. 
\end{theorem}

\subsection{Orbits of $\Mat_2(\mathbb{Z}/p\mathbb{Z})$ under conjugation by elements of $\GL_2(\mathbb{Z}/p\mathbb{Z})$ }
\cite[Table 3.1]{Sutherland} lists representatives for all the distinct conjugacy classes of $\GL_2(\mathbb{Z}/p\mathbb{Z})$. Table \ref{tableorbit} below uses Theorem \ref{theoremConjClass} to extend \cite[Table 3.1]{Sutherland} to include the representatives of the similarity classes of $\Mat_2(\mathbb{Z}/p\mathbb{Z})$. By Lemma \ref{lemmaKerEltConj}, representatives of the conjugacy classes of elements of $\ker \varphi$ can be given as $I + Ap$, where $A$ is one of the matrices in Table \ref{tableorbit}. 

\begin{table}[h] \label{tableorbit}
\caption{Representatives of the Similarity Classes of $\Mat_2(\mathbb{Z}/p\mathbb{Z})$}
\begin{tabular}{|c|c|c|}
	\hline
	Representative & $\det$ & $\trace$ \\
	\hline
	$\begin{pmatrix} w & 0 \\ 0 & w \end{pmatrix} \label{test}$ $0 \leq w < p$ & $w^2$ & $2w$  \\
	$\begin{pmatrix} w & 1 \\ 0 & w \end{pmatrix}$ $0 \leq w < p$ & $w^2$ & $2w$  \\
	$\begin{pmatrix} w & 0 \\ 0 & z \end{pmatrix}$ $0 \leq w < z < p$ & $wz$ & $w+z$  \\
	$\begin{pmatrix} w & \epsilon y \\ y & w \end{pmatrix}$ $0 < y \leq \frac{p-1}{2}$ & $w^2-\epsilon y^2$ & $2w$  \\
	\hline
\end{tabular}
\end{table}

\subsection{An Equivalent Condition for Local Conjugacy Between Subgroups of $\ker \varphi$}

Table \ref{tableorbit} yields the following observation:
\begin{lemma} \label{lemmaSimClass2}
	The similarity class of a nonscalar element $M$ of $\Mat_2(\mathbb{Z}/p\mathbb{Z})$ is uniquely determined by $\trace(M), \det(M) \in \mathbb{Z}/p\mathbb{Z}$.
\end{lemma}

We define $\chi$ below to use Lemma \ref{lemmaSimClass2} as a way to find an equivalent condition for local conjugacy between subgroups of $\ker \varphi$. 
\begin{definition}
	For a subgroup $H$ of $\ker \varphi$ and for $t,d \in \mathbb{Z}/p\mathbb{Z}$, let $\chi(H,t,d) = |\{ k \in H \mid k = I + Ap \text{ where } \trace(A) = t, \det(A) = d \}|$. \par
	Let $H_1$ and $H_2$ be subgroups of $\ker \varphi$. Say that $H_1$ and $H_2$ have \textbf{equal trace-determinant distribution} if $\chi(H_1,t,d) = \chi(H_2,t,d)$ for all $t,d \in \mathbb{Z}/p\mathbb{Z}$.
\end{definition}

\begin{proposition} \label{propLocConjKerEquiv}
	Let $H_1$ and $H_2$ be subgroups of $\ker \varphi$. Then, $H_1$ and $H_2$ are locally conjugate in $\GL_2(\mathbb{Z}/p^2\mathbb{Z})$ if and only if $H_1 \cap Z(p^2) = H_2 \cap Z(p^2)$ and $H_1$ and $H_2$ have equal trace-determinant distribution. 
	\begin{proof}
		If $H_1$ and $H_2$ are locally conjugate, then $H_1 \cap Z(p^2) = H_2 \cap Z(p^2)$ because every element of $Z(p^2)$ is the sole member of its conjugacy class. Moreover, $H_1$ and $H_2$ must have equal trace-determinant distribution by Lemma \ref{lemmaSimClass2}. \par
		Conversely, if $H_1 \cap Z(p^2) = H_2 \cap Z(p^2)$ and $H_1$ and $H_2$ have equal trace-determinant distribution, then for all $t,d \in \mathbb{Z}/p\mathbb{Z}$, 
		\begin{align*}
		&|\{ k \in H_1 \setminus Z(p^2) \mid k = I + Ap \text{ where } \trace(A) = t, \det(A) = d \}| \\
	 = &|\{ k \in H_2 \setminus Z(p^2) \mid k = I + Ap \text{ where } \trace(A) = t, \det(A) = d \}|.
		\end{align*}
		By Lemmas \ref{lemmaKerEltConj} and \ref{lemmaSimClass2}, $H_1$ and $H_2$ are locally conjugate.  
	\end{proof}
\end{proposition}

\subsection{Preliminary Results for Local Conjugacy in $\GL_2(\mathbb{Z}/p\mathbb{Z})$ among Subgroups of $\ker \varphi$}  \label{subsectionPRK}

The definition of locally conjugate subgroups yields the following result:
\begin{lemma} \label{lemmaLocConjBij}
	If $H_1, H_2 \leq \ker \varphi$ are locally conjugate in $\GL_2(\mathbb{Z}/p\mathbb{Z})$, then $\dim H_1 = \dim H_2$.
	\begin{proof}
		Locally conjugate subgroups are in bijection and $\varphi$ is a finite dimensional vector space over the finite field $\mathbb{Z}/p\mathbb{Z}$. 
	\end{proof}
\end{lemma}

Lemma \ref{lemmaLocConjCyclic} categorizes local conjugacy of finite cyclic subgroups.
\begin{lemma}\label{lemmaLocConjCyclic}
	Let $G$ be a group and let $H_1,H_2$ be finite locally conjugate subgroups of $G$. If $H_1$ is cyclic, then $H_2$ is cyclic and $H_1$ and $H_2$ are conjugate.
	\begin{proof}
		Say that $h_1$ generates $H_1$. There is some $h_2 \in H_2$ which is conjugate to $h_1$ in $G$. The orders of $h_2$ and $h_1$ are equal, $H_1$ and $H_2$ are finite and $|H_1| = |H_2|$, and so $h_2$ generates $H_2$. Therefore, $H_2$ is conjugate to $H_1$. 
	\end{proof}
\end{lemma}


\begin{lemma} \label{lemmaKer01}
	The subgroups of $\ker \varphi$ of dimension $0$ or $1$, i.e. the cyclic subgroups, are conjugate in $\GL_2(\mathbb{Z}/p\mathbb{Z})$ to one of the following subgroups of $\ker \varphi$:
	\begin{enumerate}
		\item $\langle I \rangle$
		\item $\left \langle I + \begin{pmatrix} 0 & 1 \\ 0 & 0 \end{pmatrix} p \right \rangle$
		\item $\left \langle I + \begin{pmatrix} 1 & 1 \\ 0 & 1 \end{pmatrix} p \right \rangle$
		\item $\left \langle I + \begin{pmatrix} 1 & 0 \\ 0 & d \end{pmatrix} p \right \rangle$, where $d \in \mathbb{Z}/p\mathbb{Z}$
		\item $\left \langle I + \begin{pmatrix} 0 & \epsilon \\ 1 & 0 \end{pmatrix} p \right \rangle$
		\item $\left \langle I + \begin{pmatrix} 1 & \epsilon c \\ c & 1 \end{pmatrix} p \right \rangle$, where $c \in \mathbb{Z}/p\mathbb{Z}$ and $0 < c \leq \frac{p-1}{2}$. 
	\end{enumerate}
	No two distinct subgroups among these are locally conjugate.
	\begin{proof}
		Let $H$ be a subgroup of $\ker \varphi$ generated by $h = I + Ap$ for some $A \in \Mat_2(\mathbb{Z}/p\mathbb{Z})$. $H$ can be replaced with a conjugate such that $A$ is one of the matrices in Table \ref{tableorbit}. The categorization of cyclic subgroups of $\ker \varphi$ is finished by determining an alternative generator of $H$ and conjugating $H$ if necessary. \par
		Suppose that $A = \begin{pmatrix} w & 0 \\ 0 & w \end{pmatrix}$. If $w = 0$, then $H = \langle I \rangle$. Otherwise, $H = \left \langle I + \begin{pmatrix} 1 & 0 \\ 0 & 1 \end{pmatrix} p \right \rangle$. \par
		Suppose that $A = \begin{pmatrix} w & 1 \\ 0 & w \end{pmatrix}$. If $w = 0$, then $H = \left \langle I + \begin{pmatrix} 0 & 1 \\ 0 & 0 \end{pmatrix} p \right \rangle$. Otherwise, $H$ is alternatively generated by $I + \begin{pmatrix} 1 & \frac{1}{w} \\ 0 & 1 \end{pmatrix} p$, and so $H$ is conjugate to $\left \langle I + \begin{pmatrix} 1 & 1 \\ 0 & 1 \end{pmatrix} p \right \rangle$. \par
		Suppose that $A = \begin{pmatrix} w & 0 \\ 0 & z \end{pmatrix}$, where $0 \leq w < z < p$. If $w = 0$, then $H$ is alternatively generated by $I + \begin{pmatrix} 0 & 0 \\ 0 & 1 \end{pmatrix} p$, and so $H$ is conjugate to $\left \langle I + \begin{pmatrix} 1 & 0 \\ 0 & 0 \end{pmatrix} p \right \rangle$. Otherwise, $H$ is alternatively generated by $I + \begin{pmatrix} 1 & 0 \\ 0 & \frac{z}{w} \end{pmatrix} p$. \par
		Suppose that $A = \begin{pmatrix} w & \epsilon y \\ y & w \end{pmatrix}$, where $0 < y \leq \frac{p-1}{2}$. If $w = 0$, then $H = \left \langle I + \begin{pmatrix} 0 & \epsilon \\ 1 & 0 \end{pmatrix} p \right \rangle$. Otherwise, $H$ is alternatively generated by $I + \begin{pmatrix} 1 & \epsilon \frac{y}{w} \\ \frac{y}{w} & 1 \end{pmatrix} p$. \par
		It is not difficult to see that for all $h_1 \in H_1$ and $h_2 \in H_2$ where $H_1$ and $H_2$ are distinct two groups among the ones listed in the statement of the lemma, $h_1$ and $h_2$ are not conjugate in $\GL_2(\mathbb{Z}/p^2\mathbb{Z})$ unless $h_1 = h_2 = I$. Thus, no two of the subgroups listed are locally conjugate in $\GL_2(\mathbb{Z}/p^2\mathbb{Z})$. 
	\end{proof}
\end{lemma}

\subsection{Local Conjugacy in $\GL_2(\mathbb{Z}/p\mathbb{Z})$ among Subgroups of $\ker \varphi \cap \SL(p^2)$} \label{subsectionT}

From now on, let $T$ denote $\ker \varphi \cap \SL(p^2)$. Note that $T$ is the subgroup of $\varphi$ with exactly the matrices of the form $I + Ap$, where $A \in \Mat(\mathbb{Z}/p\mathbb{Z})$ has trace $0$. Note that $T$ has dimension $3$ as $I + \begin{pmatrix} 1 & 0 \\ 0 & -1 \end{pmatrix} p$, $I + \begin{pmatrix} 0 & 1 \\ 0 & 0 \end{pmatrix} p$ and $I + \begin{pmatrix} 0 & 0 \\ 1 & 0 \end{pmatrix} p$ form a basis of $T$. Additionally, $T$ is normal in $\GL_2(\mathbb{Z}/p^2\mathbb{Z})$ because both $\varphi$ and $\SL(p^2)$ are normal in $\GL_2(\mathbb{Z}/p^2\mathbb{Z})$. A result analogous to Lemma \ref{lemmaLocConjBij} thus follows:
\begin{lemma} \label{lemmaLocConjT}
	If $H_1, H_2 \leq \GL_2(\mathbb{Z}/p^2\mathbb{Z})$ are locally conjugate in $\GL_2(\mathbb{Z}/p^2\mathbb{Z})$, then $H_1 \cap T$ and $H_2 \cap T$ are locally conjugate in $\GL_2(\mathbb{Z}/p^2\mathbb{Z})$ and $\dim (H_1 \cap T) = \dim(H_2 \cap T)$. 
\end{lemma} 

All dimension $0$ and $1$ subgroups of $\ker \varphi$ are categorized up to conjugation in Lemma \ref{lemmaKer01}. Moreover, the only dimension $3$ subgroup of $T$ is $T$ itself. Lemma \ref{lemmaT2} categorizes the 2 dimensional subgroups of $T$ up to conjugacy as well as local conjugacy among them. It will be useful to consult Lemmas \ref{lemmaDiagConj} and \ref{lemmaSkewConj} in Section \ref{sectionComputation} for several of the upcoming lemmas.

\begin{lemma}\label{lemmaT2}
	The subgroups of $T$ of dimension $2$ are conjugate in $\GL_2(\mathbb{Z}/p^2\mathbb{Z})$ to one of the following:
	\begin{enumerate}
		\item $H_1 = \left \langle I + \begin{pmatrix} 1 & 0 \\ 0 & -1 \end{pmatrix} p, I + \begin{pmatrix} 0 & 1 \\ 0 & 0 \end{pmatrix} p \right \rangle$
		\item $H_2 = \left \langle I + \begin{pmatrix} 1 & 0 \\ 0 & -1 \end{pmatrix} p, I + \begin{pmatrix} 0 & 1 \\ 1 & 0 \end{pmatrix} p \right \rangle$
		\item $H_3 = \left \langle I + \begin{pmatrix} 1 & 0 \\ 0 & -1 \end{pmatrix} p, I + \begin{pmatrix} 0 & \epsilon \\ 1 & 0 \end{pmatrix} p \right \rangle$.
	\end{enumerate}
	No two distinct subgroups among these are locally conjugate. 
\begin{proof}
	Let $H \leq T$ have dimension $2$. Suppose, for contradiction, that $\det(p(h)) \neq -a^2$ for every $h \in H$ and any nonzero $a \in \mathbb{Z}/p\mathbb{Z}$. If there is some nonidentity $h \in H$ such that $\det(p(h)) = 0$, then $h$ is conjugate to $u_1 = I + \begin{pmatrix} 0 & 1 \\ 0 & 0 \end{pmatrix} p$ by Lemmas \ref{lemmaKerEltConj} and \ref{lemmaSimClass2}. Replace $H$ with a conjugate so that $u_1 \in H$. Since $H$ is $2$ dimensional, there is some $u_2 \in H$ of the form $u_2 = I + \begin{pmatrix} a & 0 \\ c & -a \end{pmatrix} p$ for some $a,c \in \mathbb{Z}/p\mathbb{Z}$ where $a$ and $c$ are not both $0$. $a$ must be $0$ because $\det(p(u_2)) = -a^2$. $c$ is therefore nonzero and so $I + \begin{pmatrix} 0 & 0 \\ 1 & 0 \end{pmatrix} p \in H$. By extension, $I + \begin{pmatrix} 0 & 1 \\ 1 & 0 \end{pmatrix} p \in H$, but $\det \left( p \left( I + \begin{pmatrix} 0 & 1 \\ 1 & 0 \end{pmatrix} p \right) \right) = -1$, which is a contradiction.   \par
	Otherwise, $\det(p(h)) \neq -a^2$ for every nonidentity $h \in H$ and any $a \in \mathbb{Z}/p\mathbb{Z}$. For every nonidentity $h \in H$, $- \frac{\det(p(h))}{\epsilon}$ is a nonzero square in $\mathbb{Z}/p\mathbb{Z}$. Thus, $h$ is conjugate to $v_1 = I + \begin{pmatrix} 0 & \epsilon y \\ y & 0 \end{pmatrix} p$ for some nonzero $y \in \mathbb{Z}/p\mathbb{Z}$. Replace $H$ with a conjugate so that $v_1 \in H$. Since $H$ is $2$ dimensional, there is some nonidentity $v_2 \in H$ of the form $v_2 = I + \begin{pmatrix} a & 0 \\ c & -a \end{pmatrix}p$ for some $a,c \in \mathbb{Z}/p\mathbb{Z}$. However, $\det(p(v_2)) = -a^2$, which is a contradiction. \par
	Hence, there is some nonidentity $h \in H$ such that $\det(p(h)) = -a^2$ for some nonzero $a \in \mathbb{Z}/p\mathbb{Z}$. The $p$-part of $h^{\frac{1}{a}}$ has determinant $-1$ and trace $0$, and so $h^{\frac{1}{a}}$ is an element of $H$ which is conjugate to $I + \begin{pmatrix} 1 & 0 \\ 0 & -1 \end{pmatrix} p$. \par 
	Replace $H$ with a conjugate so that $I + \begin{pmatrix} 1 & 0 \\ 0 & -1 \end{pmatrix} p \in H$ and let $w_1 = I + \begin{pmatrix} 1 & 0 \\ 0 & -1 \end{pmatrix} p$. Since $H$ is $2$ dimensional, there is some nonidentity $w_2 \in H$ of the form $w_2 = I + \begin{pmatrix} 0 & b \\ c & 0 \end{pmatrix} p$. If $c = 0$, then $b \neq 0$ and $H = H_1$. If $b = 0$, then $c \neq 0$ and $H$ is conjugate to $H_1$ via $\begin{pmatrix} 0 & 1 \\ 1 & 0 \end{pmatrix}$. Now assume that $b,c \neq 0$. If $bc$ is a square, then $H$ is conjugate to $H_2$ via $\begin{pmatrix} 1 & 0 \\ 0 & \sqrt{\frac{b}{c}} \end{pmatrix}$. Otherwise, $bc$ is not a square, in which case $H$ is conjugate to $H_3$ via $\begin{pmatrix} 1 & 0 \\ 0 & \sqrt{\frac{b}{c\epsilon}} \end{pmatrix}$. \par
	It remains to show that $H_1, H_2$ and $H_3$ are not locally conjugate to one another. The $p$-parts of the elements of $H_1, H_2$ and $H_3$ are respectively of the form
	\begin{align*}
	 &x_1 \begin{pmatrix} 1 & 0 \\ 0 & -1 \end{pmatrix} + y_1 \begin{pmatrix} 0 & 1 \\ 0 & 0 \end{pmatrix}, \\
		&x_2 \begin{pmatrix} 1 & 0 \\ 0 & -1 \end{pmatrix} + y_2 \begin{pmatrix} 0 & 1 \\ 1 & 0 \end{pmatrix}, \text{ and} \\
		&x_3 \begin{pmatrix} 1 & 0 \\ 0 & -1 \end{pmatrix} + y_3 \begin{pmatrix} 0 & \epsilon \\ 1 & 0 \end{pmatrix},
	\end{align*}
	where $x_i,y_i \in \mathbb{Z}/p\mathbb{Z}$. These $p$-parts have trace $0$ and have determinants $-x_1^2, -x_2^2-y_2^2$ and $-x_3^2-\epsilon y_3^2$ respectively. Since $\mathbb{Z}/p\mathbb{Z}$ has nonsquares, there is some $x_2$ for which $x_2^2+1$ is a nonsquare. Setting $x_2$ to be such a value and $y_2 = 1$ makes $-x_2^2-y_2^2 = -(x_2^2+1)$, which shows that $H_1$ and $H_2$ are not locally conjugate by Proposition \ref{propLocConjKerEquiv}. Letting $x_3 = 0$ and $y_3 = 1$ shows that $H_1$ and $H_3$ are not locally conjugate as well. Moreover, $-x_2^2-y_2^2 = 0$ has solutions such that $(x_2,y_2) \neq (0,0)$ exactly when $-1$ is a square in $\mathbb{Z}/p\mathbb{Z}$, which is exactly when $-x_3^2-\epsilon y_3^2$ does not have solutions such that $(x_3,y_3) \neq (0,0)$. $H_2$ and $H_3$ are therefore not locally conjuguate.
\end{proof}
\end{lemma}

\subsection{Final Results for Local Conjugacy in $\GL_2(\mathbb{Z}/p\mathbb{Z})$ among Subgroups of $\ker \varphi$} \label{subsectionFRK}
Since $\dim(T) = 3$ and $\dim(\ker \varphi) = 4$, any subgroup of $\ker \varphi$ with at least $2$ dimensions must have nontrivial intersection with $T$. In particular, letting $H$ be a subgroup of $\ker \varphi$, if $\dim(H) = 2$, then $\dim(H \cap T) \geq 1$ and if $\dim(H) = 3$, then $\dim(H \cap T) \geq 2$. 

\begin{lemma} \label{lemmaKer2}
	The subgroups of $\ker \varphi$ of dimension $2$ that are not subgroups of $T$ are conjugate in $\GL_2(\mathbb{Z}/p^2\mathbb{Z})$ to one of the following:
	\begin{enumerate}
		\item $H_1 = \left \langle I + \begin{pmatrix} 0 & 1 \\ 0 & 0 \end{pmatrix} p, I + \begin{pmatrix} 0 & 0 \\ 1 & 1 \end{pmatrix} p \right \rangle$
		\item $H_2 = \left \langle I + \begin{pmatrix} 0 & 1 \\ 0 & 0 \end{pmatrix} p, I + \begin{pmatrix} 0 & 0 \\ 0 & 1 \end{pmatrix} p \right \rangle$
		\item $H_{3,d} = \left \langle I + \begin{pmatrix} 0 & 1 \\ 0 & 0 \end{pmatrix} p, I + \begin{pmatrix} 1 & 0 \\ 0 & d \end{pmatrix} p \right \rangle$, where $d \in \mathbb{Z}/p\mathbb{Z}$ is not $-1$. 
		\item $H_{4,c} = \left \langle I + \begin{pmatrix} 1 & 0 \\ 0 & -1 \end{pmatrix} p, I + \begin{pmatrix} 0 & 1 \\ c & 1 \end{pmatrix} p \right \rangle$, where $c \in \mathbb{Z}/p\mathbb{Z}$. 
		\item $H_{5} = \left \langle I + \begin{pmatrix} 1 & 0 \\ 0 & -1 \end{pmatrix} p, I + \begin{pmatrix} 0 & 0 \\ 0 & 1 \end{pmatrix} p \right \rangle$
		\item $H_{6,a,b} = \left \langle I + \begin{pmatrix} 0 & \epsilon \\ 1 & 0 \end{pmatrix} p, I + \begin{pmatrix} 1+a & -\epsilon b \\ b & 1-a \end{pmatrix} p \right \rangle$, where $a,b \in \mathbb{Z}/p\mathbb{Z}$. 
	\end{enumerate}
	In particular, $H_2$ and $H_{3,0}$ are nontrivially locally conjugate. For $d,d' \in \mathbb{Z}/p\mathbb{Z}$ where $d,d' \neq -1$, $H_{3,d}$ and $H_{3,d'}$ are nontrivially locally conjugate if $d \neq d'$ and $dd' = 1$. For $a_1,a_2,b_1,b_2 \in \mathbb{Z}/p\mathbb{Z}$, $H_{6,a_1,b_1}$ and $H_{6,a_2,b_2}$ are conjugate if $a_1^2 - \epsilon b_1^2 = a_2^2 - \epsilon b_2^2$. All other pairs of distinct subgroups listed above are not locally conjugate.
	\begin{proof}
		Let $H$ be a $2$ dimensional subgroup of $\ker \varphi$ that is not a subgroup of $T$. $H \cap T$ has dimension $1$. Replace $H$ with a conjugate so that $H \cap T$ is one of the subgroups of $T$ as listed in Lemma \ref{lemmaKer01}, i.e. $H \cap T = \left \langle u \right \rangle$, where $u = I + \begin{pmatrix} 0 & 1 \\ 0 & 0 \end{pmatrix} p$, $I + \begin{pmatrix} 1 & 0 \\ 0 & -1 \end{pmatrix} p$, or $I + \begin{pmatrix} 0 & \epsilon \\ 1 & 0 \end{pmatrix} p$. \par
		Suppose that $u = I + \begin{pmatrix} 0 & 1 \\ 0 & 0 \end{pmatrix} p$. Choose a nonidentity element $h \in H$ to be of the form $h = I + \begin{pmatrix} a & 0 \\ c & d \end{pmatrix} p$, i.e. $H = \langle u,h \rangle$. Since $\dim(H \cap T) =1 $, $a+d \neq 0$. If $a = 0$, then $d \neq 0$. $h$ can be replaced with $h^{\frac{1}{d}}$ so that $H$ is still $\langle u,h \rangle$ and $d = 1$. If $c = 0$ as well, then $H = H_2$. Otherwise, $H$ is conjugate to $H_1$ via $\begin{pmatrix} c & 0 \\ 0 & 1 \end{pmatrix}$. If $a \neq 0$, then $h$ can be replaced with $h^{\frac{1}{a}}$ so that $H = \langle u,h \rangle$, $a = 1$ and $a+d \neq 0$. If $c = 0$, then $H = H_{3,d}$. If $c \neq 0$, then $H$ is conjugate to $H_1$ via $\begin{pmatrix} \frac{c}{d+1} & - \frac{1}{d+1} \\ 0 & 1 \end{pmatrix}$. \par
		Suppose that $u = I + \begin{pmatrix} 1 & 0 \\ 0 & -1 \end{pmatrix} p$. Choose a nonidentity element $h \in H$ to be of the form $I + \begin{pmatrix} 0 & b \\ c & d \end{pmatrix} p$. Since $\trace(p(h)) = d \neq 0$, replacing $h$ with $h^{\frac{1}{d}}$ makes $d = 1$. If $b = c = 0$, then $H = H_5$. If $b \neq 0$, then $H$ is conjugate to $H_{4,\frac{c}{b}}$ via $\begin{pmatrix} 1 & 0 \\ 0 & b \end{pmatrix}$. Otherwise, $b = 0$ and $c \neq 0$, but conjugating $H$ via $\begin{pmatrix} 0 & 1 \\ 1 & 0 \end{pmatrix}$ reduces $H$ to the case where $ b \neq 0$. \par
		Suppose that $u = I + \begin{pmatrix} 0 & \epsilon \\ 1 & 0 \end{pmatrix} p$. Choose a nonidentity element $h \in H$ so that $\trace(p(h)) = 2$, i.e. $h$ is of the form $h = I + \begin{pmatrix} 1+a & b' \\ c' & 1-a \end{pmatrix} p$. Letting $b = \frac{-b'+\epsilon c'}{2\epsilon}$, compute
		\begin{align*}
			h u^{ - \frac{b'+\epsilon c'}{2\epsilon}} &= \left( I + \begin{pmatrix} 1+a & b' \\ c' & 1-a \end{pmatrix} p \right) \left(I + \begin{pmatrix} 0 & \epsilon \\ 1 & 0 \end{pmatrix} p \right)^{- \frac{b'+\epsilon c'}{2 \epsilon}} \\ 
			&= \left( I + \begin{pmatrix} 1+a & b' \\ c' & 1-a \end{pmatrix} p \right) \left(I + \begin{pmatrix} 0 & - \frac{b'+\epsilon c'}{2} \\ - \frac{b'+\epsilon c'}{2 \epsilon}  & 0 \end{pmatrix} p \right) \\
			&= \left( I + \begin{pmatrix} 1+a & \frac{b'-\epsilon c'}{2} \\ \frac{-b'+\epsilon c'}{2\epsilon} & 1-a \end{pmatrix} p \right) \\
			&= I + \begin{pmatrix} 1+a & -b\epsilon \\ b & 1-a \end{pmatrix} p.
		\end{align*}
		Replacing $h$ with $I + \begin{pmatrix} 1+a & -b \epsilon \\ b & 1-a \end{pmatrix} p$ shows that $H = H_{6,a,b}$. \par
		It remains to determine local conjugacy among the listed subgroups. If $H$ and $H'$ are locally conjugate and among the subgroups listed, then $H \cap T$ and $H' \cap T$ must be locally conjugate by Lemma \ref{lemmaLocConjT}. Thus, $H \cap T$ and $H' \cap T$ are equal due to the the way in which they were chosen in the beginning of the proof. In particular, $H_1,H_2,H_{3,d}$ are not locally conjugate to $H_{4,c}, H_5, H_{6,a,b}$ and $H_{4,c},H_5$ are not locally conjugate to $H_{6,a,b}$. For an element $h \in \ker \varphi$, we will respectively call $\det(p(h))$ and $\trace(p(h))$ simply the determinant and trace of $h$ for the rest of the proof. \par 
		The elements of $H_2$ all have zero determinant, and so $H_2$ is not locally conjugate to $H_1$ or $H_{3,d}$ where $d \neq 0$. Elements of $H_2$ are of the form $I + \begin{pmatrix} 0 & x \\ 0 & y \end{pmatrix} p$ and elements of $H_{3,0}$ are of the form $I + \begin{pmatrix} y & x \\ 0 & 0 \end{pmatrix} p$ where $x,y \in \mathbb{Z}/p\mathbb{Z}$. By Lemmas \ref{lemmaKerEltConj} and \ref{lemmaSimClass2}, $I + \begin{pmatrix} 0 & x \\ 0 & y \end{pmatrix} p$ is conjugate to $I + \begin{pmatrix} y & x \\ 0 & 0 \end{pmatrix} p$, and so $H_2$ and $H_{3,0}$ are locally conjugate. Suppose, for contradiction, that $H_2$ and $H_{3,0}$ are conjugate, say via $g \in \GL_2(\mathbb{Z}/p\mathbb{Z})$.
		In this case, $H_2 \cap T = H_{3,0} \cap T = g (H_2 \cap T) g^{-1}$ because $T$ is normal. Using that $H_{3,0} \cap T = H_2 \cap T = \left \langle I + \begin{pmatrix} 0 & 1 \\ 0 & 0 \end{pmatrix} p \right \rangle$, it is not difficult to see that $g$ must be upper triangular. However, $H_2 = gH_2g^{-1}$ in this case. Hence, $H_2$ and $H_{3,0}$ are nontrivially locally conjugate. This fully categorizes local conjugacy of $H_2$ with the other subgroups. \par
		An element of $H_1$ with trace $1$ can have any determinant, whereas a trace $1$ element of $H_{3,d}$ can only have determinant $\frac{d}{(d+1)^2}$. $H_1$ and $H_{3,d}$ are therefore not locally conjugate by Proposition \ref{propLocConjKerEquiv}. This fully categorizes local conjugacy of $H_1$ with the other subgroups. \par
		$H_{3,0}$ is not locally conjugate to $H_{3,d}$ where $d \neq 0$ because $H_{3,0}$ has only elements of determinant $0$ whereas $H_{3,d}$ has elements of nonzero determinant. Moreover, $H_{3,-1}$ is not locally conjugate to $H_{3,d}$ where $d \neq -1$ because the latter has elements of nonzero trace, whereas the former does not. Let $d_1, d_2 \in \mathbb{Z}/p\mathbb{Z}$ such that $d_1, d_2 \neq 0,-1$. For $t \in \mathbb{Z}/p\mathbb{Z}$ and $i=1,2$, the trace $t$ elements of $H_{3,d_i}$ are of the form $I + \begin{pmatrix} \frac{t}{d_i+1} & x \\ 0 & \frac{td_i}{d_i+1} \end{pmatrix} p$ where $x \in \mathbb{Z}/p\mathbb{Z}$. Such an element has determinant $\frac{t^2d_i}{(d_i+1)^2}$. Thus, $H_{3,d_1}$ and $H_{3,d_2}$ are locally conjugate exactly when $H_{3,d_1} \cap Z(p^2) = H_{3,d_2} \cap Z(p^2)$ and $\frac{d_1}{(d_1+1)^2} = \frac{d_2}{(d_2+1)^2}$ by Proposition \ref{propLocConjKerEquiv}. The latter condition is equivalent to  
		\begin{align*}
			0 &= d_1(d_2+1)^2 - d_2(d_1+1)^2 \\
				&= (d_1d_2-1)(d_2-d_1),
		\end{align*}	
		i.e. $d_1 = d_2$ or $d_1d_2 = 1$. One can check that the latter condition implies the former condition. Hence, $H_{3,d_1}$ and $H_{3,d_2}$ are locally conjugate exactly when $d_1 = d_2$ or when $d_1d_2 = 1$. \par
		Suppose, for contradiction, that $H_{3,d_1}$ and $H_{3,d_2}$ are conjugate but $d_1 \neq d_2$. Let $g \in \GL_2(\mathbb{Z}/p\mathbb{Z})$ satisfy $gH_{3,d_1} g^{-1} = H_{3,d_2}$. Using that $H_{3,d_1} \cap T = H_{3,d_2} \cap T = \left \langle I + \begin{pmatrix} 0 & 1 \\ 0 & 0 \end{pmatrix} p \right \rangle$, one can deduce that $g$ must be upper triangular, but then $gH_{3,d_1} g^{-1} = H_{3,d_1}$, which contradicts $H_{3,d_1} \neq H_{3,d_2}$. Hence, $H_{3,d_1}$ and $H_{3,d_2}$ are nontrivially locally conjugate if $d_1 \neq d_2$ and $d_1d_2 = 1$. This fully categorizes local conjugacy of $H_{3,d}$ with the other subgroups. \par
		The trace $1$ elements of $H_{4,c}$ are of the form $I + \begin{pmatrix} x & 1 \\ c & 1-x \end{pmatrix} p$ for $x \in \mathbb{Z}/p\mathbb{Z}$, and such an element has determinant $-x^2+x-c$. Compute $-x^2+x-c = -\left( x - \frac{1}{2} \right)^2 - c + \frac{1}{4}$, and so $-x^2+x-c$ takes the value $-c+\frac{1}{4}$ exactly once and all other values in $\mathbb{Z}/p\mathbb{Z}$ exactly $2$ or $0$ times. Therefore, if $c,c' \in \mathbb{Z}/p\mathbb{Z}$ are distinct, then $H_{4,c}$ and $H_{4,c'}$ are not locally conjugate. \par
		Note that $H_{4,c} \cap Z(p^2) = \langle I \rangle$, whereas $H_5 \cap Z(p^2) \neq \langle I \rangle$. Thus, $H_{4,c}$ and $H_5$ are not locally conjugate. This fully categorizes local conjugacy of $H_{4,c}$ and $H_5$ with the other subgroups. \par
		Suppose $a_1,b_1,a_2,b_2 \in \mathbb{Z}/p\mathbb{Z}$ satisfy $a_1^2-\epsilon b_1^2 = a_2^2 - \epsilon b_2^2$. Conjugating $H_{6,a,b}$ via $\begin{pmatrix} -\sqrt{\epsilon} & -\epsilon \\ -\sqrt{\epsilon} & \epsilon \end{pmatrix} \in \GL_2(\mathbb{F}_{p^2})$ results in the group 
		\[
		\left \langle I + \begin{pmatrix} \sqrt{\epsilon} & 0 \\ 0 & -\sqrt{\epsilon} \end{pmatrix} p, I + \begin{pmatrix} 1 & a+b\sqrt{\epsilon} \\ a-b\sqrt{\epsilon} & 1 \end{pmatrix} p \right \rangle. \footnote{See Lemma \ref{lemmaEpsilonConj}}
		\]
		Further conjugating this group by $\begin{pmatrix} \alpha & 0 \\ 0 & \delta \end{pmatrix} \in \GL_2(\mathbb{F}_{p^2})$ results in 
		\[
		\left \langle I + \begin{pmatrix} \sqrt{\epsilon} & 0 \\ 0 & -\sqrt{\epsilon} \end{pmatrix} p, I + \begin{pmatrix} 1 & (a + b\sqrt{\epsilon}) \frac{\alpha}{\delta} \\ (a-b\sqrt{\epsilon}) \frac{\delta}{\alpha} & 1 \end{pmatrix} p \right \rangle.
		\]
		Therefore, $H_{6,a_1,b_1}$ is conjugate to $H_{6,a_2,b_2}$ via 
		\[
		\begin{pmatrix} -\sqrt{\epsilon} & -\epsilon \\ -\sqrt{\epsilon} & \epsilon \end{pmatrix}^{-1} \begin{pmatrix} a_2+b_2\sqrt{\epsilon} & 0 \\ 0 & a_1+b_1\sqrt{\epsilon} \end{pmatrix} \begin{pmatrix} -\sqrt{\epsilon} & -\epsilon \\ -\sqrt{\epsilon} & \epsilon \end{pmatrix},
		\]
		which is a scalar multiple of
		\[
		\begin{pmatrix} (a_1+a_2)^2-(b_1+b_2)^2 \epsilon & 2(a_2b_1-a_1b_2) \\ \frac{2(a_2b_1-a_1b_2)}{\epsilon} & (a_1+a_2)^2-(b_1+b_2)^2\epsilon \end{pmatrix}.
		\]
		$H_{6,a_1,b_1}$ and $H_{6,a_2,b_2}$ are thus conjugate in $\GL_2(\mathbb{Z}/p^2\mathbb{Z})$. \par 
		The trace $1$ elements of $H_{6,a,b}$ are of the form $I + \begin{pmatrix} \frac{1+a}{2} & -\frac{\epsilon b}{2} + \epsilon x \\ \frac{b}{2} + x & \frac{1-a}{2} \end{pmatrix} p$, whose determinant is $\frac{1-a^2+\epsilon b^2}{4} - \epsilon x^2$. This expression takes the value $\frac{1-a^2+\epsilon b^2}{4}$ exactly once and all values of $\mathbb{Z}/p\mathbb{Z}$ exactly two or zero times. Therefore, if $a_1,b_1,a_2,b_2 \in \mathbb{Z}/p\mathbb{Z}$ satisfy $a_1^2-\epsilon b_1^2 \neq a_2^2 - \epsilon b_2^2$, then $H_{6,a_1,b_1}$ and $H_{6,a_2,b_2}$ are not locally conjugate. This fully categorizes local conjugacy of $H_{6,a,b}$ with the other subgroups.
	\end{proof}
\end{lemma}
 
\begin{lemma} \label{lemmaKer3}
The subgroups of $\ker \varphi$ of dimension $3$ that are not subgroups of $T$ are conjugate in $\GL_2(\mathbb{Z}/p^2\mathbb{Z})$ to one of the following:
\begin{enumerate}
	\item $H_1 = \left \langle I + \begin{pmatrix} 1 & 0 \\ 0 & -1 \end{pmatrix} p, I + \begin{pmatrix} 0 & 1 \\ 0 & 0 \end{pmatrix} p, I + \begin{pmatrix} 0 & 0 \\ 1 & 1 \end{pmatrix} p \right \rangle$ \\
	\item $H_2 = \left \langle I + \begin{pmatrix} 1 & 0 \\ 0 & -1 \end{pmatrix} p, I + \begin{pmatrix} 0 & 1 \\ 0 & 0 \end{pmatrix} p, I + \begin{pmatrix} 0 & 0 \\ 0 & 1 \end{pmatrix} p \right \rangle$ \\
	\item $H_{3,c} = \left \langle I + \begin{pmatrix} 1 & 0 \\ 0 & -1 \end{pmatrix} p, I + \begin{pmatrix} 0 & 1 \\ 1 & 0 \end{pmatrix} p, I + \begin{pmatrix} 0 & 0 \\ c & 1 \end{pmatrix} p \right \rangle$, where $c \in \mathbb{Z}/p\mathbb{Z}$ with $0 \leq c \leq \frac{p-1}{2}$.
	\item $H_{4,c} = \left \langle I + \begin{pmatrix} 1 & 0 \\ 0 & -1 \end{pmatrix} p, I + \begin{pmatrix} 0 & \epsilon \\ 1 & 0 \end{pmatrix} p, I + \begin{pmatrix} 0 & 0 \\ c & 1 \end{pmatrix} p \right \rangle$, where $c \in \mathbb{Z}/p\mathbb{Z}$ with $0 \leq c \leq \frac{p-1}{2}$. 
\end{enumerate}
No two distinct subgroups among these are locally conjugate. 
\begin{proof}
	Let $H$ be a $3$ dimensional subgroup of $\ker \varphi$ that is not a subgroup of $T$. $H \cap T$ has dimension $2$. Replace $H$ with a conjugate so that $H \cap T$ is one of the subgroups of $T$ as listed in Lemma \ref{lemmaT2}. In any of these cases, a third basis element $h$ of $H$ can be chosen to be of the form $h = I + \begin{pmatrix} 0 & 0 \\ c & d \end{pmatrix} p$, where $d \neq 0$. Replace $h$ with $h^{\frac{1}{d}}$ so that $d = 1$. \par
	Suppose that $H \cap T  = \left \langle I + \begin{pmatrix} 1 & 0 \\ 0 & -1 \end{pmatrix} p, I + \begin{pmatrix} 0 & 1 \\ 0 & 0 \end{pmatrix} p \right \rangle$. If $c \neq 0$, then $H$ is conjugate to $H_1$ via $\begin{pmatrix} c & 0 \\ 0 & 1 \end{pmatrix}$. Otherwise, $H = H_2$. \par
	Suppose that $H \cap T = \left \langle I + \begin{pmatrix} 1 & 0 \\ 0 & -1 \end{pmatrix} p, I + \begin{pmatrix} 0 & 1 \\ 1 & 0 \end{pmatrix} p \right \rangle$. If $0 \leq c \leq \frac{p-1}{2}$, then $H$ is $H_{3,c}$. Otherwise, $H$ is conjugate to $H_{3,-c}$ via $\begin{pmatrix} 1 & 0 \\ 0 & -1 \end{pmatrix}$. \par
	Suppose that $H \cap T = \left \langle I + \begin{pmatrix} 1 & 0 \\ 0 & -1 \end{pmatrix} p, I + \begin{pmatrix} 0 & \epsilon \\ 1 & 0 \end{pmatrix} p \right \rangle$. If $0 \leq c \leq \frac{p-1}{2}$, then $H$ is $H_{3,c}$. Otherwise, $H$ is conjugate to $H_{4,-c}$ via $\begin{pmatrix} 1 & 0 \\ 0 & -1 \end{pmatrix}$. \par
	It remains to determine local conjugacy among the listed subgroups. Similarly as in Lemma \ref{lemmaKer2}, if $H$ and $H'$ are among the listed subgroups, then $H \cap T = H' \cap T$. In particular, $H_1$ and $H_2$ are not locally conjugate to $H_{3,c}$ and $H_{4,c}$ and $H_{3,c}$ is not locally conjugate to $H_{4,c}$. \par
	$H_1$ and $H_2$ are not locally conjugate because $H_1 \cap Z(p^2) = \langle I \rangle$, whereas $H_2 \cap Z(p^2) \neq \langle I \rangle$. This fully categorizes local conjugacy of $H_1$ and $H_2$ with the other subgroups.  \par
	$H_{3,0}$ and $H_{3,c}$ where $0 < c \leq \frac{p-1}{2}$ are not locally conjugate because $H_{3,0} \cap Z(p^2) \neq H_{3,c} \cap Z(p^2)$. Let $c_1,c_2 \in \mathbb{Z}/p\mathbb{Z}$ be distinct such that $0 < c_1,c_2 \leq \frac{p-1}{2}$. Elements $h_i \in H_{3,c_i}$ are of the form $h_i = I + \begin{pmatrix} x & y \\ y + c_i z & -x + z \end{pmatrix} p$ where $x,y,z \in \mathbb{Z}/p\mathbb{Z}$. The trace of such an element is $z$. Note that scaling $x,y,z$ by $r \in \mathbb{Z}/p\mathbb{Z}$ multiplies the $\trace(p(h_i))$ by $r$ and $\det(p(h_i))$ by $r^2$. Therefore, local conjugacy of $H_{3,c_i}$ is determined by the determinants of its trace $1$ elements and $H_{3,c_i} \cap Z(p^2)$. Fixing $z = 1$, we compute $\det(p(h_i)) = -x^2+x-y^2-c_i y$. For each $d \in \mathbb{Z}/p\mathbb{Z}$, the number of solutions $(x,y)$ to $d = \det(p(h_i))$ is the number of solutions $(x',y')$ to $x'^2+y'^2 = -d+\frac{1+c_i^2}{4}$, where $x',y' \in \mathbb{Z}/p\mathbb{Z}$ are parametrized as $x' = x - \frac{1}{2}$ and $y' = y + \frac{c_i}{2}$. \par
	For nonzero $K \in \mathbb{Z}/p\mathbb{Z}$, $x'^2$ and $K-y'^2$ each take $\frac{p+1}{2}$ distinct values over $x' \in \mathbb{Z}/p\mathbb{Z}$ and $y' \in \mathbb{Z}/p\mathbb{Z}$ respectively. By the pigeonhole principle, $x'^2+y'^2 = K$ has at least one solution $(x',y')$. In particular, $(\pm x', \pm y')$ yields at least two distinct solutions to $x'^2+y'^2 = K$ even if one of $x'$ and $y'$ is $0$. If $p \equiv 1 \pmod{4}$, then the equation $x'^2+y'^2 = K$ is equivalent to $(x'+jy)(x'-jy) = K$ where $j \in \mathbb{Z}/p\mathbb{Z}$ satisfies $j^2 = -1$. The number of solutions to $x'^2+y'^2 = K$ in this case is therefore $p-1$. \par
	The number of solutions to $\det(p(h_i)) = \frac{1+c_i^2}{4}$ is the number of solutions to $x'^2+y'^2 = 0$, which is $1$ if $p \equiv 3 \pmod{4}$ and $2p-1$ if $p \equiv 1 \pmod{4}$. Since $c_1^2 \neq c_2^2$, the number of solutions to $\det(p(h_2)) = \frac{1+c_1^2}{4}$ is at least $2$ and exactly $p-1$ if $p \equiv 1 \pmod{4}$. $H_{3,c_1}$ and $H_{3,c_2}$ are thus not locally conjugate. Similarly, $H_{4,c_1}$ and $H_{4,c_2}$ are not locally conjugate. 
\end{proof}
\end{lemma}
Local conjugacy in $\ker \varphi$ can be summarized as follows:

\begin{proposition} \label{propNTLCKer}
Let $H_1, H_2 \leq \GL_2(\mathbb{Z}/p^2\mathbb{Z})$ be nontrivially locally conjugate in $\GL_2(\mathbb{Z}/p^2\mathbb{Z})$. $H_1$ and $H_2$ are conjugate, in some order, to the following subgroups for some $d \in \mathbb{Z}/p\mathbb{Z}$ such that $d \neq \pm 1$:
\begin{align*}
	\left \langle I + \begin{pmatrix} 1 & 0 \\ 0 & d \end{pmatrix} p, I + \begin{pmatrix} 0 & 1 \\ 0 & 0 \end{pmatrix}p \right \rangle, \left \langle I + \begin{pmatrix} d & 0 \\ 0 & 1 \end{pmatrix} I + \begin{pmatrix} 0 & 1 \\ 0 & 0 \end{pmatrix} p \right \rangle.
\end{align*}
\begin{proof}
	$\dim H_1 = \dim H_2$ and $\dim(H_1 \cap T) = \dim(H_2 \cap T)$ by Lemmas \ref{lemmaLocConjBij} and \ref{lemmaLocConjT}. The claim follows from Lemmas \ref{lemmaKer01}, \ref{lemmaT2}, \ref{lemmaKer2} and \ref{lemmaKer3}.
\end{proof}
\end{proposition}

\section{Local Conjugacy in $\GL_2(\mathbb{Z}/p\mathbb{Z})$}

Dickson \cite{Dickson} classifies the subgroups of $\GL_2(\mathbb{Z}/p\mathbb{Z})$ based on their images in $\PGL_2(\mathbb{Z}/p\mathbb{Z})$:
\begin{proposition} \label{propSubGpGL2p}
	Let $p$ be an odd prime and let $G$ be a subgroup of $\GL_2(\mathbb{Z}/p\mathbb{Z})$ with image $H$ in $\PGL_2(\mathbb{Z}/p\mathbb{Z})$. If $G$ contains an element of order $p$ then $G \subseteq B(p)$ or $\SL_2(\mathbb{Z}/p\mathbb{Z}) \subseteq G$. Otherwise, one of the following holds:
	\begin{enumerate}
		\item $H$ is cyclic and a conjugate of $G$ lies in $C_s(p)$ or $C_{ns}(p)$.
		\item $H$ is dihedral and a conjugate of $G$ lies in $N(C_s(p))$ or $N(C_{ns}(p))$, but no conjugate of $G$ lies in $C_s(p)$ or $C_{ns}(p)$. 
		\item $H$ is isomorphic to $A_4, S_4$ or $A_5$ and no conjugate of $G$ lies in $N(C_s(p))$ or $N(C_{ns}(p))$. 
	\end{enumerate}
	\begin{proof}
		See \cite[Section 2]{Serre1971/72} or \cite[Lemma 2]{Swinnerton-Dyer1973}. 
	\end{proof}
\end{proposition}

Sutherland \cite{Sutherland} uses this classification to identify local conjugacy among the subgroups of $\GL_2(\mathbb{Z}/p\mathbb{Z})$. 

\begin{theorem} \label{theoremLocConjGL2p}
	Let $H_1, H_2 \leq \GL_2(\mathbb{Z}/p\mathbb{Z})$ be nontrivially locally conjugate in $\GL_2(\mathbb{Z}/p\mathbb{Z})$. $H_1$ and $H_2$ are, in some order, conjugate to the following groups:
	\begin{align*}
		\left \langle D, \begin{pmatrix} 1 & 1 \\ 0 & 1 \end{pmatrix} \right \rangle, \left \langle D', \begin{pmatrix} 1 & 1 \\ 0 & 1 \end{pmatrix} \right \rangle,
	\end{align*}
	where $D \leq C_s(p)$, $D' = \begin{pmatrix} 0 & 1 \\ 1 & 0 \end{pmatrix} D \begin{pmatrix} 0 & 1 \\ 1 & 0 \end{pmatrix}^{-1}$\footnote{Since $D \leq C_s(p^2)$, $\begin{pmatrix} z & 0 \\ 0 & w \end{pmatrix} \in D'$ if and only if $\begin{pmatrix} w & 0 \\ 0 & z \end{pmatrix} \in D$ by Lemma \ref{lemmaSkewConj}.} and $D \neq D'$, i.e. there is some $\begin{pmatrix} w & 0 \\ 0 & z \end{pmatrix} \in D$ such that $\begin{pmatrix} z & 0 \\ 0 & w \end{pmatrix} \not \in D$. 
	\begin{proof}
		See \cite[Lemma 3.6, Corollary 3.30]{Sutherland}.
	\end{proof}
\end{theorem} 

\section{Splitting}
A special case of the Schur-Zassenhaus Theorem, which is stated below in Theorem \ref{theoremSchurZassenhaus}, gives a sufficient condition for certain pairs of subgroups of $\GL_2(\mathbb{Z}/p^2\mathbb{Z})$ to be conjugate. Recall that a Hall subgroup of a finite group is a subgroup whose order is relatively prime to its index.
\begin{theorem}[Schur-Zassenhaus]\label{theoremSchurZassenhaus}
	If $K$ is an abelian normal Hall subgroup of a finite group $G$, then there is a splitting $\psi: G/K \rightarrow G$ which is unique up to conjugation. 
	\begin{proof}
		See \cite[Theorem 7.39, 7.40]{Rotman}
	\end{proof}
\end{theorem}

\begin{proposition}\label{propSchurZassenhaus}
	Suppose that $H_1, H_2 \leq \GL_2(\mathbb{Z}/p^2\mathbb{Z})$. If $H_1 \cap \ker \varphi = H_2 \cap \ker \varphi$, $\varphi(H_1) = \varphi(H_2)$ and $p$ does not divide $\left| \varphi(H_i) \right|$, then $H_1$ and $H_2$ are conjugate in $\GL_2(\mathbb{Z}/p^2\mathbb{Z})$.
	\begin{proof}
		For $i = 1,2$, consider the short exact sequence
		\begin{align*}
			1 \rightarrow \ker \varphi \rightarrow \varphi^{-1} \left( \varphi(H_i) \right) \rightarrow \varphi(H_i) \rightarrow 1.
		\end{align*}
		Since $\ker \varphi$ is a $4$ dimensional $\mathbb{Z}/p\mathbb{Z}$ vector space, $\ker \varphi$ is abelian and $\left| \ker \varphi \right| = p^4$. Moreover, $p$ does not divide $\left| \varphi(H_i) \right|$ by assumption, and so the Schur-Zassenhaus Theorem yields a splitting $\psi: \varphi(H_i) \rightarrow \varphi^{-1} \left( \varphi (H_i) \right)$, which is unique up to conjugation. Similarly, there is a splitting $\psi_i: \varphi(H_i) \rightarrow H_i$ that arise from the short exact sequence
		\begin{align*}
			1 \rightarrow H_i \cap \ker \varphi \rightarrow H_i \rightarrow \varphi(H_i) \rightarrow 1.
		\end{align*}
		Since $H_i$ is a subgroup of $\varphi^{-1} \left( \varphi(H_i) \right)$, $\psi_i$ and $\psi$ are both splittings into $\varphi^{-1} \left( \varphi (H_i) \right)$. $\psi_i$ and $\psi$ are thus conjugate in $ \varphi^{-1} \left( \varphi (H_i) \right)$ and by extension, $\psi_1$ is conjugate to $\psi_2$ via some $g \in \varphi^{-1} \left( \varphi (H_i) \right)$. One can express $H_i$ as the internal semidirect product $H_i = \left( H_i \cap \ker \varphi \right) \rtimes \varphi (H_i)$. Conjugating the expression for $H_1$ yields $gH_1g^{-1} = \left( g \left(H_1 \cap \ker \varphi \right) g^{-1} \right) \rtimes \varphi (H_2)$. By Lemma \ref{lemmaUsedMuch}, $g(H_1 \cap \ker \varphi) g^{-1} = H_1 \cap \ker \varphi$ because $\varphi(g) \in \varphi \left( \varphi^{-1} \left( \varphi \left( H_i \right) \right) \right) = \varphi(H_1)$. Since $\varphi(H_1) = \varphi(H_2)$ by assumption, $gH_1g^{-1} = H_2$ and, so $H_1$ and $H_2$ are conjugate to each other.
	\end{proof}
\end{proposition}

\section{Computational Facts} \label{sectionComputation}

This section lists algebraic computations and facts resulting from such computations that are used in previous sections. \par
	From this point on and unless stated otherwise, $R$ will denote the ring $(\mathbb{Z}_p[\sqrt{\epsilon}])/(p^2\mathbb{Z}_p[\sqrt{\epsilon}])$, which is isomorphic to $(\mathbb{Z}/p^2\mathbb{Z})[\sqrt{\epsilon}]$. The elements of $R$ are identifiable as the sums $a + b\sqrt{\epsilon}$ where $a,b \in \mathbb{Z}/p^2\mathbb{Z}$. Let $\tilde{\varphi}$ denote the natural homomorphism 
	$$\tilde{\varphi}: \GL_2(R) \rightarrow \GL_2(\mathbb{Z}_p[\sqrt{\epsilon}]/p\mathbb{Z}_p[\sqrt{\epsilon}]) \simeq \GL_2(\mathbb{F}_{p^2}).$$
	In particular, $\tilde{\varphi}$ is an extension of $\varphi$. \par
	There is a splitting $(\mathbb{Z}/p\mathbb{Z})^\times \rightarrow (\mathbb{Z}/p^2\mathbb{Z})^\times$ given by $\overline{x} \mapsto x^p$, where $\overline{x} \in (\mathbb{Z}/p\mathbb{Z})^\times$ and $x$ is any lift of $x$ in $\mathbb{Z}/p^2\mathbb{Z}$. This map is well defined because $(x+ap)^p \equiv x^p \pmod{p^2}$ for all $x,a \in \mathbb{Z}/p^2\mathbb{Z}$ by the bionamial theorem. It is a splitting as $x^p \equiv \overline{x} \pmod{p}$.  Similarly, there is a splitting $(\mathbb{Z}/p\mathbb{Z}[\sqrt{\epsilon}])^\times \simeq \mathbb{F}_{p^2}^\times \rightarrow R^\times$ given by $\overline{x} \mapsto x^{p^2}$ where $\overline{x} \in \mathbb{F}_{p^2}^\times$ and $x$ is any lift of $x$ in $R$.  \par
	Let $S$ denote the image of the map $(\mathbb{Z}/p\mathbb{Z}[\sqrt{\epsilon}])^\times \rightarrow R^\times$. For $x \in \mathbb{F}_{p^2}$, we will often abuse notation and let $x$ also denote the lift of $x \in \mathbb{F}_{p^2}$ in $S$. In particular, given $w,x,y,z,a,b,c,d \in \mathbb{Z}/p\mathbb{Z}$, write
	\begin{align*}
		\begin{pmatrix} w & 0 \\ 0 & w \end{pmatrix} + \begin{pmatrix} a & b \\ c & d \end{pmatrix} p \\
		\begin{pmatrix} w & 0 \\ 0 & z \end{pmatrix} + \begin{pmatrix} a & b \\ c & d \end{pmatrix} p \\
		\begin{pmatrix} 0 & x \\ y & 0 \end{pmatrix} + \begin{pmatrix} a & b \\ c & d \end{pmatrix} p \\
	\end{align*}
	to denote some elements of $\GL_2(\mathbb{Z}/p^2\mathbb{Z})$; $w,x,y,z$ as written above are elements of $S \cap \mathbb{Z}/p^2\mathbb{Z}$. Likewise, given $w,x,y,z,a,b,c,d \in \mathbb{F}_{p^2}$, write
	\begin{align*}
		\begin{pmatrix} w & 0 \\ 0 & w \end{pmatrix} + \begin{pmatrix} a & b \\ c & d \end{pmatrix} p \\
		\begin{pmatrix} w & 0 \\ 0 & z \end{pmatrix} + \begin{pmatrix} a & b \\ c & d \end{pmatrix} p \\
		\begin{pmatrix} 0 & x \\ y & 0 \end{pmatrix} + \begin{pmatrix} a & b \\ c & d \end{pmatrix} p \\
	\end{align*}
	to denote some elements of $\GL_2(R)$. 

\begin{lemma} \label{lemmaScalarPower}
	For $w \in (\mathbb{Z}/p\mathbb{Z})^\times$ and $a,b,c,d \in \mathbb{Z}/p\mathbb{Z}$,
	\begin{align*}
			\left( \begin{pmatrix} w & 0 \\ 0 & w \end{pmatrix}  + \begin{pmatrix} a & b \\ c & d \end{pmatrix} p \right)^{p-1} = I - \frac{1}{w}\begin{pmatrix} a & b \\ c & d \end{pmatrix} p
	\end{align*}
	and
	\begin{align*}
		\left( \begin{pmatrix} w & 0 \\ 0 & w \end{pmatrix}  + \begin{pmatrix} a & b \\ c & d \end{pmatrix} p \right)^p = \begin{pmatrix} w & 0 \\ 0 & w \end{pmatrix}.
	\end{align*}
	For $w \in R^\times$ and $a,b,c,d \in R$,
	\begin{align*}
			\left( \begin{pmatrix} w & 0 \\ 0 & w \end{pmatrix}  + \begin{pmatrix} a & b \\ c & d \end{pmatrix} p \right)^{p^2-1} = I - \frac{1}{w}\begin{pmatrix} a & b \\ c & d \end{pmatrix} p
	\end{align*}
	and
	\begin{align*}
			\left( \begin{pmatrix} w & 0 \\ 0 & w \end{pmatrix}  + \begin{pmatrix} a & b \\ c & d \end{pmatrix} p \right)^{p^2} = \begin{pmatrix} w & 0 \\ 0 & w \end{pmatrix}.
	\end{align*}
	\begin{proof}
		Since $\begin{pmatrix} w & 0 \\ 0 & w \end{pmatrix}$ multiplicatively commutes with $\begin{pmatrix} a & b \\ c & d \end{pmatrix}$, the Binomial Theorem applies. 
	\end{proof}
\end{lemma}

\begin{corollary} \label{corScalarMember}
	Let $H$ be a subgroup of $\GL_2(\mathbb{Z}/p^2\mathbb{Z})$. Then, $H$ contains $h = \begin{pmatrix} w & 0 \\ 0 & w \end{pmatrix} + \begin{pmatrix} a & b \\ c & d \end{pmatrix} p$, where $w \in (\mathbb{Z}/p\mathbb{Z})^\times$ and $a,b,c,d \in \mathbb{Z}/p\mathbb{Z}$, if and only if $H$ contains $h_1 = \begin{pmatrix} w & 0 \\ 0 & w \end{pmatrix}$ and $h_2 = I + \begin{pmatrix} a & b \\ c & d \end{pmatrix} p$. 
	\begin{proof}
		Suppose that $h \in H$. $h_1 \in H$ by Lemma \ref{lemmaScalarPower} and $hh_1^{-1} = I + \frac{1}{w} \begin{pmatrix} a & b \\ c & d \end{pmatrix} p$. Therefore, $(hh_1^{-1})^{w} = I + \begin{pmatrix} a & b \\ c & d \end{pmatrix} p = h_2$, and so $h_2 \in H$. \par
		Conversely, if $h_1,h_2 \in H$, then $H$ contains $h_1h_2^{\frac{1}{w}} = h_1 \left( I + \frac{1}{w} \begin{pmatrix} a & b \\ c & d \end{pmatrix} p \right) = h$.  
	\end{proof}
\end{corollary}

\begin{lemma} \label{lemmaDiagPower}
	For $w,z \in (\mathbb{Z} / p\mathbb{Z})^\times$ and $a,b,c,d \in \mathbb{Z}/p\mathbb{Z}$ such that $w \neq z$, 
	\begin{align*}
		\left( \begin{pmatrix} w & 0 \\ 0 & z \end{pmatrix} + \begin{pmatrix} a & b \\ c & d \end{pmatrix} p \right)^{p-1} = I + \begin{pmatrix} -\frac{a}{w} & 0 \\ 0 & -\frac{d}{z} \end{pmatrix} p
	\end{align*}
	and
	\begin{align*}
		\left( \begin{pmatrix} w & 0 \\ 0 & z \end{pmatrix} + \begin{pmatrix} a & b \\ c & d \end{pmatrix} p \right)^{p} = \begin{pmatrix} w & 0 \\ 0 & z \end{pmatrix} + \begin{pmatrix} 0 & b \\ c & 0 \end{pmatrix} p
	\end{align*}
	For $w,z \in R^\times$ and $a,b,c,d \in \mathbb{Z}/p\mathbb{Z}$ such that $w \neq z$, 
	\begin{align*}
		\left( \begin{pmatrix} w & 0 \\ 0 & z \end{pmatrix} + \begin{pmatrix} a & b \\ c & d \end{pmatrix} p \right)^{p^2-1} = I + \begin{pmatrix} -\frac{a}{w} & 0 \\ 0 & -\frac{d}{z} \end{pmatrix} p
	\end{align*}
	and
	\begin{align*}
		\left( \begin{pmatrix} w & 0 \\ 0 & z \end{pmatrix} + \begin{pmatrix} a & b \\ c & d \end{pmatrix} p \right)^{p^2} = \begin{pmatrix} w & 0 \\ 0 & z \end{pmatrix} + \begin{pmatrix} 0 & b \\ c & 0 \end{pmatrix} p
	\end{align*}
	\begin{proof}
		Suppose that $w,z \in (\mathbb{Z}/p\mathbb{Z})^\times$ and $a,b,c,d \in \mathbb{Z}/p\mathbb{Z}$. By expanding,
		\begin{align*}
			\left( \begin{pmatrix} w & 0 \\ 0 & z \end{pmatrix} + \begin{pmatrix} a & b \\ c & d \end{pmatrix} p \right)^{p-1} &= \begin{pmatrix} w & 0 \\ 0 & z \end{pmatrix}^{p-1} + \left( \sum_{k=0}^{p-2} \begin{pmatrix} w & 0 \\ 0 & z \end{pmatrix}^k \begin{pmatrix} a & b \\ c & d \end{pmatrix} \begin{pmatrix} w & 0 \\ 0 & z \end{pmatrix}^{p-2-k} \right) p \\
			&= I + \sum_{k=0}^{p-2} \begin{pmatrix} aw^{p-2} & bw^k z^{p-2-k}  \\ cw^{p-2-k} z^k & dz^{p-2} \end{pmatrix} p.
		\end{align*}
		Since $w \not\equiv z \pmod{p}$, $\sum_{k=0}^{p-2} w^kz^{p-2-k}$ and $\sum_{k=0}^{p-2} w^{p-2-k}z^k$ are both geometric series evaluating to $0$. Therefore, 
		\begin{align*}
			\left( \begin{pmatrix} w & 0 \\ 0 & z \end{pmatrix} + \begin{pmatrix} a & b \\ c & d \end{pmatrix} p \right)^{p-1} &= I + \begin{pmatrix} -aw^{p-2} & 0 \\ 0 & -dz^{p-2} \end{pmatrix} p \\
							&= I + \begin{pmatrix} - \frac{a}{w} & 0 \\ 0 & -\frac{d}{z} \end{pmatrix} p.
		\end{align*}
		From here, 
		\begin{align*}
			\left( \begin{pmatrix} w & 0 \\ 0 & z \end{pmatrix} + \begin{pmatrix} a & b \\ c & d \end{pmatrix} p \right)^{p} = \begin{pmatrix} w & 0 \\ 0 & z \end{pmatrix} + \begin{pmatrix} 0 & b \\ c & 0 \end{pmatrix} p
		\end{align*}
		can be immediately calculated. The claims made for $w,z \in R^\times$ and $a,b,c,d, \in R$ such that $w \neq z$ can be proved similarly.
	\end{proof}
\end{lemma}

\begin{corollary} \label{corDiagMember}
	Let $H$ be a subgroup of $\GL_2(R)$. $H$ contains $h = \begin{pmatrix} w & 0 \\ 0 & z \end{pmatrix} + \begin{pmatrix} a & b \\ c & d \end{pmatrix} p$, where $w,z \in (\mathbb{Z}/p\mathbb{Z})^\times$ and $a,b,c,d \in \mathbb{Z}/p\mathbb{Z}$ with $w \neq z$, if and only if $H$ contains $h_1 = I + \begin{pmatrix} -\frac{a}{w} & 0 \\ 0 & -\frac{d}{z} \end{pmatrix} p$ and $\begin{pmatrix} w & 0 \\ 0 & z \end{pmatrix} + \begin{pmatrix} 0 & b \\ c & 0 \end{pmatrix} p$
	\begin{proof}
		This is immediate from Lemma \ref{lemmaDiagPower}.
	\end{proof}
\end{corollary}

\begin{lemma} \label{lemmaEpsilonConj}
	For $w,y,a,b,c,d \in \mathbb{Z}/p\mathbb{Z}$, 
	\begin{enumerate}
		\item $\begin{pmatrix} -\sqrt{\epsilon} & -\epsilon \\ -\sqrt{\epsilon} & \epsilon \end{pmatrix} \begin{pmatrix} a & b \\ c & d \end{pmatrix} \begin{pmatrix} -\sqrt{\epsilon} & -\epsilon \\ -\sqrt{\epsilon} & \epsilon \end{pmatrix}^{-1} = \frac{1}{2} \begin{pmatrix} (a+d) + \left( \frac{b}{\epsilon} + c \right) \sqrt{\epsilon} & (a-d) + \left( - \frac{b}{\epsilon} + c \right) \sqrt{\epsilon} \\ (a-d) + \left( \frac{b}{\epsilon} - c \right) \sqrt{\epsilon} & (a+d) + \left( -\frac{b}{\epsilon} - c \right) \sqrt{\epsilon} \end{pmatrix}$
		\item $\begin{pmatrix} -\sqrt{\epsilon} & -\epsilon \\ -\sqrt{\epsilon} & \epsilon \end{pmatrix} \begin{pmatrix} w & \epsilon y \\ y & w \end{pmatrix} \begin{pmatrix} -\sqrt{\epsilon} & -\epsilon \\ -\sqrt{\epsilon} & \epsilon \end{pmatrix}^{-1} = \begin{pmatrix} w + \sqrt{\epsilon} y & 0 \\ 0 & w - \sqrt{\epsilon} y \end{pmatrix}$
		\item $\begin{pmatrix} -\sqrt{\epsilon} & -\epsilon \\ -\sqrt{\epsilon} & \epsilon \end{pmatrix} \begin{pmatrix} w & \epsilon y \\ - y & -w \end{pmatrix} \begin{pmatrix} -\sqrt{\epsilon} & -\epsilon \\ -\sqrt{\epsilon} & \epsilon \end{pmatrix}^{-1} = \begin{pmatrix} 0 & w- \sqrt{\epsilon} y \\ w + \sqrt{\epsilon} y & 0 \end{pmatrix}$
	\end{enumerate}
\end{lemma}

\begin{lemma} \label{lemmaDiagConj}
	Let $R$ be a ring. If $\begin{pmatrix} w & 0 \\ 0 & z \end{pmatrix} \in \GL_2(R)$ and $\begin{pmatrix} a & b \\ c & d \end{pmatrix} \in \Mat_2(R)$, then
	\begin{align*}
		\begin{pmatrix} w & 0 \\ 0 & z \end{pmatrix} \begin{pmatrix} a & b \\ c & d \end{pmatrix} \begin{pmatrix} w & 0 \\ 0 & z \end{pmatrix}^{-1} = \begin{pmatrix} a & b \frac{w}{z} \\ c \frac{z}{w} & d \end{pmatrix}. 
	\end{align*}
\end{lemma}

\begin{lemma} \label{lemmaSkewConj}
	Let $R$ be a ring. If $\begin{pmatrix} 0 & x \\ y & 0 \end{pmatrix} \in \GL_2(R)$ and $\begin{pmatrix} a & b \\ c & d \end{pmatrix} \in \Mat_2(R)$, then 
	\begin{align*}
		\begin{pmatrix} 0 & x \\ y & 0 \end{pmatrix} \begin{pmatrix} a & b \\ c & d \end{pmatrix} \begin{pmatrix} 0 & x \\ y & 0 \end{pmatrix}^{-1} = \begin{pmatrix} d & c \frac{x}{y} \\ b \frac{y}{x} & a \end{pmatrix}. 
	\end{align*}
	In particular,
	\begin{align*}
		\begin{pmatrix} 0 & 1 \\ 1 & 0 \end{pmatrix} \begin{pmatrix} a & b \\ c & d \end{pmatrix} \begin{pmatrix} 0 & 1 \\ 1 & 0 \end{pmatrix}^{-1} = \begin{pmatrix} d & c \\ b &  a \end{pmatrix}.
	\end{align*}
\end{lemma}

\begin{lemma} \label{lemmaTPower}
	For $a,b,c,d \in \mathbb{Z}/p\mathbb{Z}$ and $n \in \mathbb{Z}$, 
	\begin{align*}
		\left( \begin{pmatrix} 1 & 1 \\ 0 & 1 \end{pmatrix} + \begin{pmatrix} a & b \\ c & d \end{pmatrix} p \right)^n = \begin{pmatrix} 1 & n \\ 0 & 1 \end{pmatrix} + \begin{pmatrix} an + \frac{c(n-1)n}{2} & \frac{(a+d+c(n-1))(n-1)n}{2} - c \sum_{k=0}^{n-1} k^2 + bn \\ cn & dn + \frac{c(n-1)n}{2} \end{pmatrix} p
	\end{align*}
	\begin{proof}
		Expand
		\begin{align*}
			\left( \begin{pmatrix} 1 & 1 \\ 0 & 1 \end{pmatrix} + \begin{pmatrix} a & b \\ c & d \end{pmatrix} p \right)^n &= \begin{pmatrix} 1 & 1 \\ 0 & 1 \end{pmatrix}^n + \sum_{k=0}^{n-1} \begin{pmatrix} 1 & 1 \\ 0 & 1 \end{pmatrix}^k \begin{pmatrix} a & b \\ c & d \end{pmatrix} \begin{pmatrix} 1 & 1 \\ 0 & 1 \end{pmatrix}^{n-1-k} p \\
			&= \begin{pmatrix} 1 & n \\ 0 & 1 \end{pmatrix} + \sum_{k=0}^{n-1} \begin{pmatrix} 1 & k \\ 0 & 1 \end{pmatrix} \begin{pmatrix} a & b \\ c & d \end{pmatrix} \begin{pmatrix} 1 & n-1-k \\ 0 & 1 \end{pmatrix} p \\
			&= \begin{pmatrix} 1 & n \\ 0 & 1 \end{pmatrix} + \sum_{k=0}^{n-1} \begin{pmatrix} a+ck & (a+ck)(n-1-k) + b+dk \\ c & c(n-1-k)+d \end{pmatrix} p \\
			&= \begin{pmatrix} 1 & n \\ 0 & 1 \end{pmatrix} + \sum_{k=0}^{n-1} \begin{pmatrix}  a+ck & a(n-1) -ak + ck(n-1) - ck^2 + b + dk \\ c & c(n-1)-ck + d \end{pmatrix} p.
		\end{align*}
		Since $p$ is an odd prime, $\sum_{k=0}^{n-1} k = \frac{(n-1)n}{2}$, and so
		\begin{align*}
			\left( \begin{pmatrix} 1 & 1 \\ 0 & 1 \end{pmatrix} + \begin{pmatrix} a & b \\ c & d \end{pmatrix} p \right)^n &= \begin{pmatrix} 1 & n \\ 0 & 1 \end{pmatrix} + \begin{pmatrix} an + \frac{c(n-1)n}{2} & \frac{(a+d+c(n-1))(n-1)n}{2} - c \sum_{k=0}^{n-1} k^2 + bn \\ cn & dn + \frac{c(n-1)n}{2} \end{pmatrix} p
		\end{align*}
	\end{proof}
\end{lemma}

\begin{lemma} \label{lemmaTConj}
	Let $R$ be a ring. For $a,b,c,d \in R$, 
	\begin{align*}
		\begin{pmatrix} 1 & 1 \\ 0  & 1 \end{pmatrix} \begin{pmatrix} a & b \\ c & d \end{pmatrix} \begin{pmatrix} 1 & 1 \\ 0 & 1 \end{pmatrix}^{-1} = \begin{pmatrix} a+c & -a+b-c+d \\ c & -c+d \end{pmatrix}.
	\end{align*}
\end{lemma}

\begin{lemma} \label{lemmaTMember}
	Let $H \leq \GL_2(\mathbb{Z}/p^2\mathbb{Z})$ with $t = \begin{pmatrix} 1 & 1 \\ 0 & 1 \end{pmatrix} \in \varphi(H)$. Suppose that $k = I + \begin{pmatrix} a & b \\ c & d \end{pmatrix} p \in H$. 
	\begin{enumerate}
		\item If $a \neq d$, then $I + \begin{pmatrix} 0 & 1 \\ 0 & 0 \end{pmatrix} p \in H$. \label{lemmaTMemberP1}
		\item If $c \neq 0$, then $I + \begin{pmatrix} 0 & 1 \\ 0 & 0 \end{pmatrix} p, I + \begin{pmatrix} 1 & 0 \\ 0 & -1 \end{pmatrix} p \in H$. 
	\end{enumerate}
	\begin{proof}
		\begin{enumerate}
				\item By Lemma \ref{lemmaUsedMuch}, $tkt^{-1} \in H$. Let $k' = tkt^{-1}k^{-1}$, which must be in $H$ as well. Using Lemma \ref{lemmaTConj}, compute
			\begin{align*}
				k' = tkt^{-1}k^{-1} &= \left( I + \begin{pmatrix} a+c & -a+b-c+d \\ c & -c+d \end{pmatrix} p \right) \left( I - \begin{pmatrix} a & b \\ c & d \end{pmatrix} p \right) \\
											&= I + \begin{pmatrix} c & -a-c+d \\ 0 & -c \end{pmatrix}p.
			\end{align*}
			If $c = 0$, then a power of $k'$ is $I + \begin{pmatrix} 0 & 1 \\ 0 & 0 \end{pmatrix} p$ because $a \neq d$ by assumption. If $c \neq 0$, then $I + \begin{pmatrix} 0 & 1 \\ 0 & 0 \end{pmatrix} p \in H$ because 
			\begin{align*}
				tk't^{-1}k'^{-1} = I + \begin{pmatrix} 0 & -2c \\ 0 & 0 \end{pmatrix} p
			\end{align*}
			is an element of $H$. 
				\item By \ref{lemmaTMemberP1}, $I + \begin{pmatrix} 0 & 1 \\ 0 & 0 \end{pmatrix} p \in H$. Furthermore, since $tkt^{-1}k^{-1} = I + \begin{pmatrix} c & -a-c+d \\ 0 & -c \end{pmatrix} p \in H$, $I + \begin{pmatrix} c & 0 \\ 0 & -c \end{pmatrix} p \in H$ as well. Thus, $I + \begin{pmatrix} 1 & 0 \\ 0 & -1 \end{pmatrix} p \in H$. 
		\end{enumerate}
	\end{proof}
\end{lemma}

\begin{lemma} \label{lemmaTMemberPG3}
	Let $H \leq \GL_2(\mathbb{Z}/p^2\mathbb{Z})$ with $t = \begin{pmatrix} 1 & 1 \\ 0 & 1 \end{pmatrix} \in \varphi(H)$. If $p > 3$, then $I + \begin{pmatrix} 0 & 1 \\ 0 & 0 \end{pmatrix} p \in H$. 
	\begin{proof}
		There is some element $h \in H$ of the form $\begin{pmatrix} 1 & 1 \\ 0 & 1 \end{pmatrix} + \begin{pmatrix} a & b \\ c & d \end{pmatrix} p$. By Lemma \ref{lemmaTPower}, the $p$th power of $h$ is
		\begin{align*}
			h^p &= \begin{pmatrix} 1 & p \\ 0 & 1 \end{pmatrix} + \begin{pmatrix} ap + \frac{c(p-1)p}{2} & \frac{(a+d+c(p-1))(p-1)p}{2} - c \sum_{k=0}^{p-1} k^2 + bp \\ cp & dp + \frac{c(p-1)p}{2} \end{pmatrix} p \\
					&= \begin{pmatrix} 1 & p \\ 0 & 1 \end{pmatrix} + \begin{pmatrix} 0 & -c\sum_{k=0}^{p-1} k^2 \\ 0 & 0 \end{pmatrix} p.
		\end{align*}
		Since $p > 3$, $\sum_{k=0}^{p-1} k^2 = \frac{(p-1)p(2p-1)}{6}$, which is $0$ modulo $p$. Therefore, $I + \begin{pmatrix} 0 & 1 \\ 0 & 0 \end{pmatrix} p \in H$, 
	\end{proof}
\end{lemma}

\begin{lemma} \label{lemmaTProdOne}
	For $a,b,c,d,\beta \in \mathbb{Z}/p\mathbb{Z}$,
	\begin{align*}
	 \left( \begin{pmatrix} 1 & 1 \\ 0 & 1 \end{pmatrix} + \begin{pmatrix} a & b \\ c & d \end{pmatrix} p \right) \left( I + \begin{pmatrix} 0 & \beta \\ 0 & 0 \end{pmatrix} p \right) = \begin{pmatrix} 1 & 1 \\ 0 & 1 \end{pmatrix} + \begin{pmatrix} a & b + \beta \\ c & d \end{pmatrix} p.
	\end{align*}
\end{lemma}

\begin{lemma} \label{lemmaTProdTwo}
	For $a,b,c,d, \alpha, \delta \in \mathbb{Z}/p\mathbb{Z}$, 
	\begin{align*}
		\left( \begin{pmatrix} 1 & 1 \\ 0 & 1 \end{pmatrix} + \begin{pmatrix} a & b \\ c & d \end{pmatrix} p \right) \left( I + \begin{pmatrix} \alpha & 0 \\ 0 & \delta \end{pmatrix} p \right) = \begin{pmatrix} 1 & 1 \\ 0 & 1 \end{pmatrix} + \begin{pmatrix} a + \alpha & b + \delta \\ c & d + \delta \end{pmatrix} p.
	\end{align*}
\end{lemma}

\begin{lemma} \label{lemmaTProdThree}
	For $a,b,c,d, \alpha, \gamma, \delta \in \mathbb{Z}/p\mathbb{Z}$,
	\begin{align*}
		\left( \begin{pmatrix} 1 & 1 \\ 0 & 1 \end{pmatrix} + \begin{pmatrix} a & b \\ c & d \end{pmatrix} p \right) \left( I + \begin{pmatrix} \alpha & 0 \\ \gamma & \delta \end{pmatrix} p \right) = \begin{pmatrix} 1 & 1 \\ 0 & 1 \end{pmatrix} + \begin{pmatrix} a + \alpha + \gamma & b + \delta \\ c + \gamma & d + \delta \end{pmatrix} p.
	\end{align*}
\end{lemma}

\begin{lemma} \label{lemmaConvolutedBp}
	Let $H \leq \GL_2(\mathbb{Z}/p^2\mathbb{Z})$. If $\tau = \begin{pmatrix} 1 & 1 \\ 0 & 1 \end{pmatrix} + \begin{pmatrix} a & b \\ c & d \end{pmatrix} \in H$ and $h = \begin{pmatrix} w & 0 \\ 0 & z \end{pmatrix} \in H$ where $w \not\equiv z \pmod{p}$, then $a = d$ or there is some element of $H$ of the form
	\begin{align*}
		I + \begin{pmatrix} \alpha & 0 \\ 0 & \delta \end{pmatrix} p 
	\end{align*}
	where $\alpha, \delta \in \mathbb{Z}/p\mathbb{Z}$ are unequal.
	\begin{proof}
		Using Lemmas \ref{lemmaDiagConj} and \ref{lemmaTPower}, compute
		\begin{align*}
			h\tau h^{-1} \tau^{-\frac{w}{z}} &= \left( \begin{pmatrix} 1 & \frac{w}{z} \\ 0 & 1 \end{pmatrix} + \begin{pmatrix} a & b \frac{w}{z} \\ c \frac{z}{w} & d \end{pmatrix} p \right) \left( \begin{pmatrix} 1 & -\frac{w}{z} \\ 0 & 1 \end{pmatrix} + \begin{pmatrix} -a\frac{w}{z} + \frac{c \left( \frac{w}{z} +1 \right) \frac{w}{z} }{2} & * \\ - c \frac{w}{z} & -d \frac{w}{z} + \frac{c \left( \frac{w}{z} + 1 \right) \frac{w}{z} }{2} \end{pmatrix} p \right) \\
			&= I + \left( \begin{pmatrix} * & * \\ c \frac{z}{w} & * \end{pmatrix} + \begin{pmatrix} * & * \\ -c \frac{w}{z} & * \end{pmatrix} \right) p \\
			&= I + \begin{pmatrix} * & * \\ c \left( \frac{z}{w} - \frac{w}{z} \right) & * \end{pmatrix} p.
		\end{align*}
		If $c \neq 0$ and $\frac{w}{z} \neq \pm 1$, then $c \left( \frac{z}{w} - \frac{w}{z} \right) \neq 0$ and so $I + \begin{pmatrix} 1 & 0 \\ 0 & -1 \end{pmatrix} p \in H$ by Lemma \ref{lemmaTMember}.
		If $c = 0$, then compute
		\begin{align*}
			h\tau h^{-1} \tau^{-\frac{w}{z}} &= \left( \begin{pmatrix} 1 & \frac{w}{z} \\ 0 & 1 \end{pmatrix} + \begin{pmatrix} a & b \frac{w}{z} \\ 0 & d \end{pmatrix} p \right) \left( \begin{pmatrix} 1 & -\frac{w}{z} \\ 0 & 1 \end{pmatrix} + \begin{pmatrix} -a\frac{w}{z} & * \\ 0 & -d \frac{w}{z} \end{pmatrix} p \right) \\
			&= I + \left( \begin{pmatrix} a & -a \frac{w}{z} + b \frac{w}{z} \\ 0 & d \end{pmatrix} + \begin{pmatrix} -a \frac{w}{z} & * \\ 0 & -d \frac{w}{z} \end{pmatrix} \right) p \\
			&= I + \begin{pmatrix} a \left( 1 - \frac{w}{z} \right) & * \\ 0 & d \left( 1 - \frac{w}{z} \right) \end{pmatrix} p.
			\end{align*}
			If $a \neq d$, then $I + \begin{pmatrix} 0 & 1 \\ 0 & 0 \end{pmatrix} p \in H$ by Lemma \ref{lemmaTMember} and since $\frac{w}{z} \neq 1$ by assumption, $I + \begin{pmatrix} a \left( 1 - \frac{w}{z} \right) & 0 \\ 0 & d \left( 1 - \frac{w}{z} \right) \end{pmatrix} p \in H$ with $a \left( 1 - \frac{w}{z} \right) \neq d \left( 1 - \frac{w}{z} \right)$. \par
			If $\frac{w}{z} = -1$, then compute
		\begin{align*}
			h\tau h^{-1} \tau^{-\frac{w}{z}} &= \left( \begin{pmatrix} 1 & -1 \\ 0 & 1 \end{pmatrix} + \begin{pmatrix} a & -b \\ -c & d \end{pmatrix} p \right) \left( \begin{pmatrix} 1 & 1 \\ 0 & 1 \end{pmatrix} + \begin{pmatrix} a & * \\ c & d \end{pmatrix} p \right) \\
			&= I + \begin{pmatrix} 2a-c & * \\ 0 & 2d-c \end{pmatrix} p.
		\end{align*}
		If $a \neq d$, then $I + \begin{pmatrix} 2a-c & 0 \\ 0 & 2d-c \end{pmatrix} p \in H$ with $2a-c \neq 2d-c$. 
	\end{proof}
\end{lemma}

\section{The Center case} \label{sectionCenter}

This section categorizes local conjugacy for subgroups $H$ of $\GL_2(\mathbb{Z}/p^2\mathbb{Z})$ whose images under $\varphi$ lie in $Z(p)$. To do so, we first understand the structure of such $H$. Let $\sigma: Z(p) \rightarrow Z(p^2)$ be the splitting given by $\begin{pmatrix} \overline{w} & 0 \\ 0 & \overline{w} \end{pmatrix} \mapsto \begin{pmatrix} w & 0 \\ 0 & w \end{pmatrix}$ where $\overline{w} \in (\mathbb{Z}/p\mathbb{Z})
^\times$ and $w \in S$ is the lift of $\overline{w}$. For all $h \in \varphi(H)$, $\sigma(h)$ is an element of $H$ by Lemma \ref{lemmaScalarPower}. It is not difficult to see that $H$ is the direct product $\sigma(\varphi(H)) \times (H \cap \ker \varphi)$ by Corollary \ref{corScalarMember}. 

\begin{lemma}\label{lemmaCenter}
	Let $H_1, H_2 \leq \GL_2(\mathbb{Z}/p^2\mathbb{Z})$ such that $\varphi(H_i) \leq Z(p)$ for $i = 1,2$. Then, $H_1$ and $H_2$ are nontrivially locally conjugate in $\GL_2(\mathbb{Z}/p^2\mathbb{Z})$ if and only if $\varphi(H_1) = \varphi(H_2)$ and $H_1 \cap \ker \varphi$ and $H_2 \cap \ker \varphi$ are nontrivially locally conjugate in $\GL_2(\mathbb{Z}/p^2\mathbb{Z})$. 
	\begin{proof}
		If $H_1$ and $H_2$ are nontrivially locally conjugate, then $\varphi(H_1)$ and $\varphi(H_2)$ are locally conjugate in $\GL_2(\mathbb{Z}/p\mathbb{Z})$ by Proposition \ref{propNecessaryTwo}. Since $\varphi(H_i) \leq Z(p)$, $\varphi(H_1)$ must equal $\varphi(H_2)$. Furthermore, $H_1 \cap \ker \varphi$ and $H_2 \cap \ker \varphi$ are locally conjugate in $\GL_2(\mathbb{Z}/p^2\mathbb{Z})$ by Proposition \ref{propNecessaryOne}. Suppose, for contradiction, that $H_1 \cap \ker \varphi$ and $H_2 \cap \ker \varphi$ are conjugate in $\GL_2(\mathbb{Z}/p^2\mathbb{Z})$. Replace $H_1$ with a conjugate so that $H_1 \cap \ker \varphi = H_2 \cap \ker \varphi$. This conjugation preserves $\varphi(H_1)$. Since $\varphi(H_i) \leq Z(p)$ and $|Z(p)| = p-1$, $H_1$ and $H_2$ are conjugate by Proposition \ref{propSchurZassenhaus}. This contradicts that $H_1$ and $H_2$ are nontrivially locally conjugate. Hence, $H_1 \cap \ker \varphi$ and $H_2 \cap \ker \varphi$ are nontrivially locally conjugate. \par
		Conversely, suppose that $\varphi(H_1) = \varphi(H_2)$ and $H_1 \cap \ker \varphi$ and $H_2 \cap \ker \varphi$ are nontrivially locally conjugate in $\GL_2(\mathbb{Z}/p^2\mathbb{Z})$. Since $H_i = \sigma(\varphi(H_1) \times (H_i \cap \ker \varphi)$, and since $\sigma(\varphi(H_1))$ consists only of scalar matrices, it is not difficult to see that $H_1$ and $H_2$ are locally conjugate. They are not conjugate because any conjugation from $H_1$ to $H_2$ yields a conjugation from $H_1 \cap \ker \varphi$ to $H_2 \cap \ker \varphi$. Hence, $H_1$ and $H_2$ are nontrivially locally conjugate. 
	\end{proof}
\end{lemma}

\section{The Cartan cases}\label{sectionCartan} 
	This section categorizes the subgroups, up to conjugation, of $\GL_2(\mathbb{Z}/p^2\mathbb{Z})$ whose images under $\varphi$ are subgroups of $C_s(p)$ or $C_{ns}(p)$ but not subgroups of $Z(p)$. \par
	The notation below will make the Cartan split and Cartan nonsplit cases similar to each other.
	\begin{definition}
	Let $H \leq \GL_2(R)$. Say that $H$ is of \textbf{type $C_s$} if $H \leq \GL_2(\mathbb{Z}/p^2\mathbb{Z})$ and $\varphi(H) \leq C_s(p)$. Say that $H$ is of \textbf{type $N(C_{s})$} if $H \leq \GL_2(\mathbb{Z}/p^2\mathbb{Z})$ and $\varphi(H) \leq N(C_{s}(p))$. \par
	Let $H'$ be the conjugate of $H$ via $\begin{pmatrix} -\sqrt{\epsilon} & -\epsilon \\ -\sqrt{\epsilon} & \epsilon \end{pmatrix}^{-1}$. Say that $H$ is of \textbf{type $C_{ns}$} if $H' \leq \GL_2(\mathbb{Z}/p^2\mathbb{Z})$ and $\tilde{\varphi}(H') \leq C_{ns}(p)$. Say that $H$ is of \textbf{type $N(C_{ns})$} if $H' \leq \GL_2(\mathbb{Z}/p^2\mathbb{Z})$ and $\tilde{\varphi}(H') \leq N(C_{ns}(p))$. \par
	Say that $H$ is in \textbf{diagonalized form} if $H$ is of type $N(C_s)$ or of type $N(C_{ns})$. 
	\end{definition}
	
	Suppose $H \leq \GL_2(R)$ is of type $N(C_{ns})$. The elements of $\tilde{\varphi}(H)$ must be of the form $\begin{pmatrix} w + \sqrt{\epsilon} y & 0 \\ 0 & w - \sqrt{\epsilon} y \end{pmatrix}$ or of the form $\begin{pmatrix} 0 & w - \sqrt{\epsilon} y \\ w + \sqrt{\epsilon} y & 0 \end{pmatrix}$ by Lemma \ref{lemmaEpsilonConj}. From now on, let $K$ denote the conjugate of $\ker \varphi$ via $\begin{pmatrix} -\sqrt{\epsilon} & -\epsilon \\ -\sqrt{\epsilon} & \epsilon \end{pmatrix}$. Note that $K$ is a $4$ dimensional $\mathbb{Z}/p\mathbb{Z}$ vector space just as $\ker \varphi$ is. By Lemma \ref{lemmaEpsilonConj}, 
	\begin{align*}
		K = \left \langle I + \begin{pmatrix} 1 & 0 \\ 0 & 1 \end{pmatrix} p, I + \begin{pmatrix} 0 & 1 \\ 1 & 0 \end{pmatrix} p, I + \begin{pmatrix} \sqrt{\epsilon} & 0 \\ 0 & -\sqrt{\epsilon} \end{pmatrix} p, I + \begin{pmatrix} 0 & \sqrt{\epsilon} \\ -\sqrt{\epsilon} & 0 \end{pmatrix} p \right \rangle.
	\end{align*}
	Note that $H \cap \ker \tilde{\varphi}$ is a subgroup of $K$ when $H$ is of type $N(C_{ns})$. On the other hand, $H \cap \ker \tilde{\varphi}$ is a subgroup of $\ker \varphi$ when $H$ is of type $N(C_{s})$.  \par
	
	Lemma \ref{lemmaEltConjRPSquared} and Corollaries \ref{corOne} and \ref{corTwo} below show that local conjugacy in $\GL_2(\mathbb{Z}/p^2\mathbb{Z})$ between two subgroups $H_1, H_2$ of $\GL_2(\mathbb{Z}/p^2\mathbb{Z})$ such that $\varphi(H_i) \leq N(C_{ns}(p))$ is equivalently determined by local conjugacy between $H_1'$ and $H_2'$ in $\GL_2(R)$, where $H_i'$ is the conjugate of $H_i$ via $\begin{pmatrix} -\sqrt{\epsilon} & -\epsilon \\ -\sqrt{\epsilon} & \epsilon \end{pmatrix}$.

\begin{lemma}\label{lemmaEltConjRPSquared}
	If $g_1,g_2 \in \GL_2(\mathbb{Z}/p^2\mathbb{Z})$ are conjugate via $g \in \GL_2(R)$, then $g_1$ and $g_2$ are conjugate via some $g' \in \GL_2(\mathbb{Z}/p^2\mathbb{Z})$, whose value is only dependent on $g$. 
	\begin{proof}
		Express $g$ in the form
		\begin{align*}
			g = \begin{pmatrix} \alpha_1 + \alpha_2 \sqrt{\epsilon} & \beta_1 + \beta_2 \sqrt{\epsilon} \\ \gamma_1 + \gamma_2 \sqrt{\epsilon} & \delta_1 + \delta_2 \sqrt{\epsilon} \end{pmatrix},
		\end{align*}
		where $\alpha_i,\beta_i,\gamma_i,\delta_i \in \mathbb{Z}/p^2\mathbb{Z}$ for $i = 1,2$. Since the entries of $g_1$ and $g_2$ are in $\mathbb{Z}/p^2\mathbb{Z}$, 
		\begin{align*}
			\begin{pmatrix} \alpha_i & \beta_i \\ \gamma_i & \delta_i \end{pmatrix} h_1 = h_2 \begin{pmatrix} \alpha_i & \beta_i \\ \gamma_i & \delta_i \end{pmatrix}. 
		\end{align*}
		Therefore, if $\begin{pmatrix} \alpha_i & \beta_i \\ \gamma_i & \delta_i \end{pmatrix}$ is invertible for $i = 1$ or $2$, then $h_2 = g' h_1 g'^{-1}$, where $g' = \begin{pmatrix} \alpha_i & \beta_i \\ \gamma_i & \delta_i \end{pmatrix}$. Otherwise, $\alpha_i \delta_i - \beta_i \gamma_i = 0$ for $i=1,2$, in which case
		\begin{align*}
			\det(g) &= (\alpha_1 \delta_1 - \beta_1 \gamma_1) + (\alpha_2 \delta_2 - \beta_2 \gamma_2) \epsilon + (\alpha_1 \delta_2 + \alpha_2 \delta_1 - \beta_1 \gamma_2 - \beta_2 \gamma_1) \sqrt{\epsilon} \\
			&= (\alpha_1 \delta_2 + \alpha_2 \delta_1 - \beta_1 \gamma_2 - \beta_2 \gamma_1) \sqrt{\epsilon}.
		\end{align*}
		Since $g$ is invertible, $\alpha_1 \delta_2 + \alpha_2 \delta_1 - \beta_1 \gamma_2 - \beta_2 \gamma_1 \neq 0$. Let $g' = \begin{pmatrix} \alpha_1 + \alpha_2 & \beta_1 + \beta_2 \\ \gamma_1 + \gamma_2 & \delta_1 + \delta_2 \end{pmatrix}$. Note that $g'$ is invertible because
		\begin{align*}
			\det(g') &= (\alpha_1 \delta_1 + \alpha_1 \delta_2 + \alpha_2 \delta_1 + \alpha_2 \delta_2) - (\beta_1 \gamma_1 + \beta_1 \gamma_2 + \beta_2 \gamma_1 + \beta_2 \gamma_2) \\
			&= (\alpha_1 \delta_1 - \beta_1 \gamma_1) + (\alpha_2 \delta_2 - \beta_2 \gamma_2) + (\alpha_1 \delta_2 + \alpha_2 \delta_1 - \beta_1 \gamma_2 - \beta_2 \gamma_1) \\
			&= (\alpha_1 \delta_2 + \alpha_2 \delta_1 - \beta_1 \gamma_2 - \beta_2 \gamma_1) \\
			&\neq 0.
		\end{align*}
		Furthermore, $g'g_1 = g_2 g'$, and so $g_2 = g'g_1g'^{-1}$ as desired. 
	\end{proof}
\end{lemma}

The idea behind Lemma \ref{lemmaEltConjRPSquared} can be immediately extended to the following corollaries:
\begin{corollary} \label{corOne}
	Two subgroups of $\GL_2(\mathbb{Z}/p^2\mathbb{Z})$ that are conjugate in $\GL_2(R)$ are conjugate in $\GL_2(\mathbb{Z}/p^2\mathbb{Z})$.
\end{corollary}
\begin{corollary} \label{corTwo}
	Two subgroups of $\GL_2(\mathbb{Z}/p^2\mathbb{Z})$ that are locally conjugate in $\GL_2(R)$ are locally conjugate in $\GL_2(\mathbb{Z}/p^2\mathbb{Z})$.
\end{corollary}

\begin{lemma} \label{lemmaCartanKer}
	Let $H \leq \GL_2(R)$ be in diagonialized form. Suppose that there is some $h \in H$ such that $\tilde{\varphi}(h) = \begin{pmatrix} w & 0 \\ 0 & z \end{pmatrix}$ where $w \neq z$. Then, $H$ has some element $k = I + \begin{pmatrix} a & b \\ c & d \end{pmatrix} p$ if and only if $H$ has both $l = I + \begin{pmatrix} a & 0 \\ 0 & d \end{pmatrix} p$ and $m = I + \begin{pmatrix} 0 & b \\ c & 0 \end{pmatrix} p$.
	\begin{proof}
		Since $k = lm$, all three of $l,m,k$ are in $H$ if two of them are. In particular, if $l,m, \in H$, then $k \in H$. \par
		Conversely, assume that $k \in H$. Use Lemma \ref{lemmaDiagConj} to compute
		\begin{align*}
			hkh^{-1}k^{-1} = I + \begin{pmatrix} 0 & b \left( \frac{w}{z} - 1 \right) \\ c \left( \frac{z}{w} - 1 \right) & 0 \end{pmatrix} p.
		\end{align*}
		If $w = -z$, then $hkh^{-1}k^{-1} = I + \begin{pmatrix} 0 & -2b \\ -2c & 0 \end{pmatrix} p$. Moreover, $m$ is a power of $hkh^{-1}k^{-1}$, and so $m \in H$. Thus, $l \in H$ as well. \par
		Now assume that $w \neq \pm z$. If $b = 0$, then $hkh^{-1}k^{-1} = I + \begin{pmatrix} 0 & 0 \\ c \left( \frac{z}{w} - 1 \right) & 0 \end{pmatrix} p$. If $H$ is of type $N(C_{ns})$ as well, then $c = 0$ by Lemma \ref{lemmaEpsilonConj}, and there is nothing to prove. Otherwise, $H$ is of type $N(C_s)$, in which case $\frac{z}{w} - 1$ is in $\mathbb{Z}/p^2\mathbb{Z}$ and nonzero modulo $p$, and so a power of $hkh^{-1}k^{-1}$ is $m$ and we are done. The case where $c = 0$ is similar. \par
		Assume that $b,c \neq 0$. Use Lemma \ref{lemmaDiagConj} to compute
		\begin{align*}
			h(hkh^{-1}k^{-1})h^{-1} = I + \begin{pmatrix} 0 & b \left( \frac{w}{z}-1 \right) \frac{w}{z} \\ c \left( \frac{z}{w} - 1 \right) \frac{w}{z} & 0 \end{pmatrix} p
		\end{align*}
		to see that $hkh^{-1}k^{-1}$ and $h(hkh^{-1}k^{-1})h^{-1}$ are linearly independent as $\mathbb{Z}/p\mathbb{Z}$-vectors. Moreover, both have $p$-parts whose diagonal entries are $0$. Note that the subspaces of $K$ and $\ker \varphi$ consisting of the matrices whose $p$-parts have $0$ as their diagonal entries are both $2$ dimensional. Therefore, $m \in H$ as desired. 
	\end{proof}
\end{lemma}

For $H \leq \GL_2(R)$ in diagonalized form, define $\Delta_H$ as the set of elements of $H \cap \ker \tilde{\varphi}$ whose $p$-parts are diagonal. In other words, 
\begin{align*}
	\Delta_H =\begin{cases} \left \{ I + \begin{pmatrix} a & 0 \\ 0 & d \end{pmatrix} p \in H \cap \ker \varphi \right \} & \text{if } H \text{ is of type } N(C_s) \\
													\left \{ I + \begin{pmatrix} a & 0 \\ 0 & d \end{pmatrix} p \in H \cap K \right \} & \text{ if } H \text{ is of type } N(C_{ns}).
						\end{cases}
\end{align*}
Similarly, define
\begin{align*}
	\Delta_H^{\perp} = \begin{cases} \left \{ I + \begin{pmatrix} 0 & b \\ c & 0 \end{pmatrix} p \in H \cap \ker \varphi \right \} & \text{if } H \text{ is of type } N(C_s) \\
													\left \{ I + \begin{pmatrix} 0 & b \\ c & 0 \end{pmatrix} p \in H \cap K \right \} & \text{ if } H \text{ is of type } N(C_{ns}).
						\end{cases}
\end{align*}

Lemma \ref{lemmaCartanKer} then yields the following: 
\begin{corollary} \label{corDiagGpKer}
	Let $H \leq \GL_2(R)$ be in diagonalized form. Suppose that there is some $h \in H$ such that $\tilde{\varphi}(h) = \begin{pmatrix} w & 0 \\ 0 & z \end{pmatrix}$ where $w \neq z$. 
	\begin{enumerate}
		\item The groups $\Delta_H$ and $\Delta_H^\perp$ together generate $H \cap \ker \tilde{\varphi}$.
		\item Suppose that $w \neq \pm z$. If $H$ is of type $N(C_s)$, then $\Delta_H^\perp$ is one of the following:
		\begin{enumerate}
			\item $\left \langle I + \begin{pmatrix} 0 & 1 \\ 0 & 0 \end{pmatrix} p, I + \begin{pmatrix} 0 & 0 \\ 1 & 0 \end{pmatrix} p \right \rangle$
			\item $\left \langle I + \begin{pmatrix} 0 & 1 \\ 0 & 0 \end{pmatrix} p \right \rangle$
			\item $\left \langle I + \begin{pmatrix} 0 & 0 \\ 1 & 0 \end{pmatrix} p \right \rangle$
			\item $\langle I \rangle$. 
		\end{enumerate}
		If $H$ is of type $N(C_{ns})$, then $\Delta_H^\perp$ is one of the following:
		\begin{enumerate}
			\item $\left \langle I + \begin{pmatrix} 0 & 1 \\ 1 & 0 \end{pmatrix} p, I + \begin{pmatrix} 0 & -\sqrt{\epsilon}  \\ \sqrt{\epsilon} & 0 \end{pmatrix} p \right \rangle$
			\item $\langle I \rangle$.
		\end{enumerate}
	\end{enumerate}
	\begin{proof}
	\begin{enumerate}
		\item This is immediate from Lemma \ref{lemmaCartanKer}.
		\item Suppose that $H$ is of type $N(C_{s})$. If there is some $k = I + \begin{pmatrix} 0 & b \\ c & 0 \end{pmatrix}p \in \Delta_H^\perp$ where $b,c \neq 0$, then $hkh^{-1} = I + \begin{pmatrix} 0 & b \frac{w}{z} \\ c \frac{z}{w} & 0 \end{pmatrix} p$ is linearly independent to $k$ because $bc \frac{w}{z} \neq bc \frac{z}{w}$. Thus, $\Delta_H^\perp = \left \langle I + \begin{pmatrix} 0 & 1 \\ 0 & 0 \end{pmatrix} p, I + \begin{pmatrix} 0 & 0 \\ 1 & 0 \end{pmatrix} p \right \rangle$. If no such $k$ is in $\Delta_H^\perp$, then $\Delta_H^\perp$ is one of $\left \langle I + \begin{pmatrix} 0 & 1 \\ 0 & 0 \end{pmatrix} p \right \rangle$, $\left \langle I + \begin{pmatrix} 0 & 0 \\ 1 & 0 \end{pmatrix} p \right \rangle$ and $\langle I \rangle$.  \par
		Suppose that $H$ is of type $N(C_{ns})$. If there is some nonidentity $k = I + \begin{pmatrix} 0 & b \\ c & 0 \end{pmatrix} p$, then one of $b$ or $c$ is nonzero. Since $b$ and $c$ are of the form $\alpha + \sqrt{\epsilon} \delta$ and $\alpha - \sqrt{\epsilon} \delta$ for some $\alpha,\delta \in \mathbb{Z}/p\mathbb{Z}$ by Lemma \ref{lemmaEpsilonConj}, $b$ and $c$ are both nonzero. Similarly as in the last paragraph, $hkh^{-1}$ and $k$ are linearly independent, and so $\Delta_H^\perp$ must be $2$-dimensional and hence equal to $\left \langle I + \begin{pmatrix} 0 & 1 \\ 1 & 0 \end{pmatrix} p, I + \begin{pmatrix} 0 & -\sqrt{\epsilon}  \\ \sqrt{\epsilon} & 0 \end{pmatrix} p \right \rangle$.
	\end{enumerate}
	\end{proof}
\end{corollary}

\subsection{The Cartan Cases}

This section categorizes the subgroups $H$ of $\GL_2(R)$ of type $C_s$ or of type $C_{ns}$ such that $\tilde{\varphi}(H) \not\leq Z(p)$ up to conjugation to understand local conjugacy among such subgroups. For $H \leq \GL_2(R)$ let $D_{H}$ denote the group consisting of all elements of $H$ that are diagonal matrices. 
\begin{proposition} \label{propCartanConj}
	Let $H \leq \GL_2(R)$ be of type $C_s$ or of type $C_{ns}$ such that $\tilde{\varphi}(H) \not\leq Z(p)$. There is a conjugate $H'$ of $H$ satisfying the following:
	\begin{enumerate}
		\item $H'$ is of type $C_s$ if $H$ is of type $C_s$ and $H'$ is of type $C_{ns}$ if $H$ is of type $C_{ns}$
		\item $\tilde{\varphi}(H') = \tilde{\varphi}(H)$
		\item $H'$ is the internal semidirect product $\Delta^\perp_{H'} \rtimes D_{H'}$
	\end{enumerate}
	\begin{proof}
		Suppose that $H$ is of type $C_s$, i.e. $\varphi(H)$ is a subgroup of $C_s(p)$. The image of $C_s(p)$ in $\PGL_2(p)$ is cyclic and so $C_s(p)$ is generated by two elements, one of which is of the form $\begin{pmatrix} w & 0 \\ 0 & z \end{pmatrix}$ where $w \neq z$ and the other of which is of the form $\begin{pmatrix} w_0 & 0 \\ 0 & w_0 \end{pmatrix}$. By Corollaries \ref{corScalarMember} and \ref{corDiagMember}, $H$ has elements of the form $h_0 = \begin{pmatrix} w_0 & 0 \\ 0 & w_0 \end{pmatrix}$ and $h = \begin{pmatrix} w & 0 \\ 0 & z \end{pmatrix} + \begin{pmatrix} 0 & b \\ c & 0 \end{pmatrix} p$. Note that $H$ is generated by $h,h_0$ and $H \cap \ker \varphi$. Let $h' = \begin{pmatrix} w & 0 \\ 0 & z \end{pmatrix}$ and let $H' = \langle h', h_0, H \cap \ker \varphi \rangle$. \par
		We show that $H' \cap \ker \varphi = H \cap \ker \varphi$. Clearly, $H' \cap \ker \varphi \supseteq H \cap \ker \varphi$, and so it suffices to show that $H' \cap \ker \varphi \subseteq H \cap \ker \varphi$. The elements of $H'$ are of the form $m = \prod_{i=1}^n m_i$ where $m_i = h'$, $m_i = h_0$, or $m_i \in H \cap \ker \varphi$ for each $i$. By Lemmas \ref{lemmaConjDefined} and \ref{lemmaUsedMuch}, if $m_i \in H \cap \ker \varphi$, then there is some $m_i' \in H \cap \ker \varphi$ such that $m_i h' = h' m'_i$. Moreover, $h_0$ is a scalar matrix and hence commutes with $h'$ and all elements of $H \cap \ker \varphi$, and so $m$ is alternatively of the form $m = h'^{n_1} h_0^{n_2} \prod_{i=1}^{n_3} k_i$, where $k_i \in H \cap \ker \varphi$ for each $i$. Due to the way that $w,z,w_0 \in (\mathbb{Z}/p^2\mathbb{Z})^\times$ are chosen based on $w,z,w_0 \in (\mathbb{Z}/p\mathbb{Z})^\times$, $h'^{n_1} h_0^{n_2} = I$ whenever $\varphi(h')^{n_1} \varphi(h_0)^{n_2} = I$.\footnote{Recall $S$ as defined in Section \ref{sectionComputation}}  Therefore, if $m \in H' \cap \ker \varphi$, then $m \in H \cap \ker \varphi$, and so $H' \cap \ker \varphi = H \cap \ker \varphi$ as desired. \par
		Since $\varphi(H) = \varphi(H')$ and $|C_s(p)|$ is indivisible by $p$, $H$ and $H'$ are conjugate by Proposition \ref{propSchurZassenhaus}. Note that $H' \cap \ker \varphi$ is generated by $\Delta_{H'}$ and $\Delta^\perp_{H'}$ by Corollary \ref{corDiagGpKer}. Thus, $H'$ is generated by $h', h_0, \Delta_{H'}$ and $\Delta^\perp_{H'}$. It is not difficult to see that $\Delta^\perp_{H'}$ is a normal subgroup of $H'$, that $H' = (\Delta^\perp_{H'}) (\langle h', h_0, \Delta_{H'} \rangle)$ and that $\Delta^\perp_{H'} \cap \langle h',h_0, \Delta_{H'} \rangle = \langle I \rangle$. Therefore, $H' = \Delta^\perp_{H'} \rtimes \langle h',h_0, \Delta_{H'} \rangle$. Furthermore, $D_{H'} = \langle h', h_0, \Delta_{H'} \rangle$, and so $H' = \Delta^\perp_{H'} \rtimes D_{H'}$. \par
		The case where $H$ is of type $C_{ns}$ is similar. 
	\end{proof}
\end{proposition}

Before proceeding, we introduce a definition which will be useful for understanding the structure of nontrivially locally conjugate subgroups of $\GL_2(\mathbb{Z}/p^2\mathbb{Z})$.
\begin{definition}
	Let $R$ be a ring and let $D_1$ and $D_2$ be two subgroups of $\GL_2(R)$ consisting only of diagonal matrices. Say that $D_1$ and $D_2$ are \textbf{diagonal swaps} if $\begin{pmatrix} w & 0 \\ 0 & z \end{pmatrix} \in D_1$ exactly when $\begin{pmatrix} z & 0 \\ 0 & w \end{pmatrix} \in D_2$. Equivalently, $D_1$ and $D_2$ are conjugate to each other via $\begin{pmatrix} 0 & 1 \\ 1 & 0 \end{pmatrix}$. 
\end{definition}
It will turn out that up to conjugation, pairs of nontrivially locally conjugate subgroups of $\GL_2(\mathbb{Z}/p^2\mathbb{Z})$ are generated by some equal generators along with two subgroups of $C_s(p^2)$ which are unequal diagonal swaps. 

\begin{proposition} \label{propCartanChoice}
	Let $H_1, H_2 \leq \GL_2(R)$ both be of type $C_s$ or of type $C_{ns}$. Suppose that $H_1$ and $H_2$ are locally conjugate in $\GL_2(R)$, $\tilde{\varphi}(H_i) \not\leq Z(p)$, $\tilde{\varphi}(H_1) = \tilde{\varphi}(H_2)$ and $H_i = \Delta^\perp_{H_i} \rtimes D_{H_i}$.
	\begin{enumerate}
		\item The groups $D_{H_1}$ and $D_{H_2}$ are equal or are diagonal swaps.
		\item If $H_1$ and $H_2$ are nontrivially locally conjugate and $D_{H_1} = D_{H_2}$, then they are of type $C_s$ and one of $\Delta^\perp_{H_1}$ and $\Delta^\perp_{H_2}$ is $\left \langle I + \begin{pmatrix} 0 & 1 \\ 0 & 0 \end{pmatrix} p \right \rangle$ and the other is $\left \langle I + \begin{pmatrix} 0 & 0 \\ 1 & 0 \end{pmatrix} \right \rangle p$. 
	\end{enumerate}
	\begin{proof}
	\begin{enumerate}
		\item Suppose that $H_1$ and $H_2$ are both of type $C_s$. Elements of $H_i$ are expressible as the product $k h$ for unique $k \in \Delta^\perp_{H_i}$ and $h \in D_{H_i}$. Fix $h \in D_{H_i}$ so that $\varphi(h) = \begin{pmatrix} w & 0 \\ 0 & z \end{pmatrix}$, where $w,z \in (\mathbb{Z}/p\mathbb{Z})^\times$ are unequal. Express $h$ in the form $h = \begin{pmatrix} w & 0 \\ 0 & z \end{pmatrix} + \begin{pmatrix} a & 0 \\ 0 & d \end{pmatrix} p$. For $k \in \Delta^\perp_{H_i}$, express $k$ in the form $k = I + \begin{pmatrix} 0 & \beta \\ \gamma & 0 \end{pmatrix} p$. Compute 
		\begin{align*}
			kh &= \left( I + \begin{pmatrix} 0 & \beta \\ \gamma & 0 \end{pmatrix} \right) \left( \begin{pmatrix} w & 0 \\ 0 & z \end{pmatrix} + \begin{pmatrix} a & 0 \\ 0 & d \end{pmatrix} p \right) \\
			   &= \begin{pmatrix} w & 0 \\ 0 & z \end{pmatrix} + \begin{pmatrix} a & \beta z \\ \gamma w & d \end{pmatrix} p
		\end{align*}
		and 
		\begin{align*}
			\trace(kh) &= (w+z) + (a+d) p = \trace(h) \\
			\det(kh) &= (w+ap)(z+dp) = \det(h).
		\end{align*}
		By Theorem \ref{theoremConjClass}, $kh$ and $h$ are in the same conjugacy class, i.e. the conjugacy class of $kh$ does not depend on $k$. The above calculation also shows that the only matrices that could be elements of $D_{H_i}$ and conjugate to $h$ are $h$ itself and $\begin{pmatrix} z & 0 \\ 0 & w \end{pmatrix} + \begin{pmatrix} d & 0 \\ 0 & a \end{pmatrix} p$. \par
		Suppose that there is some $\begin{pmatrix} w & 0 \\ 0 & z \end{pmatrix} \in \varphi(H_1)$ such that $\begin{pmatrix} z & 0 \\ 0 & w \end{pmatrix} \not\in \varphi(H_1)$. Since $H_1 = \Delta^{\perp}_{H_1} \rtimes D_{H_1}$, $\begin{pmatrix} w & 0 \\ 0 & z \end{pmatrix} \in H_1$ by Corollary \ref{corDiagMember}. If $I + \begin{pmatrix} a & 0 \\ 0 & d \end{pmatrix} \in H_1$, then 
		\begin{align*}
			\begin{pmatrix} w & 0 \\ 0 & z \end{pmatrix} \left( I + \begin{pmatrix} a & 0 \\ 0 & d \end{pmatrix} p \right) = \begin{pmatrix} w & 0 \\ 0 & z \end{pmatrix} + \begin{pmatrix} aw & 0 \\ 0 & dz \end{pmatrix} p 
		\end{align*}
		is an element of $H_1$. Since $\begin{pmatrix} z & 0 \\ 0 & w \end{pmatrix} \not \in \varphi(H_2)$, $\begin{pmatrix} w & 0 \\ 0 & z \end{pmatrix} + \begin{pmatrix} aw & 0 \\ 0 & dz \end{pmatrix} p \in H_2$, and so $I + \begin{pmatrix} a & 0 \\ 0 & d \end{pmatrix} p \in H_2$ as well. Therefore, $D_{H_1} \subseteq D_{H_2}$ and $D_{H_2} \subseteq D_{H_1}$ by symmetry. \par
		Now suppose that for all $\begin{pmatrix} w & 0 \\ 0 & z \end{pmatrix} \in \varphi(H_1)$, $\begin{pmatrix} z & 0 \\ 0 & w \end{pmatrix} \in \varphi(H_1)$ as well. For any $I + \begin{pmatrix} a & 0 \\ 0 & d \end{pmatrix} p \in H_1$, $\begin{pmatrix} w & 0 \\ 0 & z \end{pmatrix} + \begin{pmatrix} aw & 0 \\ 0 & dz \end{pmatrix} p \in H_1$, and so $\begin{pmatrix} w & 0 \\ 0 & z \end{pmatrix} + \begin{pmatrix} aw & 0 \\ 0 & dz \end{pmatrix} p$ or $\begin{pmatrix} z & 0 \\ 0 & w \end{pmatrix} + \begin{pmatrix} dz & 0 \\ 0 & aw \end{pmatrix} p$ is in $H_2$. In the former case, $I + \begin{pmatrix} a & 0 \\ 0 & d \end{pmatrix} p \in H_2$. In the latter, $I + \begin{pmatrix} d & 0 \\ 0 & a \end{pmatrix} p \in H_2$. If $\Delta_{H_1}$ is $2$ dimensional, i.e. it is generated by $I + \begin{pmatrix} 1 & 0 \\ 0 & 1 \end{pmatrix} p$ and $I + \begin{pmatrix} 1 & 0 \\ 0 & -1 \end{pmatrix} p$, then $\Delta_{H_2}$ contains both $I + \begin{pmatrix} 1 & 0 \\ 0 & 1 \end{pmatrix} p$ and $I + \begin{pmatrix} 1 & 0 \\ 0 & -1 \end{pmatrix} p$, and so $\Delta_{H_1} = \Delta_{H_2}$. If $\Delta_{H_1}$ is $1$ dimensional, then say that it is generated by $I + \begin{pmatrix} a & 0 \\ 0 & d \end{pmatrix} p$. In particular, $\Delta_{H_2}$ is not $2$ dimensional. If $I + \begin{pmatrix} a & 0 \\ 0 & d \end{pmatrix} p \not\in \Delta_{H_2}$, then $I + \begin{pmatrix} d & 0 \\ 0 & a \end{pmatrix} p \in \Delta_{H_2}$. In this case, let $H'_2$ be the conjugate of $H_2$ via $\begin{pmatrix} 0 & 1 \\ 1 & 0 \end{pmatrix}$. By Lemma \ref{lemmaSkewConj}, $I + \begin{pmatrix} a & 0 \\ 0 & d \end{pmatrix} p \in \Delta_{H'_2}$, and so $\Delta_{H_1} = \Delta_{H'_2}$. If $\Delta_{H_1}$ is $0$ dimensional, then so is $\Delta_{H_2}$, concluding the case where $H_1$ and $H_2$ are both of type $C_s$. \par
		The case where $H_1$ and $H_2$ are both of type $C_{ns}$ is similar. 
		
		\item Suppose that $H_1$ and $H_2$ are nontrivially locally conjugate and $D_{H_1} = D_{H_2}$. Further suppose that $H_1$ and $H_2$ are both of type $C_s$. By Proposition \ref{propNecessaryOne} and Lemma \ref{lemmaLocConjBij}, $\dim( \Delta^{\perp}_{H_1}) = \dim( \Delta^{\perp}_{H_2} )$. Note that $\dim(\Delta^{\perp}_{H_i})$ is not $2$ or $0$ because $H_1 = H_2$ otherwise. Hence, $\dim(\Delta^{\perp}_{H_i}) = 1$. Say that $I + \begin{pmatrix} 0 & b_i \\ c_i & 0 \end{pmatrix} p$ generates $\Delta^{\perp}_{H_i}$. In particular, $b_i$ or $c_i$ is nonzero. By Lemma \ref{lemmaDiagConj}, there is some diagonal matrix $\begin{pmatrix} \alpha_i & 0 \\ 0 & \delta_i \end{pmatrix} $ such that
		\begin{align*}
			\begin{pmatrix} \alpha_i & 0 \\ 0 & \delta_i \end{pmatrix} \left( I + \begin{pmatrix} 0 & b_i \\ c_i & 0 \end{pmatrix} p \right) \begin{pmatrix} \alpha_i & 0 \\ 0 & \delta_i \end{pmatrix}^{-1} = \begin{cases} I + \begin{pmatrix} 0 & 1 \\ 1 & 0 \end{pmatrix} &\text{if } b_ic_i \neq 0 \text{ and } b_ic_i \text{ is a square} \\
			I + \begin{pmatrix} 0 & \epsilon \\ 1 & 0 \end{pmatrix} &\text{if } b_ic_i \text{ is not a square} \\
			I + \begin{pmatrix} 0 & 1 \\ 0 & 0 \end{pmatrix} &\text{if } b_i \neq 0 \\
			I + \begin{pmatrix} 0 & 0 \\ 1 & 0 \end{pmatrix} &\text{if } c_i \neq 0. \end{cases}
		\end{align*}
		By the categorization of subgroups of $\ker \varphi$ up to conjugacy in $\GL_2(\mathbb{Z}/p^2\mathbb{Z})$ and techniques, particularly those involving Proposition \ref{propLocConjKerEquiv}, used to obtain these categorization as discussed in Sections \ref{subsectionPRK}, \ref{subsectionT} and \ref{subsectionFRK}, $b_1c_1$ and $b_2c_2$ must be both zero, both nonzero square or both nonsquares. Replace $H_i$ by its conjugate via $\begin{pmatrix} \alpha_i & 0 \\ 0 & \delta_i \end{pmatrix}$. $D_{H_i}$ is preserved in this conjugation because diagonal matrices multiplicatively commute. Therefore, for $H_1$ and $H_2$ to be nontrivially locally conjugate, one of $D_{H_1}$ and $D_{H_2}$ must be $\left \langle I + \begin{pmatrix} 0 & 1 \\ 0 & 0 \end{pmatrix} p \right \rangle$ and the other must be $\left \langle I + \begin{pmatrix} 0 & 0 \\ 1 & 0 \end{pmatrix} p \right \rangle$. This must have have been true before the conjugation as well. \par
		Now suppose that $H_1$ and $H_2$ are both of type $C_{ns}$. Note that $\Delta^\perp_{H_i}$ consists only of elements of the form $I + \begin{pmatrix} 0 & a - \sqrt{\epsilon} c \\ a + \sqrt{\epsilon} c & 0 \end{pmatrix} p$. Similarly as in the type $C_s$ case, $\dim(\Delta^\perp_{H_1}) = \dim(\Delta^\perp_{H_2}) = 1$. Say that $I + \begin{pmatrix} 0 & a_i + \sqrt{\epsilon} c_i \\ a_i - \sqrt{\epsilon} c_i & 0 \end{pmatrix} p$ generates $\Delta^\perp_{H_i}$. At least one of $a_i$ and $c_i$ is nonzero modulo $p$, and so $\begin{pmatrix} a_2 + \sqrt{\epsilon} c_2 & 0 \\ 0 & a_1 + \sqrt{\epsilon} c_1 \end{pmatrix}$ is invertible. Conjugating $H_1$ by $\begin{pmatrix} a_2 + \sqrt{\epsilon} c_2 & 0 \\ 0 & a_1 + \sqrt{\epsilon} c_1 \end{pmatrix}$ yields $H_2$ by Lemma \ref{lemmaDiagConj}, and so $H_1$ and $H_2$ are conjugate to begin with, a contradiction. Hence, $H_1$ and $H_2$ must both be of type $C_{s}$. 
	\end{enumerate}
	\end{proof}
\end{proposition}

Corollary \ref{corCNS} below categorizes local conjugacy for the subgroups of $\GL_2(\mathbb{Z}/p^2\mathbb{Z})$ whose images under $\varphi$ are contained in $C_{ns}(p)$ but not $Z(p)$.

\begin{corollary}\label{corCNS}
	Let $H_1, H_2 \leq \GL_2(\mathbb{Z}/p^2\mathbb{Z})$ with $\varphi(H_i) \leq C_{ns}(p)$ but $\varphi(H_i) \not\leq Z(p)$. If $H_1$ and $H_2$ are locally conjugate in $\GL_2(\mathbb{Z}/p^2\mathbb{Z})$, then they are conjugate in $\GL_2(\mathbb{Z}/p^2\mathbb{Z})$.
	\begin{proof}
		It suffices to show that $H_1$ and $H_2$ are conjugate in $\GL_2(R)$ by Corollary \ref{corOne}. Replace $H_i$ with its conjugate via $\begin{pmatrix} -\sqrt{\epsilon} & -\epsilon \\ -\sqrt{\epsilon} & \epsilon \end{pmatrix}$ so that $H_i$ is of type $C_{ns}$. Using Proposition \ref{propCartanConj}, further replace $H_i$ with a conjugate such that $H_i$ is still of type $C_{ns}$, $\tilde{\varphi}(H_i)$ is preserved and $H_i = \Delta^\perp_{H_i} \rtimes D_{H_i}$. Proposition \ref{propCartanChoice} shows that $H_1$ and $H_2$ are conjugate in $\GL_2(R)$ as desired.
	\end{proof}
\end{corollary}

Proposition \ref{propCartan} below categorizes local conjugacy for the subgroups of $\GL_2(\mathbb{Z}/p^2\mathbb{Z})$ whose images under $\varphi$ are contained in $C_s(p)$. 

\begin{proposition} \label{propCartan}
	Let $H_1, H_2 \leq \GL_2(\mathbb{Z}/p^2\mathbb{Z})$ with $\varphi(H_i) \leq C_s(p)$. Then, $H_1$ and $H_2$ are nontrivially locally conjugate if and only if they are conjugate to the groups
	\begin{align*}
		\left \langle D, I + \begin{pmatrix} 0 & 1 \\ 0 & 0 \end{pmatrix} p \right \rangle \text{ and } \left \langle D', I + \begin{pmatrix} 0 & 1 \\ 0 & 0 \end{pmatrix} p \right \rangle
	\end{align*}
	in some order, where $D$ and $D'$ are subgroups of $C_s(p^2)$, $D$ and $D'$ are diagonal swaps, and $D \neq D'$. 
	\begin{proof}
		Suppose that $H_1$ and $H_2$ are nontrivially locally conjugate. By Proposition \ref{propNecessaryTwo}, $\varphi(H_1)$ and $\varphi(H_2)$ are locally conjugate. By Theorem \ref{theoremLocConjGL2p}, $\varphi(H_1)$ and $\varphi(H_2)$ are conjugate, and so $H_2$ can be replaced with a conjugate so that $\varphi(H_1) = \varphi(H_2)$. \par
		If $\varphi(H_i) \leq Z(p)$, then $H_1 \cap \ker \varphi$ and $H_2 \cap \ker \varphi$ must be nontrivially locally conjugate. Proposition \ref{propNTLCKer} asserts that $H_1 \cap \ker \varphi$ and $H_2 \cap \ker \varphi$ are, in some order, conjugate to
		\begin{align*}
			\left \langle I + \begin{pmatrix} 1 & 0 \\ 0 & d \end{pmatrix} p, I + \begin{pmatrix} 0 & 1 \\ 0 & 0 \end{pmatrix} p \right \rangle \text{ and }\left \langle I + \begin{pmatrix} d & 0 \\ 0 & 1 \end{pmatrix} p, I + \begin{pmatrix} 0 & 1 \\ 0 & 0 \end{pmatrix} p \right \rangle.
		\end{align*}
		where $d \neq \pm 1$. The discussion in the beginning of Section \ref{sectionCenter} concludes that $H_i = \sigma(\varphi(H_i)) \times (H_i \cap \ker \varphi)$, where $\sigma: Z(p) \rightarrow Z(p^2)$ is the splitting which maps $\begin{pmatrix} w & 0 \\ 0 & w \end{pmatrix} \in Z(p)$ to $\begin{pmatrix} w & 0 \\ 0 & w \end{pmatrix} \in Z(p^2)$. Since $Z(p) \simeq (\mathbb{Z}/p\mathbb{Z})^\times$, $\sigma(\varphi(H_i))$ is cyclic. Say that $\begin{pmatrix} w_0 & 0 \\ 0 & w_0 \end{pmatrix}$ generates $\sigma(\varphi(H_i)$. $H_1$ and $H_2$ are, in some order, conjugate to
		\begin{align*}
			\left \langle \begin{pmatrix} w_0 & 0 \\ 0 & w_0 \end{pmatrix}, I + \begin{pmatrix} 1 & 0 \\ 0 & d \end{pmatrix} p, I + \begin{pmatrix} 0 & 1 \\ 0 & 0 \end{pmatrix} p \right \rangle \text{ and } \left \langle \begin{pmatrix} w_0 & 0 \\ 0 & w_0 \end{pmatrix}, I + \begin{pmatrix} d & 0 \\ 0 & 1 \end{pmatrix} p, I + \begin{pmatrix} 0 & 1 \\ 0 & 0 \end{pmatrix} p \right \rangle
		\end{align*}
		by Lemma \ref{lemmaCenter} and Proposition \ref{propNTLCKer}. Let 
		\begin{align*}
		D = \left \langle \begin{pmatrix} w_0 & 0 \\ 0 & w_0 \end{pmatrix}, I + \begin{pmatrix} 1 & 0 \\ 0 & d \end{pmatrix} p\right \rangle \qquad \text{and} \qquad D' = \left \langle \begin{pmatrix} w_0 & 0 \\ 0 & w_0 \end{pmatrix}, I + \begin{pmatrix} d & 0 \\ 0 & 1 \end{pmatrix} p \right \rangle.
		\end{align*}
		Note that $D' = \begin{pmatrix} 0 & 1 \\ 1 & 0 \end{pmatrix} D \begin{pmatrix} 0 & 1 \\ 1 & 0 \end{pmatrix}^{-1}$ and that $D \neq D'$. Furthermore, $H_1$ and $H_2$ are conjugate to 
		\begin{align*}
			\left \langle D, I + \begin{pmatrix} 0 & 1 \\ 0 & 0 \end{pmatrix} p \right \rangle \qquad \text{and} \qquad \left \langle D', I + \begin{pmatrix} 0 & 1 \\ 0 & 0 \end{pmatrix} p \right \rangle
		\end{align*}
		in some order as desired. \par
		
		Now assume that $\varphi(H_i) \not\leq Z(p)$. Through Proposition \ref{propCartanConj}, replace $H_1$ and $H_2$ so that $\varphi(H_1)$ and $\varphi(H_2)$ are still equal, $\varphi(H_i)$ is still a subgroup of $C_s(p)$ and $H_i = \Delta_{H_i}^\perp \rtimes D_{H_i}$. By Proposition \ref{propCartanChoice}, $H_2$ can be replaced with a conjugate, if necessary, so that $D_{H_1} = D_{H_2}$ as well. Proposition \ref{propCartanChoice} further asserts that one of $\Delta^{\perp}_{H_1}$ and $\Delta^{\perp}_{H_2}$ is $\left \langle I + \begin{pmatrix} 0 & 1 \\ 0 & 0 \end{pmatrix} p \right \rangle$ and that the other is $\left \langle I + \begin{pmatrix} 0 & 0 \\ 1 & 0 \end{pmatrix} p \right \rangle$. Thus, $H_1$ and the conjugate of $H_2$ via $\begin{pmatrix} 0 & 1 \\ 1 & 0 \end{pmatrix}$ are 
		\begin{align*}
			\left \langle D_{H_i}, I + \begin{pmatrix} 0 & 1 \\ 0 & 0 \end{pmatrix} p \right \rangle \qquad \text{and} \qquad \left \langle D_{H_i}, I + \begin{pmatrix} 0 & 0 \\ 1 & 0 \end{pmatrix} p \right \rangle 
		\end{align*}
		in some order. Conjugating the second of these groups by $\begin{pmatrix} 0 & 1 \\ 1 & 0 \end{pmatrix}$ yields
		\begin{align*}
			\left \langle \begin{pmatrix} 0 & 1 \\ 1 & 0 \end{pmatrix} D_{H_i} \begin{pmatrix} 0 & 1 \\ 1 & 0 \end{pmatrix}^{-1}, I + \begin{pmatrix} 0 & 1 \\ 0 & 0 \end{pmatrix} p \right \rangle.
		\end{align*}
		Note that $H_1$ and $H_2$ are conjugate if $D_{H_i} = \begin{pmatrix} 0 & 1 \\ 1 & 0 \end{pmatrix} D_{H_i} \begin{pmatrix} 0 & 1 \\ 1 & 0 \end{pmatrix}^{-1}$. It is also not difficult to see that $H_1$ and $H_2$ are nontrivially locally conjugate otherwise. 
	\end{proof}
\end{proposition}

\subsection{The Normalizer of Cartan Cases} \label{sectionNormalizerCartan}
This section categorizes the subgroups, up to conjugation, of $\GL_2(\mathbb{Z}/p^2\mathbb{Z})$ whose images via $\varphi$ are subgroups of $N(C_s(p))$ or $N(C_{ns}(p))$ that are not contained in $C_s(p)$ or $C_{ns}(p)$. 

\begin{lemma} \label{lemmaNnpm}
	Let $H \leq \GL_2(R)$ be of type $N(C_s(p))$ or $N(C_{ns}(p))$. Suppose that $\tilde{\varphi}(H)$ contains some element of the form $\begin{pmatrix} w & 0 \\ 0 & z \end{pmatrix}$ where $w \neq z$ and another element of the form $\begin{pmatrix} 0 & x \\ y & 0 \end{pmatrix}$. 
	\begin{enumerate}
		\item	Suppose that $H$ is of type $N(C_s(p))$. Then, $\Delta_H$ is generated by both, one or neither of $I + \begin{pmatrix} 1 & 0 \\ 0 & 1 \end{pmatrix} p$ and $I + \begin{pmatrix} 1 & 0 \\ 0 & -1 \end{pmatrix} p$.
		\item Suppose that $H$ is of type $N(C_{ns}(p))$. Then, $\Delta_H$ is generated by both, one or neither of $I + \begin{pmatrix} 1 & 0 \\ 0 & 1 \end{pmatrix} p, I + \begin{pmatrix} \sqrt{\epsilon} & 0 \\ 0 & -\sqrt{\epsilon} \end{pmatrix} p$.
		\item Suppose that $H$ is of type $N(C_s(p))$. Then, $\Delta^\perp_H$ is generated by both, one or neither of $I + \begin{pmatrix} 0 & x \\ y & 0 \end{pmatrix} p$ and $I + \begin{pmatrix} 0 & x \\ -y & 0 \end{pmatrix} p$. \label{lemmaNnpm2}
		\item Suppose that $H$ is of type $N(C_{ns}(p))$. Then, $\Delta^\perp_H$ is generated by both, one or neither of $I + \begin{pmatrix} 0 & x \\ y & 0 \end{pmatrix} p$ and $I + \begin{pmatrix} 0 & x \sqrt{\epsilon} \\ -y\sqrt{\epsilon} & 0 \end{pmatrix} p$.
		\item If $w \neq \pm z$, then $\dim(\Delta^\perp_H) = 2$ or $0$.
	\end{enumerate}
	\begin{proof}
	\begin{enumerate}
		\item If $I + \begin{pmatrix} a & 0 \\ 0 & d \end{pmatrix} p \in H$, then $I + \begin{pmatrix} d & 0 \\ 0 & a \end{pmatrix} p \in H$ by Lemma \ref{lemmaSkewConj}. Suppose that $\dim(\Delta_H) = 1$. Note that $\Delta_H$ must only contain elements of the form $I + \begin{pmatrix} a & 0 \\ 0 & a \end{pmatrix} p$ or only contain elements of the form $I + \begin{pmatrix} a & 0 \\ 0 & -a \end{pmatrix} p$. \par
		\item This follows from the same argument as the last part. 
		\item If $I + \begin{pmatrix} 0 & b \\ c & 0 \end{pmatrix} p \in H$, then $I + \begin{pmatrix} 0 & c \frac{x}{y} \\ b \frac{y}{x} & 0 \end{pmatrix} p \in H$ be Lemma \ref{lemmaSkewConj}. Suppose that $\dim(\Delta^\perp_H) = 1$ and assume that $b$ or $c$ is nonzero. In this case, $c^2 \frac{x}{y} = b^2 \frac{y}{x}$. Since $x$ and $y$ are nonzero and at least one of $b$ and $c$ is nonzero, both $b$ and $c$ are nonzero. Thus, $\frac{x}{y} = \pm \frac{b}{c}$. Since $b,c,x,y \in \mathbb{Z}/p\mathbb{Z}$, we are done. 
		\item Likewise, suppose that $\dim(\Delta^\perp_H) = 1$ and take $I + \begin{pmatrix} 0 & b \\ c & 0 \end{pmatrix} p \in H$ such that $b$ or $c$ is nonzero. Note that $x,y,b,c$ are of the forms 
		\begin{align*}
		x &= \alpha + \beta \sqrt{\epsilon} \\
		y &= \alpha - \beta \sqrt{\epsilon} \\
		b &= \gamma + \delta \sqrt{\epsilon} \\
		c &= \gamma - \delta \sqrt{\epsilon}
		\end{align*}
		for some $\alpha, \beta, \gamma, \delta \in \mathbb{Z}/p\mathbb{Z}$ by Lemma \ref{lemmaEpsilonConj}. If $\frac{x}{y} = \frac{b}{c}$, then $cx = by$, in which case
		\begin{align*}
			 (\alpha \gamma - \epsilon \beta \delta) + (-\alpha \delta + \beta \gamma ) \sqrt{\epsilon} &= (\gamma - \delta \sqrt{\epsilon})(\alpha + \beta \sqrt{\epsilon}) \\
			&= (\gamma + \delta \sqrt{\epsilon})(\alpha - \beta \sqrt{\epsilon}) \\
			&= (\alpha \gamma - \epsilon \beta \delta) + (\alpha \delta - \beta \gamma) \sqrt{\epsilon}.
		\end{align*}
		Thus, $\alpha \delta = \beta \gamma$, i.e. $\frac{x}{b}$ is in $\mathbb{Z}/p\mathbb{Z}$. In this case, $\Delta^\perp_H$ is generated by $I + \begin{pmatrix} 0 & x \\ y & 0 \end{pmatrix} p$. If $\frac{x}{y} = - \frac{b}{c}$, then similarly compute $\alpha \gamma = \epsilon \beta \delta$. Note that $\frac{x\sqrt{\epsilon}}{b} \in \mathbb{Z}/p\mathbb{Z}$, and so $\Delta^\perp_H$ is generated by $I + \begin{pmatrix} 0 & x \sqrt{\epsilon} \\ -y\sqrt{\epsilon} & 0 \end{pmatrix} p$.
		\item This is due to the last two parts and Lemma \ref{lemmaDiagConj}. 
	\end{enumerate}
	\end{proof}
\end{lemma}

\begin{lemma} \label{lemmaNnpmConj}
	Let $H_1, H_2 \leq \GL_2(R)$ be locally conjugate and both of type $N(C_s)$ or both of type $N(C_{ns})$ but not of types $C_s$ or $C_{ns}$. If one of $H_1$ or $H_2$ has an element of the form $\begin{pmatrix} w & 0 \\ 0 & z \end{pmatrix}$ where $w \neq \pm z$, then $H_1$ and $H_2$ are conjugate.
	\begin{proof}
		By Propositions \ref{propNecessaryTwo} and \ref{propSubGpGL2p} and Theorem \ref{theoremLocConjGL2p}, $\tilde{\varphi}(H_1)$ and $\tilde{\varphi}(H_2)$ are conjugate. Replace $H_2$ with a conjugate so that $\tilde{\varphi}(H_1) = \tilde{\varphi}(H_2)$. Note that if $H_i$ is of type $N(C_{ns})$, then such a conjugation can be done by an element of $\begin{pmatrix} -\sqrt{\epsilon} & -\epsilon \\ -\sqrt{\epsilon} & \epsilon \end{pmatrix} \GL_2(\mathbb{Z}/p\mathbb{Z}) \begin{pmatrix} -\sqrt{\epsilon} & -\epsilon \\ -\sqrt{\epsilon} & \epsilon \end{pmatrix}^{-1}$. In particular, $H_2 \cap \ker \tilde{\varphi}$ is still a subgroup of $K$ after the conjugation. \par
		By Lemma \ref{lemmaNnpm}, $\Delta^\perp_{H_i}$ is $2$ or $0$ dimensional. If $\Delta^\perp_{H_1}$ is $2$ dimensional but $\Delta^\perp_{H_2}$ is $0$ dimensional, then $\Delta_{H_1}$ must be $0$ dimensional and $\Delta_{H_2}$ must be $2$ dimensional. However, $\Delta_{H_2}$ would then have elements of nonzero trace whereas all of the elements of $\Delta_{H_1}$ has $0$ trace, which is a contradiction. Hence, $\Delta^\perp_{H_1} = \Delta^\perp_{H_2}$. Lemma \ref{lemmaNnpm} also asserts that $\Delta_{H_i}$ is generated by both, one or neither of $I + \begin{pmatrix} 1 & 0 \\ 0 & 1 \end{pmatrix} p$ and a matrix of the form $I + \begin{pmatrix} a & 0 \\ 0 & -a \end{pmatrix} p$. Note that $I + \begin{pmatrix} 1 & 0 \\ 0 & 1 \end{pmatrix} p$ is an element of $H_1$ if and only if it is an element of $H_2$. Therefore, $H_1 \cap \ker \tilde{\varphi} = H_2 \cap \ker \tilde{\varphi}$. Recall that if $H_i$ is of type $N(C_{ns})$, then its conjugate via $\begin{pmatrix} -\sqrt{\epsilon} & -\epsilon \\ -\sqrt{\epsilon} & \epsilon \end{pmatrix}^{-1}$ is a subgroup of $N(C_{ns}(p))$. Since $N(C_s(p))$ and $N(C_{ns}(p))$ both have orders which are indivisible by $p$, Proposition \ref{propSchurZassenhaus} shows that $H_1$ and $H_2$ are conjugate. 
	\end{proof}
\end{lemma}

\begin{proposition}
	Let $H_1, H_2 \leq \GL_2(R)$ be locally conjugate and both of type $N(C_s)$ or both of type $N(C_{ns})$ but not of types $C_s$ or $C_{ns}$. $H_1$ and $H_2$ are conjugate. 
	\begin{proof}
		If all of the diagonal elements of $\tilde{\varphi}(H_i)$ are scalar, then there is some $r \in R^\times$ such that $\frac{x}{y} = r$ for all $\begin{pmatrix} 0 & x \\ y & 0 \end{pmatrix} \in \tilde{\varphi}(H_i)$. Furthermore, whether or not $H_i$ is of type $N(C_s)$, it is not difficult to see that $r$ is a square in $R$. Replace $H_i$ by its conjugate via $\begin{pmatrix} \sqrt{r} & \sqrt{r} \\ 1 & -1 \end{pmatrix}^{-1}$. All elements of $\tilde{\varphi}(H_i)$ are then diagonal and it is not difficult to see that $H_i$ is in fact of type $C_s$ or of type $C_{ns}$. This possibility was already considered in Section \ref{sectionCartan}. Assume that some diagonal elements of $\tilde{\varphi}(H_i)$ are nonscalar. \par
		By Lemma \ref{lemmaNnpmConj}, it suffices to show that $H_1$ and $H_2$ are conjugate given that all nonscalar diagonal elements of $\tilde{\varphi}(H_i)$ are of the form $\begin{pmatrix} w & 0 \\ 0 & -w \end{pmatrix}$. In this case, there is some $r \in R^\times$ such that $\frac{x}{y} = \pm r$ for all $\begin{pmatrix} 0 & x \\ y & 0 \end{pmatrix} \in \tilde{\varphi(H_i)}$. Moreover, 
		\begin{align*}
			\begin{pmatrix} w & 0 \\ 0 & -w \end{pmatrix} \begin{pmatrix} 0 & x \\ y & 0 \end{pmatrix} \begin{pmatrix} w & 0 \\ 0 & -w \end{pmatrix}^{-1} \begin{pmatrix} 0 & x \\ y & 0 \end{pmatrix}^{-1} = -I
		\end{align*}
		and so $-I \in \tilde{\varphi}(H_i)$. Similarly as in Lemma \ref{lemmaNnpmConj}, $\tilde{\varphi}(H_1)$ and $\tilde{\varphi}(H_2)$ are conjugate by Propositions \ref{propNecessaryTwo} and \ref{propSubGpGL2p} and Theorem \ref{theoremLocConjGL2p}. Replace $H_2$ with a conjugate so that $\tilde{\varphi}(H_1) = \tilde{\varphi}(H_2)$. If $H_i$ is of type $N(C_{ns})$, then the replacement can be done in a way such that $H_2 \cap \ker \tilde{\varphi}$ is still a subgroup of $K$ after the conjugation. \par
		Suppose that there is some $\begin{pmatrix} w & 0 \\ 0 & -w \end{pmatrix} \in \tilde{\varphi}(H_i)$ which is not conjugate to any element of $H_i$ of the form $\begin{pmatrix} 0 & x \\ y & 0 \end{pmatrix}$. By Corollary \ref{corDiagMember}, $H_i$ has an element $h_i$ of the form $h_i = \begin{pmatrix} w & 0 \\ 0 & -w \end{pmatrix} + \begin{pmatrix} 0 & b_i \\ c_i & 0 \end{pmatrix} p$. For $k = I + \begin{pmatrix} \alpha & \beta \\ \gamma & \delta \end{pmatrix} p \in H_1$, compute 
		\begin{align*}
		h_1 k = \begin{pmatrix} w & 0 \\ 0 & -w \end{pmatrix} + \begin{pmatrix} \alpha w & \beta w + b_1 \\ - \gamma w + c_1 & - \delta w \end{pmatrix} p
		\end{align*}
		An element of $H_2$ which is conjugate to $h_1 k$ must be of the form $\begin{pmatrix} w & 0 \\ 0 & -w \end{pmatrix} + \begin{pmatrix} \alpha w & * \\ * & -\delta w \end{pmatrix} p$ or $\begin{pmatrix} -w & 0 \\ 0 & w \end{pmatrix} + \begin{pmatrix} -\delta w & * \\ * & \alpha w \end{pmatrix} p$. Thus, if $I + \begin{pmatrix} \alpha & 0 \\ 0 & \delta \end{pmatrix} p \in H_1$ then $I + \begin{pmatrix} \alpha & 0 \\ 0  & \delta \end{pmatrix} p$ or $I + \begin{pmatrix} \delta & 0 \\ 0 & \alpha \end{pmatrix} p$ is an element of $H_2$. By Lemma \ref{lemmaNnpm}, $\Delta_{H_i}$ is generated by both, one or neither of $I + \begin{pmatrix} 1 & 0 \\ 0 & 1 \end{pmatrix} p$ and $I + \begin{pmatrix} a & 0 \\ 0 & -a \end{pmatrix} p$ where $a = 1$ or $\sqrt{\epsilon}$. Thus, $\Delta_{H_1} = \Delta_{H_2}$ and so $\dim(\Delta^\perp_{H_1}) = \dim(\Delta^\perp_{H_2})$. If $\dim(\Delta^\perp_{H_i}) = 2$ or $0$, then $\Delta^\perp_{H_1} = \Delta^\perp_{H_2}$, in which case $H_1 = H_2$. \par
		Fix an element $\begin{pmatrix} 0 & x \\ y & 0 \end{pmatrix}$ of $\tilde{\varphi}(H_i)$ and suppose that $\dim(\Delta^\perp_{H_i}) = 1$. Fix $g_i \in H_i$ to be of the form $g_i = \begin{pmatrix} 0 & x \\ y & 0 \end{pmatrix} + \begin{pmatrix} a_i & b_i \\ c_i & d_i \end{pmatrix} p$. An element of $H_2$ which is conjugate to $g_1$ must be an element of $\tilde{\varphi}^{-1}\left( \begin{pmatrix} 0 & x' \\ y' & 0 \end{pmatrix} \right)$ for some $\begin{pmatrix} 0 & x' \\ y' & 0 \end{pmatrix} \in \tilde{\varphi}(H_i)$. In particular, $-xy = \det \begin{pmatrix} 0 & x \\ y & 0 \end{pmatrix} = \det \begin{pmatrix} 0 & x' \\ y' & 0 \end{pmatrix} = -x'y'$ and $\frac{y}{x} = \pm \frac{y'}{x'}$, and so $\begin{pmatrix} 0 & x' \\ y' & 0 \end{pmatrix}$ is one of $\begin{pmatrix} 0 & \pm x \\ \pm y & 0 \end{pmatrix}$ or $\begin{pmatrix} 0 & \pm j x \\ \mp j y & 0 \end{pmatrix}$ where $j$ is a square root of $-1$. If $H_i$ is of type $N(C_s)$, then $j \in \mathbb{Z}/p\mathbb{Z}$. Moreover, using Lemma \ref{lemmaEpsilonConj} shows that $j \not\in \mathbb{Z}/p\mathbb{Z}$ if $H_i$ is of type $N(C_{ns})$. \par
		Suppose that one of $\begin{pmatrix} 0 & \pm j x \\ \mp j y & 0 \end{pmatrix}$ is an element of $\tilde{\varphi}(H_i)$. Since $-I \in \tilde{\varphi}(H_i)$, both of them are elements of $\tilde{\varphi}(H_i)$. Moreover, 
		\begin{align*}
			\begin{pmatrix} 0 & \pm j x \\ \mp j y & 0 \end{pmatrix} \begin{pmatrix} 0 & \frac{1}{y} \\ \frac{1}{x} & 0 \end{pmatrix} = \begin{pmatrix} \pm j & 0 \\ 0 & \mp j \end{pmatrix},
		\end{align*}
		and so for all $\begin{pmatrix} 0 & x' \\ y' & 0 \end{pmatrix} \in \tilde{\varphi}(H_i)$, $\begin{pmatrix} 0 & \pm j x' \\ \mp j y' & 0 \end{pmatrix}$ are elements of $\tilde{\varphi}(H_i)$. Note that Lemma \ref{lemmaDiagConj} shows that $\tilde{\varphi}(H_i)$ is preserved under conjugation via $\begin{pmatrix} j & 0 \\ 0 & 1 \end{pmatrix}$. Further suppose that $\Delta^\perp_{H_1} \neq \Delta^\perp_{H_2}$. If $H_i$ is of type $N(C_s)$, then say that, without loss of generality, $\Delta^\perp_{H_1}$ is generated by $I + \begin{pmatrix} 0 & x \\ y & 0 \end{pmatrix} p$ and that $\Delta^\perp_{H_2}$ is generated by $I + \begin{pmatrix} 0 & x \\ -y & 0\end{pmatrix} p$ using Lemma \ref{lemmaNnpm}. Replace $H_2$ with its conjugate via $\begin{pmatrix} j & 0 \\ 0 & 1 \end{pmatrix}$. Since $j \in \mathbb{Z}/p\mathbb{Z}$, $\Delta^\perp_{H_1}$ and $\Delta^\perp_{H_2}$ are now equal. By Proposition \ref{propSchurZassenhaus}, $H_1$ and $H_2$ are conjugate. Similarly, if $H_i$ is of type $N(C_{ns})$, then one can replace $H_2$ with its conjugate via $\begin{pmatrix} j & 0 \\ 0 & 1 \end{pmatrix}$. Since $j\sqrt{\epsilon}$ is an element of $\mathbb{Z}/p\mathbb{Z}$, $\Delta^\perp_{H_1}$ and $\Delta^\perp_{H_2}$ are now equal and so $H_1$ and $H_2$ are conjugate by Proposition \ref{propSchurZassenhaus}. \par
		Suppose that neither of $\begin{pmatrix} 0 & \pm j x \\ \mp j y & 0 \end{pmatrix}$ is in $\tilde{\varphi}(H_i)$. For $k = I + \begin{pmatrix} \alpha & \beta \\ \gamma & \delta \end{pmatrix} p \in H_1$, compute
		\begin{align*}
			g_1k = \begin{pmatrix} 0 & x \\ y & 0 \end{pmatrix} + \begin{pmatrix} a_i + \gamma x & b_i + \delta x \\ c_i + \alpha y & d_i + \beta y \end{pmatrix}p .
		\end{align*}
		If $\Delta^\perp_{H_i}$ is generated by $I + \begin{pmatrix} 0 & x \\ y & 0 \end{pmatrix} p$, then $\trace(g_1k)$ is of the form $(a_i + d_i + 2nxy u)p$ where $ n \in \mathbb{Z}/p\mathbb{Z}$ and $u = 1$ or $\sqrt{\epsilon}$, in particular, when $\beta = nx u$ and $\gamma = ny u$. Otherwise, $\trace(g_1k)$ can only have trace $(a_1+d_1)p$. Since the only elements of $\tilde{\varphi}(H_i)$ that are conjugate to $\begin{pmatrix} 0 & x \\ y & 0 \end{pmatrix}$ are $\begin{pmatrix} 0 & \pm x \\ \pm y & 0 \end{pmatrix}$, one can see that $\Delta^\perp_{H_1} = \Delta^\perp_{H_2}$. Thus, $H_1$ and $H_2$ are conjugate by Proposition \ref{propSchurZassenhaus}. This concludes the case where there is some $\begin{pmatrix} w & 0 \\ 0 & -w \end{pmatrix} \in \tilde{\varphi}(H_i)$ which is not conjugate to any element of $\tilde{\varphi}(H_i)$ of the form $\begin{pmatrix} 0 & x \\ y & 0 \end{pmatrix}$. \par
		Now assume that every $\begin{pmatrix} w & 0 \\ 0 & -w \end{pmatrix} \in \tilde{\varphi}(H_i)$ is conjugate to an element of $\tilde{\varphi}(H_i)$ of the form $\begin{pmatrix} 0 & x \\ y & 0 \end{pmatrix}$. In particular, $xy = w^2$. Replace $H_i$ with its conjugate via $\begin{pmatrix} \frac{w}{x} & 0 \\ 0 & \frac{x}{w} \end{pmatrix}$ so that $\begin{pmatrix} 0 & w \\ w & 0 \end{pmatrix} \in \tilde{\varphi}(H_i)$ instead of $\begin{pmatrix} 0 & x \\ y & 0 \end{pmatrix}$. \par
		Suppose that $H_i$ is of type $N(C_s(p))$. By Lemma \ref{lemmaNnpm}, some matrices among $I + \begin{pmatrix} 1 & 0 \\ 0 & 1 \end{pmatrix} p$, $I + \begin{pmatrix} 1 & 0 \\ 0 & -1 \end{pmatrix} p$, $I + \begin{pmatrix} 0 & 1 \\ 1 & 0 \end{pmatrix} p$ and $I + \begin{pmatrix} 0 & 1 \\ -1 & 0 \end{pmatrix} p$ together generate $H_i \cap \ker \varphi$. Since $I + \begin{pmatrix} 1 & 0 \\ 0 & 1 \end{pmatrix} p$ is an element of $H_1$ if and only if it is an element of $H_2$, $H_1 \cap \ker \varphi$ and $H_2 \cap \ker \varphi$ have the same number of matrices among $I + \begin{pmatrix} 1 & 0 \\ 0 & -1 \end{pmatrix} p$, $I + \begin{pmatrix} 0 & 1 \\ 1 & 0 \end{pmatrix} p$ and $I + \begin{pmatrix} 0 & 1 \\ -1 & 0 \end{pmatrix} p$. \par 
		Suppose that $H_i \cap \ker \varphi$ has only one among the three matrices. If $H_1 \cap \ker \varphi$ has $I + \begin{pmatrix} 1 & 0 \\ 0 & -1 \end{pmatrix} p$ and that $H_2 \cap \ker \varphi$ has $I + \begin{pmatrix} 0 & 1 \\ 1 & 0 \end{pmatrix} p$, then replace $H_2$ with its conjugate via $\begin{pmatrix} 1 & 1 \\ -1 & 1 \end{pmatrix}$. Compute
		\begin{align*}
			\begin{pmatrix} 1 & 1 \\ -1 & 1 \end{pmatrix} \begin{pmatrix} w & 0 \\ 0 & -w \end{pmatrix} \begin{pmatrix} 1 & 1 \\ -1 & 1 \end{pmatrix}^{-1} &= \begin{pmatrix} 0 & -w \\ -w & 0 \end{pmatrix} \\
			\begin{pmatrix} 1 & 1 \\ -1 & 1 \end{pmatrix} \begin{pmatrix} 0 & w \\ w & 0 \end{pmatrix} \begin{pmatrix} 1 & 1 \\ -1 & 1 \end{pmatrix}^{-1} & = \begin{pmatrix} w & 0 \\ 0 & -w \end{pmatrix}
		\end{align*}
		and recall that $-I \in \varphi(H_i)$. Therefore $H_1 \cap \ker \varphi = H_2 \cap \ker \varphi$ and $\varphi(H_2)$ is preserved after the conjugation. Now say that $H_1 \cap \ker \varphi$ has $I + \begin{pmatrix} 1 & 0 \\ 0 & -1 \end{pmatrix} p$ and that $H_2 \cap \ker \varphi$ has $I + \begin{pmatrix} 0 & 1 \\ -1 & 0 \end{pmatrix} p$. Whether or not $I + \begin{pmatrix} 1 & 0 \\ 0 & 1 \end{pmatrix} p \in \varphi(H_i)$, it is not difficult to see that $-1$ must be a square in $\mathbb{Z}/p\mathbb{Z}$ using Proposition \ref{propLocConjKerEquiv}. Again, say that $j \in \mathbb{Z}/p\mathbb{Z}$ is a square root of $-1$. If neither of $\begin{pmatrix} 0 & \pm j w \\ \mp j w & 0 \end{pmatrix}$ is in $\varphi(H_i)$, then compute
		\begin{align*}
			\begin{pmatrix} 1 & j \\ 1 & -j \end{pmatrix} \begin{pmatrix} w & 0 \\ 0 & -w \end{pmatrix} \begin{pmatrix} 1 & j \\ 1 & -j \end{pmatrix}^{-1} &= \begin{pmatrix} 0 & w \\ w & 0 \end{pmatrix} \\
			\begin{pmatrix} 1 & j \\ 1 & -j \end{pmatrix} \begin{pmatrix} 0 & w \\ w & 0 \end{pmatrix} \begin{pmatrix} 1 & j \\ 1 & -j \end{pmatrix}^{-1} &= \begin{pmatrix} 0 & wj \\ -wj & 0 \end{pmatrix}.
		\end{align*}
		Thus, replacing $H_2$ with its conjugate via $\begin{pmatrix} 1 & j \\ 1 & -j \end{pmatrix}$ makes all diagonal elements of $\varphi(H_2)$ scalar, which is a case that was already discussed. If $\begin{pmatrix} 0 & \pm j x \\ \mp j y & 0 \end{pmatrix} \in \varphi(H_i)$, then further compute
		\begin{align*}
			\begin{pmatrix} 1 & j \\ 1 & -j \end{pmatrix} \begin{pmatrix} 0 & j \\ -j & 0 \end{pmatrix} \begin{pmatrix} 1 & j \\ 1 & -j \end{pmatrix}^{-1} = \begin{pmatrix} 1 & 0 \\ 0 & -1 \end{pmatrix}.
		\end{align*}
		Replacing $H_2$ with its conjugate via $\begin{pmatrix} 1 & j \\ 1 & -j \end{pmatrix}$ therefore makes $H_1 \cap \ker \varphi$ and $H_2 \cap \ker \varphi$ equal and preserves $\varphi(H_2)$. Thus, $H_1$ and $H_2$ are conjugate, no matter which of the three matrices they have. \par
		Suppose that $H_i \cap \ker \varphi$ has two among the matrices $I + \begin{pmatrix} 1 & 0 \\ 0 & -1 \end{pmatrix} p$, $I + \begin{pmatrix} 0 & 1 \\ 1 & 0 \end{pmatrix} p$ and $I + \begin{pmatrix} 0 & 1 \\ - 1 & 0 \end{pmatrix} p$. If $I + \begin{pmatrix} 1 & 0 \\ 0 & -1 \end{pmatrix} p, I + \begin{pmatrix} 0 & 1 \\ 1 & 0 \end{pmatrix} p \in H_1 \cap \ker \varphi$ and $I + \begin{pmatrix} 1 & 0 \\ 0 & -1 \end{pmatrix} p, I + \begin{pmatrix} 0 & 1 \\ -1 & 0 \end{pmatrix} p \in H_2 \cap \ker \varphi$, then $-1$ is a square modulo $p$. Similarly as before, we can assume that $\begin{pmatrix} 0 & \pm j w \\ \mp j w & 0 \end{pmatrix} \in \varphi(H_i)$. Replacing $H_2$ with its conjugate via $\begin{pmatrix} 1 & j \\ 1 & -j \end{pmatrix}$ makes $H_1 \cap \ker \varphi$ and $H_2 \cap \ker \varphi$ equal and preserves $\varphi(H_2)$, and so $H_1$ and $H_2$ are conjugate. If $I + \begin{pmatrix} 1 & 0 \\ 0 & -1 \end{pmatrix} p, I + \begin{pmatrix} 0 & 1 \\ -1 & 0 \end{pmatrix} p \in H_1 \cap \ker \varphi$ and $I + \begin{pmatrix} 0 & 1 \\ 1 & 0 \end{pmatrix}, I + \begin{pmatrix} 0 & 1 \\ -1 & 0 \end{pmatrix} p \in H_2 \cap \ker \varphi$, then compute
		\begin{align*}
			\begin{pmatrix} 1 & 1 \\ -1 & 1 \end{pmatrix} \begin{pmatrix} 0 & 1 \\ -1 & 0 \end{pmatrix} \begin{pmatrix} 1 & 1 \\ -1 & 1 \end{pmatrix}^{-1} = \begin{pmatrix} 0 & 1 \\ -1 & 0 \end{pmatrix}.
		\end{align*}
		Replacing $H_2$ with its conjugate via $\begin{pmatrix} 1 & 1 \\ - 1 & 1 \end{pmatrix}$ makes $H_1 \cap \ker \varphi = H_2 \cap \ker \varphi$ and preserves $\varphi(H_2)$, and so $H_1$ and $H_2$ are conjugate. \par
		If $H_i \cap \ker \varphi$ has all three of $I + \begin{pmatrix} 1 & 0 \\ 0 & -1 \end{pmatrix} p, I + \begin{pmatrix} 0 & 1 \\ 1 & 0 \end{pmatrix}p$ and $I + \begin{pmatrix} 0 & 1 \\ -1 & 0 \end{pmatrix} p$, then $H_1 \cap \ker \varphi = H_2 \cap \ker \varphi$, in which case $H_1$ and $H_2$ are conjugate.
	\end{proof}
\end{proposition}

\section{The Borel case} \label{sectionBorel}
	This section categorizes the subgroups, up to conjugation, of $\GL_2(\mathbb{Z}/p^2\mathbb{Z})$ whose images under $\varphi$ are subgroups of $B(p)$ containing $\begin{pmatrix} 1 & 1 \\ 0 & 1 \end{pmatrix}$.  \par
	As with Lemma \ref{lemmaTMemberPG3}, the types of such subgroups that can arise differ between the case when $p > 3$ and the case when $p = 3$. We consider the case where $p > 3$ first. \par
	For $H \leq \GL_2(\mathbb{Z}/p^2\mathbb{Z})$, let $C_H$ denote the image of $\varphi(H) \cap C_s(p)$ under the splitting $C_s(p) \hookrightarrow C_s(p^2)$ which maps $\begin{pmatrix} w & 0 \\ 0 & z \end{pmatrix} \in C_s(p)$ to $\begin{pmatrix} w & 0 \\ 0 & z \end{pmatrix} \in C_s(p^2)$. 

\begin{lemma} \label{lemmaBpConjOne}
	Suppose that $p > 3$. Let $H \leq \GL_2(\mathbb{Z}/p^2\mathbb{Z})$ with $\varphi(H) \leq B(p)$ and $\begin{pmatrix} 1 & 1 \\ 0 & 1 \end{pmatrix} \in \varphi(H)$. There is a conjugate $H'$ of $H$ such that $H' = \langle \tau, H' \cap \ker \varphi \rangle \rtimes C_{H'}$, where $\tau \in \GL_2(\mathbb{Z}/p^2\mathbb{Z})$ is a matrix of the form $\begin{pmatrix} 1 & 1 \\ 0 & 1 \end{pmatrix} + \begin{pmatrix} a & 0 \\ c & a \end{pmatrix} p$ for some $a,c \in \mathbb{Z}/p\mathbb{Z}$. 
	\begin{proof}
	
		\cite[Lemma 3.3]{Sutherland} shows that $\varphi(H)$ is expressible as the internal semidirect product $\left \langle \begin{pmatrix} 1 & 1 \\ 0 & 1 \end{pmatrix} \right \rangle \rtimes \langle \varphi(H) \cap C_s(p) \rangle$. Suppose that $\varphi(H) \cap C_s(p) \leq Z(p)$. For every element $h$ of $H$, $\varphi(h)$ is expressible as $\varphi(h) = \begin{pmatrix} 1 & 1 \\ 0 & 1 \end{pmatrix}^l z$ for some $l \in \mathbb{Z}$ and $z \in \varphi(H) \cap C_s(p)$. Say that $z = \begin{pmatrix} w & 0 \\ 0 & w \end{pmatrix}$, in which case $\begin{pmatrix} w & 0 \\ 0 & w \end{pmatrix} \in H$ by Corollary \ref{corScalarMember}. Fix any $\tau \in H$ such that $\varphi(\tau) = \begin{pmatrix} 1 & 1 \\ 0 & 1 \end{pmatrix}$. $h$ is then expressible as $\tau^l \begin{pmatrix} w & 0 \\ 0 & w \end{pmatrix} k = \tau^l k \begin{pmatrix} w & 0 \\ 0 & w \end{pmatrix}$ for some $k \in H \cap \ker \varphi$. Moreover, it is not difficult to see that $\langle \tau, H \cap \ker \varphi \rangle$ is normal in $H$ and that $\langle \tau, H \cap \ker \varphi \rangle \cap C_{H} = \langle I \rangle$. Therefore, $H$ is the internal semidirect product $H = \langle \tau, H \cap \ker \varphi \rangle \rtimes C_{H}$. Express $\tau$ as $\tau = \begin{pmatrix} 1 & 1 \\ 0 & 1 \end{pmatrix} + \begin{pmatrix} a_0 & b_0 \\ c & d_0 \end{pmatrix} p$. Calculate 
		\begin{align*}
			\det(\tau) = 1 + (a_0 + d_0 - \gamma) p \qquad \text{and} \qquad \trace(\tau) = 2 + (a_0 + d_0) p
		\end{align*}
		By Theorem \ref{theoremConjClass}, $\tau$ is conjugate to $\tau' = \begin{pmatrix} 1 & 1 \\ 0 & 1 \end{pmatrix} + \begin{pmatrix} a & 0 \\ c & a \end{pmatrix} p$, where $a = \frac{a_0+d_0}{2}$. Say that $\tau' = g\tau g^{-1}$ for $g \in \GL_2(\mathbb{Z}/p^2\mathbb{Z})$ and let $H' = gHg^{-1}$. One can then express $H'$ as the internal semidirect product $\langle \tau, H' \cap \ker \varphi \rangle \rtimes C_{H'}$. \par
		
		Now assume that $\varphi(H) \cap C_s(p) \not \leq Z(p)$. By Proposition \ref{propCartanConj}, replace $H$ with a conjugate so that $\varphi(H) \cap C_s(p)$ is preserved and $H \cap \varphi^{-1}(C_s(p)) = \Delta^{\perp}_H \rtimes D_H$. Say that the conjugation is done via $g \in \GL_2(\mathbb{Z}/p^2\mathbb{Z})$. The conjugation preserves $\varphi(H) \cap C_s(p)$, a diagonal subgroup of $\GL_2(\mathbb{Z}/p^2\mathbb{Z})$ which does not lie in the center. It is not difficult to show that $\varphi(g)$ is therefore an element of $N(C_s(p))$. Furthermore, if $\varphi(H)$ is not a subgroup of $B(p)$, then $\begin{pmatrix} 0 & 1 \\ 1 & 0 \end{pmatrix} \varphi(H) \begin{pmatrix} 0 & 1 \\ 1 & 0 \end{pmatrix}^{-1}$ is a subgroup of $B(p)$. If necessary, replace $H$ with such a conjugate so that $\varphi(H) \leq B(p)$ and $H \cap \varphi^{-1}(C_s(p))$ is still $\Delta^{\perp}_H \rtimes D_H$. In particular, the image of the splitting $C_s(p) \hookrightarrow C_s(p^2)$ restricted to $\varphi(H) \cap C_s(p)$ maps isomorphically onto $C_H$. Similarly as before, $H$ is the internal semidirect product $\langle \tau, H \cap \ker \varphi \rangle \rtimes C_H$ for any $\tau \in H$ satisfying $\varphi(\tau) = t$. \par
		Since $p > 3$ by assumption, Lemma \ref{lemmaTMemberPG3} asserts that $I + \begin{pmatrix} 0 & 1 \\ 0 & 0 \end{pmatrix} p \in H$. Suppose that $\tau$ is of the form $\tau = \begin{pmatrix} 1 & 1 \\ 0 & 1 \end{pmatrix} + \begin{pmatrix} a & b \\ c & a \end{pmatrix} p$. By Lemma \ref{lemmaTProdOne}, there is some $\tau' \in H$ of the form $\tau' = \begin{pmatrix} 1 & 1 \\ 0 & 1 \end{pmatrix} + \begin{pmatrix} a & 0 \\ c & a \end{pmatrix} p$, in which case $H = \langle \tau', H \cap \ker \varphi \rangle \rtimes C_H$ and we are done.  \par
		Now assume that $\tau$ is of the form $\tau = \begin{pmatrix} 1 & 1 \\ 0 & 1 \end{pmatrix} + \begin{pmatrix} a & b \\ c & d \end{pmatrix} p$ where $a \neq d$. By Lemma \ref{lemmaConvolutedBp}, there is some $I + \begin{pmatrix} \alpha & 0 \\ 0 & \delta \end{pmatrix} p \in H$ where $\alpha \neq \delta$. Moreover,
		\begin{align*}
			\tau \left( I + \begin{pmatrix} \alpha & 0 \\ 0 & \delta \end{pmatrix} p \right)^{ - \frac{a-d}{\alpha - \delta}} = \begin{pmatrix} 1 & 1 \\ 0 & 1 \end{pmatrix} + \begin{pmatrix} a - \alpha \frac{a-d}{\alpha - \delta} & * \\ c & d - \delta \frac{a-d}{\alpha- \delta} \end{pmatrix} p
		\end{align*}
		by Lemma \ref{lemmaTProdTwo}. One can further show that $a - \alpha \frac{a-d}{\alpha - \delta} = d - \delta \frac{a-d}{\alpha- \delta}$. By Lemma \ref{lemmaTProdOne}, there is some $\tau' \in H$ of the form $\tau' = \begin{pmatrix} 1 & 1 \\ 0 & 1 \end{pmatrix} + \begin{pmatrix} a' & 0 \\ c' & a' \end{pmatrix} p$, which yields the desired result. 
	\end{proof}
\end{lemma}

\begin{lemma} \label{lemmaKerBp}
	Suppose that $p > 3$. Let $H \leq \GL_2(\mathbb{Z}/p\mathbb{Z})$ such that $\varphi(H) \leq B(p)$ and $\begin{pmatrix} 1 & 1 \\ 0 & 1 \end{pmatrix} \in \varphi(H)$. $H \cap \ker \varphi$ is one of the following:
	\begin{enumerate}
		\item $\left \langle I + \begin{pmatrix} a & 0 \\ 0 & d \end{pmatrix} p, I + \begin{pmatrix} 0 & 1 \\ 0 & 0 \end{pmatrix} p \right \rangle$ for some $a,d \in \mathbb{Z}/p\mathbb{Z}$ 
		\item $\left \langle I + \begin{pmatrix} 1 & 0 \\ 0 & 0 \end{pmatrix} p, I + \begin{pmatrix} 0 & 0 \\ 0 & 1 \end{pmatrix} p, I + \begin{pmatrix} 0 & 1 \\ 0 & 0 \end{pmatrix} p \right \rangle$
		\item $\left \langle I + \begin{pmatrix} a & 0 \\ c & a \end{pmatrix} p, I + \begin{pmatrix} 1 & 0 \\ 0 & -1 \end{pmatrix} p, I + \begin{pmatrix} 0 & 1 \\ 0 & 0 \end{pmatrix} p \right \rangle$ for some $a,c \in \mathbb{Z}/p\mathbb{Z}$ such that $c \neq 0$
		\item $\ker \varphi$. 
	\end{enumerate}
	\begin{proof}
		By Lemma \ref{lemmaTMemberPG3}, $I + \begin{pmatrix} 0 & 1 \\ 0 & 0 \end{pmatrix} p \in H$. Suppose that there is some element $k \in H$ of the form $k = I + \begin{pmatrix} a & b \\ c & d \end{pmatrix} p$ such that $c \neq 0$. By Lemma \ref{lemmaTMember}, $I + \begin{pmatrix} 1 & 0 \\ 0 & -1 \end{pmatrix} p \in H$, and so $I + \begin{pmatrix} \frac{a+d}{2} & 0 \\ c & \frac{a+d}{2} \end{pmatrix} p \in H$ as well. Therefore, $H$ contains $\left \langle I + \begin{pmatrix} a' & 0 \\ c & a' \end{pmatrix} p, I + \begin{pmatrix} 1 & 0 \\ 0 & -1 \end{pmatrix} p, I + \begin{pmatrix} 0 & 1 \\ 0 & 0 \end{pmatrix} p \right \rangle$, where $a' = \frac{a+d}{2}$, and so $H = \left \langle I + \begin{pmatrix} a' & 0 \\ c & a' \end{pmatrix} p, I + \begin{pmatrix} 1 & 0 \\ 0 & -1 \end{pmatrix} p, I + \begin{pmatrix} 0 & 1 \\ 0 & 0 \end{pmatrix} p \right \rangle$ or $H = \ker \varphi$. \par
		Now assume that every $k \in H$ is of the form $k = I + \begin{pmatrix} a & b \\ 0 & d \end{pmatrix} p$. It is not difficult to see that 
		\begin{align*}
			H = \left \langle I + \begin{pmatrix} a & 0 \\ 0 & d \end{pmatrix} p, I + \begin{pmatrix} 0 & 1 \\ 0 & 0 \end{pmatrix} p \right \rangle, \text{ where } a,d \in \mathbb{Z}/p\mathbb{Z}
		\end{align*}
			or
		\begin{align*}
			H = \left \langle I + \begin{pmatrix} 1 & 0 \\ 0 & 0 \end{pmatrix} p, I + \begin{pmatrix} 0 & 0 \\ 0 & 1 \end{pmatrix} p, I + \begin{pmatrix} 0 & 1 \\ 0 & 0 \end{pmatrix} p \right \rangle,
		\end{align*}
	\end{proof}
\end{lemma}

\begin{proposition} \label{propBpConjTwo}
	Suppose that $p > 3$. Let $H \leq \GL_2(\mathbb{Z}/p^2\mathbb{Z})$ with $\varphi(H) \leq \left \langle \begin{pmatrix} 1 & 1 \\ 0 & 1 \end{pmatrix} \right \rangle$. $H$ is conjugate to some $H' \leq \GL_2(\mathbb{Z}/p^2\mathbb{Z})$ which is of the form $H' = \langle \tau, H' \cap \ker \varphi \rangle$, where one of the following holds:
	\begin{enumerate}
		\item $H' \cap \ker \varphi$ is $\ker \varphi$ and $\tau = \begin{pmatrix} 1 & 1 \\ 0 & 1 \end{pmatrix}$.
		\item $H' \cap \ker \varphi$ is one of 
			\begin{enumerate}
				\item $\left \langle I + \begin{pmatrix} 1 & 0 \\ \gamma & 1 \end{pmatrix} p, I + \begin{pmatrix} 1 & 0 \\ 0 & -1 \end{pmatrix} p, I + \begin{pmatrix} 0 & 1 \\ 0 & 0 \end{pmatrix} p \right \rangle$ where $\gamma \neq 0$
				\item $T$
			\end{enumerate}
			and $\tau$ is one of 
			\begin{enumerate}
				\item $\begin{pmatrix} 1 & 1 \\ 0 & 1 \end{pmatrix}$
				\item $\begin{pmatrix} 1 & 1 \\ 0 & 1 \end{pmatrix} + \begin{pmatrix} 1 & 0 \\ 0 & 1 \end{pmatrix} p$
			\end{enumerate}
		\item $H' \cap \ker \varphi$ does not contain any elements of the form $I + \begin{pmatrix} \alpha & \beta \\ \gamma & \delta \end{pmatrix} p$ where $\gamma \neq 0$ and $\tau$ is one of
			\begin{enumerate}
				\item $\begin{pmatrix} 1 & 1 \\ 0 & 1 \end{pmatrix}$
				\item $\begin{pmatrix} 1 & 1 \\ 0 & 1 \end{pmatrix} + \begin{pmatrix} 1 & 0 \\ 0 & 1 \end{pmatrix}p$
				\item	$\begin{pmatrix} 1 & 1 \\ 0 & 1 \end{pmatrix} + \begin{pmatrix} a & 0 \\ 1 & a \end{pmatrix} p$ for some $a \in \mathbb{Z}/p\mathbb{Z}$
				\item $\begin{pmatrix} 1 & 1 \\ 0 & 1 \end{pmatrix} + \begin{pmatrix} a & 0 \\ \epsilon & a \end{pmatrix} p$ for some $a \in \mathbb{Z}/p\mathbb{Z}$. 
			\end{enumerate}.
	\end{enumerate}
	\begin{proof}
		Let $H'$ be a conjugate of $H$ where $H' = \langle \tau, H' \cap \ker \varphi \rangle$ with $\tau = \begin{pmatrix} 1 & 1 \\ 0 & 1 \end{pmatrix} + \begin{pmatrix} a & 0 \\ c & a \end{pmatrix} p$ using Lemma \ref{lemmaBpConjOne}. If $H' \cap \ker \varphi = \ker \varphi$, then clearly $\begin{pmatrix} 1 & 1 \\ 0 & 1 \end{pmatrix} \in H'$, so assume that $H' \cap \ker \varphi \neq \ker \varphi$. \par
		Suppose that $H' \cap \ker \varphi = \left \langle I + \begin{pmatrix} \alpha & 0 \\ \gamma & \alpha \end{pmatrix} p, I + \begin{pmatrix} 1 & 0 \\ 0 & -1 \end{pmatrix}, I + \begin{pmatrix} 0 & 1 \\ 0 & 0 \end{pmatrix} p \right \rangle$ where $\gamma \neq 0$. Alternatively, we can choose either $\alpha = 1$ or $H' \cap \ker \varphi = T$. By using Lemmas \ref{lemmaTProdThree}, \ref{lemmaTProdTwo} and \ref{lemmaTProdOne} in that order, one sees that $H'$ has an element $\tau'$ of the form $\tau' =\begin{pmatrix} 1 & 1 \\ 0 & 1 \end{pmatrix} + \begin{pmatrix} a' & 0 \\ 0 & a' \end{pmatrix} p$. If $a' = 0$, then we are done. Otherwise,
		\begin{align*}
			\tau'^{\frac{1}{a'}} = \begin{pmatrix} 1 & \frac{1}{a'} \\ 0 & 1 \end{pmatrix} + \begin{pmatrix} 1 & 0 \\ 0 & 1 \end{pmatrix} p
		\end{align*}
		by Lemma \ref{lemmaTPower}. Replace $H'$ with its conjugate via $\begin{pmatrix} a' & 0 \\ 0 & 1 \end{pmatrix}$. By Lemma \ref{lemmaDiagConj}, $\begin{pmatrix} 1 & 1 \\ 0 & 1 \end{pmatrix} + \begin{pmatrix} 1 & 0 \\ 0 & 1 \end{pmatrix} p$ is an element of $H'$. Therefore, $H'$
		\begin{align*}
			\left \langle \begin{pmatrix} 1 & 1 \\ 0 & 1 \end{pmatrix} + \begin{pmatrix} 1 & 0 \\ 0 & 1 \end{pmatrix} p, H' \cap \ker \varphi \right \rangle
		\end{align*}
		and $H' \cap \ker \varphi$ is $T$ or $\left \langle I + \begin{pmatrix} 1 & 0 \\ \gamma' & 1 \end{pmatrix} p, I + \begin{pmatrix} 1 & 0 \\ 0 & -1 \end{pmatrix}, I + \begin{pmatrix} 0 & 1 \\ 0 & 0 \end{pmatrix} p \right \rangle $ for some $\gamma' \neq 0$. \par
		Now further assume that $H' \cap \ker \varphi$ is not of the form 
		\begin{align*}
		H' \cap \ker \varphi = \left \langle I + \begin{pmatrix} \alpha & 0 \\ \gamma & \alpha \end{pmatrix} p, I + \begin{pmatrix} 1 & 0 \\ 0 & -1 \end{pmatrix}, I + \begin{pmatrix} 0 & 1 \\ 0 & 0 \end{pmatrix} p \right \rangle
		\end{align*}
		where $\gamma \neq 0$. By Lemma \ref{lemmaKerBp}, $H' \cap \ker \varphi$ does not contain any elements of the form $I + \begin{pmatrix} \alpha & \beta \\ \gamma & \delta \end{pmatrix} p$ where $\gamma \neq 0$. If $c = 0$, then proceed just as in the last paragraph to replace $H'$ with a conjugate of the form $H' = \langle \tau, H' \cap \ker \varphi \rangle$ where $\tau$ is $\begin{pmatrix} 1 & 1 \\ 0 & 1 \end{pmatrix}$ or $\begin{pmatrix} 1 & 1 \\ 0 & 1 \end{pmatrix} + \begin{pmatrix} 1 & 0 \\ 0 & 1 \end{pmatrix} p$. Assume that $c \neq 0$. If $c$ is a square, then replace $H'$ with its conjugate via $\begin{pmatrix} \sqrt{c} & 0 \\ 0 & \frac{1}{\sqrt{c}} \end{pmatrix}$ so that $\begin{pmatrix} 1 & \sqrt{c} \\ 0 & 1 \end{pmatrix} + \begin{pmatrix} a & 0 \\ \sqrt{c} & a \end{pmatrix} p \in H'$ by Lemma \ref{lemmaDiagConj}. By Lemma \ref{lemmaTPower}, raising this element to the $\frac{1}{\sqrt{c}}$th power yields a matrix of the form
		\begin{align*}
			\begin{pmatrix} 1 & 1 \\ 0 & 1 \end{pmatrix} + \begin{pmatrix} a' & * \\ 1 & a' \end{pmatrix} p,
		\end{align*}
		and so $\begin{pmatrix} 1 & 1 \\ 0 & 1 \end{pmatrix} + \begin{pmatrix} a' & 0 \\ 1 & a' \end{pmatrix} p \in H'$. If $c$ is a nonsquare, then one can similarly replace $H'$ with a conjugate so that a matrix of the form $\begin{pmatrix} 1 & 1 \\ 0 & 1 \end{pmatrix} + \begin{pmatrix} a' & 0 \\ \epsilon & a' \end{pmatrix} p \in H'$ and so $H'$ is therefore of the desired form.
	\end{proof}
\end{proposition}

\begin{lemma} \label{lemmaBpParts}
	Let $H_1, H_2 \leq \GL_2(\mathbb{Z}/p^2\mathbb{Z})$ such that $\begin{pmatrix} 1 & 1 \\ 0 & 1 \end{pmatrix} \in \varphi(H_i)$. Suppose that $H_1$ and $H_2$ are locally conjugate and that $H_i$ is of the form $\langle \tau_i, H_i \cap \ker \varphi \rangle \rtimes C_{H_i}$ for some $\tau_i \in \varphi^{-1} \left( \begin{pmatrix} 1 & 1 \\ 0 & 1 \end{pmatrix} \right)$. Then, $\langle \tau_1, H_1 \cap \ker \varphi \rangle$ and $\langle \tau_2, H_2 \cap \ker \varphi \rangle$ are locally conjugate and $C_{H_1}$ and $C_{H_2}$ are equal or conjugate via $\begin{pmatrix} 0 & 1 \\ 1 & 0 \end{pmatrix}$.
	\begin{proof}
		Since $H_1$ and $H_2$ are locally conjugate, there is a bijection $f: H_1 \rightarrow H_2$ in which corresponding elements are conjugate. By Theorem \ref{theoremConjClass}, the only elements of $B(p)$ that are conjugate to $\begin{pmatrix} 1 & 1 \\ 0 & 1 \end{pmatrix}$ are of the form $\begin{pmatrix} 1 & n \\ 0 & 1 \end{pmatrix}$ where $n \neq 0$. $f$ must therefore map the elements of $\langle \tau_1, H_1 \cap \ker \varphi \rangle$ into $\langle \tau_2, H_2 \cap \ker \varphi \rangle$ and vice versa, and so these two groups are locally conjugate. \par
		Similarly, the only elements of $B(p)$ that are conjugate to $\begin{pmatrix} w & x \\ 0 & w \end{pmatrix} p$ where $x \neq 0$ are of the form $\begin{pmatrix} w & x' \\ 0 & w \end{pmatrix}$ where $x' \neq 0$. Furthermore, the only elements of $B(p)$ that are conjugate to $\begin{pmatrix} w & x \\ 0 & z \end{pmatrix} p$ where $w \neq z$ are of the form $\begin{pmatrix} w & x' \\ 0 & z \end{pmatrix}$ or $\begin{pmatrix} z & x' \\ 0 & w \end{pmatrix}$. Subgroups of $C_s(p)$ are generated by at most two elements, one of which can be chosen to be in $Z(p)$. It is then not difficult to see that $C_{H_1}$ and $C_{H_2}$ are equal or diagonal swaps.
	\end{proof}
\end{lemma}

\begin{lemma}\label{lemmaConjClassBp}
	The conjugacy class of an element $\begin{pmatrix} 1 & n \\ 0 & 1 \end{pmatrix} + \begin{pmatrix} a & b \\ c & d \end{pmatrix} p$ in $\GL_2(\mathbb{Z}/p^2\mathbb{Z})$, where $0 < n < p$, is determined completely by $a+d$ and $cn$.
	\begin{proof}
		The determinant of the matrix is $1 + (a+d-cn)p$ whereas the trace is $2+(a+d)p$. Apply Theorem \ref{theoremConjClass}. 
	\end{proof}
\end{lemma}

\begin{proposition} \label{propBorel}
	Suppose that $p > 3$. Let $H_1, H_2 \leq \GL_2(\mathbb{Z}/p\mathbb{Z}/p^2\mathbb{Z})$ be nontrivially locally conjugate with $\begin{pmatrix} 1 & 1 \\ 0 & 1 \end{pmatrix} \in \varphi(H_i)$. $H_1$ and $H_2$ are conjugate to
	\begin{align*}
		\left \langle \tau, k, D \right \rangle, \left \langle \tau, k, D' \right \rangle
	\end{align*}
	in some order, where $\tau$ is one of 
	\begin{enumerate}
		\item $\begin{pmatrix} 1 & 1 \\ 0 & 1 \end{pmatrix}$
		\item $\begin{pmatrix} 1 & 1 \\ 0 & 1 \end{pmatrix} + \begin{pmatrix} 1 & 0 \\ 0 & 1 \end{pmatrix} p$
		\item $\begin{pmatrix} 1 & 1 \\ 0 & 1 \end{pmatrix} + \begin{pmatrix} a & 0 \\ 1 & a \end{pmatrix} p$ for some $a \in \mathbb{Z}/p\mathbb{Z}$
		\item $\begin{pmatrix} 1 & 1 \\ 0 & 1 \end{pmatrix} + \begin{pmatrix} a & 0 \\ \epsilon & a \end{pmatrix} p$ for some $a \in \mathbb{Z}/p\mathbb{Z}$,
	\end{enumerate}
	$k$ is one of 
	\begin{enumerate}
		\item $I$
		\item $I + \begin{pmatrix} 0 & 0 \\ 1 & 0 \end{pmatrix} p$
	\end{enumerate}
		$D$ is a subgroup of $C_s(p)$, $D$ and $D'$ are diagonal swaps but $D \neq D'$. 
	\begin{proof}
		Replace $H_i$ with a conjugate in the form specified in Lemma \ref{lemmaBpConjOne}. In particular, $H_i = \langle \tau_i, H_i \cap \ker \varphi \rangle \rtimes C_{H_i}$ for some $\tau_i$ of the form $\tau_i = \begin{pmatrix} 1 & 1 \\ 0 & 1 \end{pmatrix} + \begin{pmatrix} a_i & 0 \\ c_i & a_i \end{pmatrix} p$. By Lemma \ref{lemmaBpParts}, $C_{H_1}$ and $C_{H_2}$ are equal or diagonal swaps, i.e. $\varphi(H_1)$ and $\varphi(H_2)$ are equal or diagonal swaps. First consider the case where $\varphi(H_i) \cap C_s(p) = \langle I \rangle$, i.e. $\varphi(H_i) = \langle t \rangle$. Replace $H_1$ and $H_2$ with conjugates listed in Proposition \ref{propBpConjTwo}. In particular, $H_i = \langle \tau_i, H_i \cap \ker \varphi \rangle$ for a $\tau_i$ listed in Proposition \ref{propBpConjTwo}. By Proposition \ref{propNecessaryOne}, $H_1 \cap \ker \varphi$ and $H_2 \cap \varphi$ are locally conjugate and so $\dim (H_1 \cap \ker \varphi) = \dim (H_2 \cap \ker \varphi)$ by Lemma \ref{lemmaLocConjBij}. If $\dim(H_i \cap \ker \varphi) = 4$, then $H_i \cap \ker \varphi = \ker \varphi$, in which case $H_1$ and $H_2$ are equal. \par
		Assume that $\dim(H_i \cap \ker \varphi) \leq 3$. If $\dim(H_i \cap \ker \varphi) = 3$, then $H_i \cap \ker \varphi$ is one of 
		\begin{enumerate}
			\item $\left \langle I + \begin{pmatrix} 1 & 0 \\ \gamma & 1 \end{pmatrix} p, I + \begin{pmatrix} 1 & 0 \\ 0 & -1 \end{pmatrix} p, I + \begin{pmatrix} 0 & 1 \\ 0 & 0 \end{pmatrix} p \right \rangle$ where $\gamma \neq 0$
			\item $T$ 
			\item $\left \langle I + \begin{pmatrix} 1 & 0 \\ 0 & 0 \end{pmatrix} p, I + \begin{pmatrix} 0 & 1 \\ 0 & 0 \end{pmatrix} p, I + \begin{pmatrix} 0 & 0 \\ 0 & 1 \end{pmatrix} p \right \rangle$.
		\end{enumerate}
		No two among these three are locally conjugate by Proposition \ref{propLocConjKerEquiv}. Suppose that $H_i \cap \ker \varphi = \left \langle I + \begin{pmatrix} 1 & 0 \\ \gamma_i & 1 \end{pmatrix} p, I + \begin{pmatrix} 1 & 0 \\ 0 & -1 \end{pmatrix} p, I + \begin{pmatrix} 0 & 1 \\ 0 & 0 \end{pmatrix} p \right \rangle$ where $\gamma_i \neq 0$. Further suppose, for contradiction, that $\tau_1 \neq \tau_2$. By Proposition \ref{propBpConjTwo}, one can set $\tau_1 = \begin{pmatrix} 1 & 1 \\ 0 & 1 \end{pmatrix}$ and $\tau_2 = \begin{pmatrix} 1 & 1 \\ 0 & 1 \end{pmatrix} + \begin{pmatrix} 1 & 0 \\ 0 & 1 \end{pmatrix} p$ without loss of generality. Note that $H_i = \langle \tau_1 \rangle \rtimes (H_i \cap \ker \varphi)$. $\tau_2$ must be conjugate to some $\tau_1^n k_1$ where $k_1 \in H_1 \cap \ker \varphi$, but this is cannot happen by Lemma \ref{lemmaConjClassBp}. Hence, $\tau_1 = \tau_2$. If $\tau_i = \begin{pmatrix} 1 & 1 \\ 0 & 1 \end{pmatrix}$, then $H_1$ is in fact conjugate to $H_2$ via $\begin{pmatrix} \gamma_1 & 0 \\ 0 & \gamma_2 \end{pmatrix}$. Now suppose that $\tau_i = \begin{pmatrix} 1 & 1 \\ 0 & 1 \end{pmatrix} + \begin{pmatrix} 1 & 0 \\ 0 & 1 \end{pmatrix} p$. The elements of $H_i$ are of the form $\tau_i^{n_i} k_i$ for some $k_i \in H_i \cap \ker \varphi$. Moreover, it is not difficult using Lemma \ref{lemmaConjClassBp} to see that multiplying $I + \begin{pmatrix} 1 & 0 \\ 0 & -1 \end{pmatrix} p$ or $I + \begin{pmatrix} 0 & 1 \\ 0 & 0 \end{pmatrix} p$ to $k_i$ preserves the conjugacy class of $\tau_i^{n_i} k_i$. Thus, local conjugacy of $H_i$ is determined by the conjugacy classes of the elements $\tau_i^{n_i} \left(I  +\begin{pmatrix} 1 & 0 \\ \gamma_i & 1 \end{pmatrix} p \right)^{m_i}$. Such an element is of the form
		\begin{align*}
			\tau_i^{n_i} \left(I + \begin{pmatrix} 1 & 0 \\ \gamma_i & 1 \end{pmatrix} p \right)^{m_i} &= \left( \begin{pmatrix} 1 & n_i \\ 0 & 1 \end{pmatrix} + \begin{pmatrix} n_i & 0 \\ 0 & n_i \end{pmatrix} p \right) \left( I + \begin{pmatrix} m_i & 0 \\ m_i \gamma_i & m_i \end{pmatrix} p \right)  \\
			&=\begin{pmatrix} 1 & n_i \\ 0 & 1 \end{pmatrix} + \left( \begin{pmatrix} n_i & 0 \\ 0 & n_i \end{pmatrix} + \begin{pmatrix} m_i + m_i n_i \gamma_i & m_in_i \\ m_i \gamma_i & m_i \end{pmatrix} \right) p \\
			&= \begin{pmatrix} 1 & n_i \\ 0 & 1 \end{pmatrix} + \begin{pmatrix} m_i + n_i + m_i n_i \gamma_i & m_i n_i \\ m_i \gamma_i & m_i + n_i \end{pmatrix} p.
		\end{align*}
		In the case that $n_i \neq 0 \pmod{p}$, Lemma \ref{lemmaConjClassBp} asserts that its conjugacy class is determined by the values of $2m_i + 2n_i + m_in_i\gamma_i$ and $m_in_i\gamma_i$. Alternatively, the conjugacy class is determined by the values of $m_i + n_i$ and $m_in_i \gamma_i$. Parametrize $m_1$ as $m_1 = rn_1$. Note that $m_1 + n_1 = (r+1)n_1$ and $m_1n_1 \gamma_1 = rn_1^2 \gamma_1$, and so there must be some $m_2$ and $n_2$ such that $m_2+n_2 = (r+1)n_1$ and $m_2n_2 = rn_1^2 \frac{\gamma_1}{\gamma_2}$, i.e. there is a solution to the quadratic equation $x^2 - (r+1) n_1 x+ rn_1^2 \frac{\gamma_1}{\gamma_2} = 0$. This happens only when the discriminant, which is $(r+1)^2n_1^2 - 4rn_1^2 \frac{\gamma_1}{\gamma_2}$, is a square in $\mathbb{Z}/p\mathbb{Z}$. When $n_1 \neq 0$, $(r+1)^2 - 4r \frac{\gamma_1}{\gamma_2}$ must be a square. Completing the square shows that 
		\begin{align*}
		(r+1)^2 - 4r \frac{\gamma_1}{\gamma_2} = \left( r + \left( 1 - 2 \frac{\gamma_1}{\gamma_2} \right) \right)^2 + 1 - \left( 1 - 2 \frac{\gamma_1}{\gamma_2} \right)^2.
		\end{align*}
		This value needs to be a square for all $r \in \mathbb{Z}/p\mathbb{Z}$. Therefore, for all squares $s \in \mathbb{Z}/p\mathbb{Z}$, $s + 1 - \left( 1 - 2 \frac{\gamma_1}{\gamma_2} \right)^2$ is also a square. Let $c = 1 - \left( 1 - 2 \frac{\gamma_1}{\gamma_2} \right)^2$ and suppose that $c \neq 0$. The case where $s = 0$ shows that $c$ must be a square, say $c = t^2$. In this case, $s/t^2 + 1$ is a square for all squares $s \in \mathbb{Z}/p\mathbb{Z}$. However, $s/t^2$ spans all squares in $\mathbb{Z}/p\mathbb{Z}$, which is a contradiction. Hence, $1 - \left( 1 - 2 \frac{\gamma_1}{\gamma_2} \right)^2 = 0$ and so $\gamma_1 = \gamma_2$. $H_1$ and $H_2$ are therefore equal. \par
		
		If $H_i \cap \ker \varphi = T$, then one can show that $\tau_1 = \tau_2$ similarly as above, in which case $H_1 = H_2$. \par
		For the remaining possibilities of $H_i \cap \ker \varphi$, $H_i \cap \ker \varphi$ has no element of the form $I + \begin{pmatrix} \alpha & \beta \\ \gamma & \delta \end{pmatrix} p$ with $\gamma \neq 0$. It is then not difficult to see that if $\tau_1$ is of the form $\tau_1 = \begin{pmatrix} 1 & 1 \\ 0 & 1 \end{pmatrix} + \begin{pmatrix} a_1 & 0 \\ c & a_1 \end{pmatrix} p$ where $c = 0,1$ or $\epsilon$, then $\tau_2$ is of the form $\tau_2 = \begin{pmatrix} 1 & 1 \\ 0 & 1 \end{pmatrix} + \begin{pmatrix} a_2 & 0 \\ c & a_2 \end{pmatrix}$ using Lemma \ref{lemmaConjClassBp}. \par
		If $H_i \cap \ker \varphi = \left \langle I + \begin{pmatrix} 1 & 0 \\ 0 & 0 \end{pmatrix} p, I + \begin{pmatrix} 0 & 1 \\ 0 & 0 \end{pmatrix} p, I + \begin{pmatrix} 0 & 0 \\ 0 & 1 \end{pmatrix} p \right \rangle$, then $\tau_i$ can be chosen to be of the form $\begin{pmatrix} 1 & 1 \\ 0 & 1 \end{pmatrix} + \begin{pmatrix} 0 & 0 \\ c & 0 \end{pmatrix} p$. Thus, $H_1 = H_2$. \par
		Suppose that $H_1 \cap \ker \varphi = \left \langle I + \begin{pmatrix} \alpha & 0 \\ 0 & \delta \end{pmatrix} p, I + \begin{pmatrix} 0 & 1 \\ 0 & 0 \end{pmatrix} p \right \rangle$ for some $\alpha,\delta \in \mathbb{Z}/p\mathbb{Z}$ which are not both $0$. Since $H_1 \cap \ker \varphi$ and $H_2 \cap \ker \varphi$ are locally conjugate, $H_2 \cap \ker \varphi = \left \langle I + \begin{pmatrix} \alpha & 0 \\ 0 & \delta \end{pmatrix} p, I + \begin{pmatrix} 0 & 1 \\ 0 & 0 \end{pmatrix} p \right \rangle$ or $\left \langle I + \begin{pmatrix} \delta & 0 \\ 0 & \alpha \end{pmatrix} p, I + \begin{pmatrix} 0 & 1 \\ 0 & 0 \end{pmatrix} p \right \rangle$ by Proposition \ref{propNTLCKer}. Further suppose that $\alpha \neq -\delta$. If $c = 1$, then replace $H_1$ by its conjugate via $\begin{pmatrix} 1 & n \\ 0 & 1 \end{pmatrix}$, where $n = a_1 \frac{\alpha - \delta}{\alpha + \delta}$. Doing so preserves $H_1 \cap \ker \varphi$ and takes $\tau_1$ to
		\begin{align*}
			\begin{pmatrix} 1 & n \\ 0 & 1 \end{pmatrix} \left( \begin{pmatrix} 1 & 1 \\ 0 & 1 \end{pmatrix} + \begin{pmatrix} a_1 & 0 \\ 1 & a_1 \end{pmatrix} p \right) \begin{pmatrix} 1 & n \\ 0 & 1 \end{pmatrix}^{-1} = \begin{pmatrix} 1 & 1 \\ 0 & 1 \end{pmatrix} + \begin{pmatrix} a_1+n & -n^2 \\ 1 & a_1-n \end{pmatrix} p. 
		\end{align*}
		Note that $\frac{a_1+n}{a_1-n} = \frac{\alpha}{\delta}$, and so $H_1$ has an element $\tau'_1$ of the form $\tau'_1 = \begin{pmatrix} 1 & 1 \\ 0 & 1 \end{pmatrix} + \begin{pmatrix} 0 & 0 \\ 1 & 0 \end{pmatrix} p$. Thus, $H_1 = \langle \tau'_1 \rangle \rtimes (H_1 \cap \ker \varphi)$. Similarly, if $\tau_1$ starts out as $\begin{pmatrix} 1 & 1 \\ 0 & 1 \end{pmatrix} + \begin{pmatrix} a_1 & 0 \\ \epsilon & a_1 \end{pmatrix} p$, then $H_1$ can be replaced by a conjugate so that $H_1 \cap \ker \varphi$ is preserved and $H_1 = \langle \tau'_1 \rangle \rtimes (H_1 \cap \ker \varphi)$, where $\tau'_1 = \begin{pmatrix} 1 & 1 \\ 0 & 1 \end{pmatrix} + \begin{pmatrix} 0 & 0 \\ \epsilon & 0 \end{pmatrix} p$. If $\tau_1 = \begin{pmatrix} 1 & 1 \\ 0 & 1 \end{pmatrix} + \begin{pmatrix} 1 & 0 \\ 0 & 1 \end{pmatrix}$, then replace $H_1$ with its conjugate via $I + \begin{pmatrix} 0 & 0 \\ \gamma & 0 \end{pmatrix} p$ where $\gamma = \frac{\delta - \alpha}{\delta + \alpha}$. Doing so preserves $H_1 \cap \ker \varphi$, while making $H_1 = \left \langle \begin{pmatrix} 1 & 1 \\ 0 & 1 \end{pmatrix} \right \rangle \rtimes (H_1 \cap \ker \varphi)$. Similarly replace $H_2$ with a conjugate. If $H_1 \cap \ker \varphi = H_2 \cap \ker \varphi$, then $H_1$ and $H_2$ are conjugate. Otherwise, it is not difficult to see that $H_1$ and $H_2$ are nontrivially locally conjugate. \par
		Now suppose that $\alpha = -\delta$. Local conjugacy of $H_i$ is determined by the conjugacy classes of the elements of $\langle \tau_i \rangle$ due to Lemma \ref{lemmaConjClassBp}. Therefore, if $\tau_1 = \begin{pmatrix} 1 & 1 \\ 0 & 1 \end{pmatrix}$, then $\tau_2 = \begin{pmatrix} 1 & 1 \\ 0 & 1 \end{pmatrix}$ and if $\tau_1 = \begin{pmatrix} 1 & 1 \\ 0 & 1 \end{pmatrix} + \begin{pmatrix} 1 & 0 \\ 0 & 1 \end{pmatrix} p$, then $\tau_2 = \begin{pmatrix} 1 & 1 \\ 0 & 1 \end{pmatrix} + \begin{pmatrix} 1 & 0 \\ 0 & 1 \end{pmatrix} p$. Lemma \ref{lemmaConjClassBp} further shows that if $\tau_1 = \begin{pmatrix} 1 & 1 \\ 0 & 1 \end{pmatrix} + \begin{pmatrix} a_1 & 0 \\ 1 & a_1 \end{pmatrix} p$ or $\begin{pmatrix} 1 & 1 \\ 0 & 1 \end{pmatrix} + \begin{pmatrix} a_1 & 0 \\ \epsilon & a_1 \end{pmatrix} p$, then $\tau_1$ is conjugate to $\tau_2$ or to $\tau_2^{-1}$. Any such conjugation must be via an element of $\varphi^{-1}(B(p))$, which preserves $H_i \cap \ker \varphi$. Again, if $H_1 \cap \ker \varphi = H_2 \cap \ker \varphi$, then $H_1$ and $H_2$ are conjugate. Otherwise, it is not difficult to see that $H_1$ and $H_2$ are nontrivially locally conjugate. \par
		If $H_i \cap \ker \varphi = \left \langle I + \begin{pmatrix} 0 & 1 \\ 0 & 0 \end{pmatrix} p \right \rangle$, then $H_i = \langle \tau_i \rangle$, and so $H_1$ and $H_2$ are conjugate. This concludes the case where $\varphi(H_i) \cap C_s(p) = \langle I \rangle$; i.e. we have showed that $H_1$ and $H_2$ are conjugate to groups of the desired form. The case where $\varphi(H_i) \cap C_s(p) \leq Z(p^2)$ is similar. \par
		Now assume that $\varphi(H_i) \cap C_s(p) \not\leq Z(p^2)$. Recall that $H_i = \langle \tau_i, H_i \cap \ker \varphi \rangle \rtimes C_{H_i}$. If $H_i \cap \ker \varphi = \ker \varphi$, then $H_1$ and $H_2$ are of the desired form; in particular, let $\tau = \begin{pmatrix} 1 & 1 \\ 0 & 1 \end{pmatrix}$, $k = I + \begin{pmatrix} 0 & 0 \\ 1 & 0 \end{pmatrix} p$ and $D = \left \langle C_{H_i}, I + \begin{pmatrix} 1 & 0 \\ 0 & 1 \end{pmatrix} p \right \rangle$. \par
		Suppose that $H_i \cap \ker \varphi = T$. If $c_i \neq 0$, then observe that $H_i$ has $\begin{pmatrix} 1 & 1 \\ 0 & 1 \end{pmatrix} + \begin{pmatrix} a_i-c_i & 0 \\ 0 & a_i \end{pmatrix} p$ by Lemma \ref{lemmaTProdThree}. One can multiply this matrix by powers of $I + \begin{pmatrix} 1 & 0 \\ 0 & -1 \end{pmatrix} p$ and $I + \begin{pmatrix} 0 & 1 \\ 0 & 0 \end{pmatrix} p$ to see that $H_i$ has a matrix of the form $\begin{pmatrix} 1 & 1 \\ 0 & 1 \end{pmatrix} + \begin{pmatrix} a'_i & 0 \\ 0 & a'_i \end{pmatrix} p$. Thus, whether or not $c_i = 0$, $H_i$ has a matrix of this form. By Lemmas \ref{lemmaBpParts} and \ref{lemmaConjClassBp}, $a_1$ and $a_2$ must be both $0$ or both nonzero. Suppose that $a'_i \neq 0$. The $\frac{1}{a'_i}$th power of this matrix has the form $\begin{pmatrix} 1 & \frac{1}{a'_i} \\ 0 & 1 \end{pmatrix} + \begin{pmatrix} 1 & 0 \\ 0 & 1 \end{pmatrix} p$. Replacing $H_i$ by its conjugate via $\begin{pmatrix} a'_i & 0 \\ 0 & 1 \end{pmatrix}$ preserves $H_i \cap \ker \varphi$ and $C_{H_i}$ and makes $H_i$ have $\begin{pmatrix} 1 & 1 \\ 0 & 1 \end{pmatrix} + \begin{pmatrix} 1 & 0 \\ 0 & 1 \end{pmatrix} p$ as an element. Therefore, $H_1$ and $H_2$ are of desired form; in particular, let $\tau = \begin{pmatrix} 1 & 1 \\ 0 & 1 \end{pmatrix}$ or $\begin{pmatrix} 1 & 1 \\ 0 & 1 \end{pmatrix} + \begin{pmatrix} 1 & 0 \\ 0 & 1 \end{pmatrix} p$, $k = I$ and $D = C_{H_i}$. \par
		
		Since $\varphi(H_i) \cap C_s(p)$ contains some element of the form $\begin{pmatrix} w & 0 \\ 0 & z \end{pmatrix}$ where $w \neq z$, $H_i \cap \ker \varphi$ is not of the form
		\begin{align*}
		\left \langle I + \begin{pmatrix} 1 & 0 \\ \gamma & 1 \end{pmatrix} p, I + \begin{pmatrix} 1 & 0 \\ 0 & -1 \end{pmatrix} p, I + \begin{pmatrix} 0 & 1 \\ 0 & 0 \end{pmatrix} p \right \rangle
		\end{align*}
		for $\gamma \neq 0$ by Corollary \ref{corDiagMember}. \par
		
		Suppose that $H_i \cap \ker \varphi$ has no element of the form $I + \begin{pmatrix} \alpha & \beta \\ \gamma & \delta \end{pmatrix} p$ with $\gamma \neq 0$. By Lemmas \ref{lemmaBpParts} and \ref{lemmaConjClassBp}, $c_1$ and $c_2$ must be both $0$, both nonzero squares, or both nonsquares. Suppose that $c_i = 0$. If $a_i \neq 0$, then replacing $H_i$ with its conjugate via $\begin{pmatrix} a_i & 0 \\ 0 & 1 \end{pmatrix}$ makes $H_i$ have $\begin{pmatrix} 1 & 1 \\ 0 & 1 \end{pmatrix} + \begin{pmatrix} 1 & 0 \\ 0 & 1 \end{pmatrix} p$ while preserving $H_i \cap \ker \varphi$ and $C_{H_i}$. If $a_i = 0$, then $\begin{pmatrix} 1 & 1 \\ 0 & 1 \end{pmatrix} \in H_i$. Suppose that $c_1$ and $c_2$ are both nonzero squares. In this case, $H_i$ has 
		\begin{align*}
			\left( \begin{pmatrix} 1 & 1 \\ 0 & 1 \end{pmatrix} + \begin{pmatrix} a_i & 0 \\ c_i & a_i \end{pmatrix} p \right)^{\frac{1}{\sqrt{c_i}}} = \begin{pmatrix} 1 & \frac{1}{\sqrt{c_i}} \\ 0 & 1 \end{pmatrix} + \begin{pmatrix} a'_i & * \\ \sqrt{c_i} & a'_i \end{pmatrix} p
		\end{align*}
		for some $a'_i \in \mathbb{Z}/p\mathbb{Z}$ and so $H_i$ also has $\begin{pmatrix} 1 & \frac{1}{\sqrt{c_i}} \\ 0 & 1 \end{pmatrix} + \begin{pmatrix} a'_i & 0 \\ \sqrt{c_i} & a'_i \end{pmatrix} p$. Replacing $H_i$ with its conjugate via $\begin{pmatrix} \sqrt{c_i} & 0 \\ 0 & 1 \end{pmatrix}$ preserves $H_i \cap \ker \varphi$ and $C_{H_i}$ and makes $H_i$ have $\begin{pmatrix} 1 & 1 \\ 0 & 1 \end{pmatrix} + \begin{pmatrix} a'_i & 0 \\ 1 & a'_i \end{pmatrix} p$. If $c_i$ is a nonsquare, then one can similarly replace $H'_i$ with a conjugate such that $H'_i \cap \ker \varphi$ and $C_{H'_i}$ are preserved and $H_i$ has an element of the form $\begin{pmatrix} 1 & 1 \\ 0 & 1 \end{pmatrix} + \begin{pmatrix} a'_i & 0 \\ \epsilon & a'_i \end{pmatrix}p$. From here, assume that $\tau_i$ has one of the following: $\begin{pmatrix} 1 & 0 \\ 0 & 1 \end{pmatrix}, \begin{pmatrix} 1 & 0 \\ 0 & 1 \end{pmatrix} + \begin{pmatrix} 1 & 0 \\ 0 & 1 \end{pmatrix} p, I + \begin{pmatrix} a_i' & 0 \\ 1 & a_i' \end{pmatrix} p, \begin{pmatrix} 1 & 1 \\ 0 & 1 \end{pmatrix} + \begin{pmatrix} a_i' & 0 \\ \epsilon & a_i' \end{pmatrix} p$.
		
		Suppose that $H_i \cap \ker \varphi = \left \langle I + \begin{pmatrix} 1 & 0 \\ 0 & 0 \end{pmatrix} p, I + \begin{pmatrix} 0 & 1 \\ 0 & 0 \end{pmatrix} p, I + \begin{pmatrix} 0 & 0 \\ 0 & 1 \end{pmatrix} p \right \rangle$. $H_i$ has $\begin{pmatrix} 1 & 1 \\ 0 & 1 \end{pmatrix} + \begin{pmatrix} 0 & 0 \\ c_i & 0 \end{pmatrix} p$. Therefore, $H_1$ and $H_2$ are of the desired form; in particular, let $\tau = \begin{pmatrix} 1 & 1 \\ 0 & 1 \end{pmatrix}$, $\begin{pmatrix} 1 & 1 \\ 0 & 1 \end{pmatrix} + \begin{pmatrix} 0 & 0 \\ 1 & 0 \end{pmatrix}p$ or $\begin{pmatrix} 1 & 1 \\ 0 & 1 \end{pmatrix} + \begin{pmatrix} 0 & 0 \\ \epsilon & 0 \end{pmatrix}p$, $k = I$ and $D = \left \langle C_{H_i}, I + \begin{pmatrix} 1 & 0 \\ 0 & 0 \end{pmatrix} p, I + \begin{pmatrix} 0 & 0 \\ 0 & 1 \end{pmatrix} p \right \rangle$. \par
		
		Suppose that $H_1 \cap \ker \varphi$ is of the form $H_1 \cap \ker \varphi = \left \langle I + \begin{pmatrix} \alpha & 0 \\ 0 & \delta \end{pmatrix} p, I + \begin{pmatrix} 0 & 1 \\ 0 & 0 \end{pmatrix} p \right \rangle$. Express $\tau_i$ in the form $\begin{pmatrix} 1 & 1 \\ 0 & 1 \end{pmatrix} + \begin{pmatrix} a_i & 0 \\ c_i & a_i \end{pmatrix} p$ where $c_i = 0, 1$ or $\epsilon$ and $a_i = 0$ or $1$ if $c_i = 0$. Since $\varphi(H_i) \cap C_s(p) \not\leq Z(p^2)$, $\varphi(H_i)$ has an element of the form $\begin{pmatrix} w & 0 \\ 0 & z \end{pmatrix}$ where $w \neq z$. Just as in Lemma \ref{lemmaConvolutedBp}, compute
		\begin{align*}
			\begin{pmatrix} w & 0 \\ 0 & z \end{pmatrix} \tau \begin{pmatrix} w & 0 \\ 0 & z \end{pmatrix}^{-1} \tau^{-\frac{w}{z}} = I + \begin{pmatrix} * & * \\ c_i \left( \frac{z}{w} - \frac{w}{z} \right) & * \end{pmatrix} p.
		\end{align*}
		By the assumption on $H_1 \cap \ker \varphi$, $c_i = 0$ or every element of $\varphi(H_i) \cap C_s(p)$ is of the form $\begin{pmatrix} w & 0 \\ 0 & \pm w \end{pmatrix}$. Suppose that $c = 0$. If $\tau_i = \begin{pmatrix} 1 & 1 \\ 0 & 1 \end{pmatrix} + \begin{pmatrix} 1 & 0 \\ 0 & 1 \end{pmatrix} p$ for $i = 1$ or $2$, then another computation in Lemma \ref{lemmaConvolutedBp} shows that $I + \begin{pmatrix} 1 & 0 \\ 0 & 1 \end{pmatrix} p \in H_i$. In this case, $\begin{pmatrix} 1 & 1 \\ 0 & 1 \end{pmatrix} \in H_i$ and $H_i \cap \ker \varphi = \left \langle I + \begin{pmatrix} 1 & 0 \\ 0 & 1 \end{pmatrix} p, I + \begin{pmatrix} 0 & 1 \\ 0 & 0 \end{pmatrix} p \right \rangle$ for both $i = 1,2$. Whether or not $\tau_i = \begin{pmatrix} 1 & 1 \\ 0 & 1 \end{pmatrix} + \begin{pmatrix} 1 & 0 \\ 0 & 1 \end{pmatrix} p$, $H_1$ and $H_2$ are of the desired form; in particular, let $\tau = \begin{pmatrix} 1 & 1 \\ 0 & 1 \end{pmatrix}$, $k = I$ and $D = \left \langle H_i \cap \ker \varphi \cap C_s(p^2), C_{H_i} \right \rangle$. \par
	Suppose that $c_i \neq 0$. Yet another computation in Lemma \ref{lemmaConvolutedBp} shows that $I + \begin{pmatrix} 1 & 0 \\ 0 & 1 \end{pmatrix} p \in H_i$ given that $2a_i \neq c_i$. If $2a_i \neq c_i$ for $i = 1$ or $i = 2$, then $H_i \cap \ker \varphi = \left \langle I + \begin{pmatrix} 1 & 0 \\ 0 & 1 \end{pmatrix} p, I + \begin{pmatrix} 0 & 1 \\ 0 & 0 \end{pmatrix} p \right \rangle$ for both $i = 1$ and $2$, and so $H_1$ and $H_2$ are of the desired form; in particular, let $\tau = \begin{pmatrix} 1 & 1 \\ 0  & 1 \end{pmatrix} + \begin{pmatrix} 0 & 0 \\ c_i & 0 \end{pmatrix} p$, $k = I$ and $D = \left \langle H_i \cap \ker \varphi \cap C_s(p^2), C_{H_i} \right \rangle$. If $2a_i = c_i$ for both $i = 1$ and $2$, then $H_1$ and $H_2$ are of the desired form as well; in particular, let $\tau = \begin{pmatrix} 1 & 1 \\ 0 & 1 \end{pmatrix} + \begin{pmatrix} \frac{c_i}{2} & 0 \\ c_i & \frac{c_i}{2} \end{pmatrix} p$, $k = I$ and $D = \left \langle H_i \cap \ker \varphi \cap C_s(p^2), C_{H_i} \right \rangle$. 
	\end{proof}
\end{proposition}

For $p = 3$, there are $40$ pairs of nontrivially locally conjugate subgroups of $\GL_2(\mathbb{Z}/p^2\mathbb{Z})$ up to conjugation. They can be expressed similarly as the pairs of nontrivially locally conjugate subgroups of $\GL_2(\mathbb{Z}/p^2\mathbb{Z})$ for $p > 3$ as described in Proposition \ref{propBorel}. In particular, the pairs are expressible in the form
\begin{align*}
	 \langle \tau, k, D \rangle, \langle \tau, k, D' \rangle,
\end{align*}
where $\tau$ satisfies $\varphi(\tau) = \begin{pmatrix} 1 & 1 \\ 0 & 1 \end{pmatrix}$, $k \leq \ker \varphi$ and $D$ and $D'$ are unequal diagonal swaps. The main difference between the cases where $p = 3$ and $p > 3$ is in Lemma \ref{lemmaTMemberPG3}: whereas $I + \begin{pmatrix} 0 & 1 \\ 0 & 0 \end{pmatrix} p \in H$ if $p > 3$ and $H$ is a subgroup of $\GL_2(\mathbb{Z}/p^2\mathbb{Z})$ with $\begin{pmatrix} 1 & 1 \\ 0 & 1 \end{pmatrix} \in \varphi(H)$, then $I + \begin{pmatrix} 0 & 1 \\ 0 & 0 \end{pmatrix} p \in H$, but if $p = 3$, then this is not necessarily true.

\section{The $\SL_2(\mathbb{Z}/p\mathbb{Z})$ case} \label{sectionSL}
This section categorizes the subgroups, up to conjugation, of $\GL_2(\mathbb{Z}/p^2\mathbb{Z})$ whose images via $\varphi$ contain $\SL_2(\mathbb{Z}/p\mathbb{Z})$. 

\begin{lemma}\label{lemmaKerSL}
 Suppose that $p > 3$. Let $H \leq \GL_2(\mathbb{Z}/p^2\mathbb{Z})$. If $\SL_2(\mathbb{Z}/p\mathbb{Z}) \leq \varphi(H)$, then $T \leq H$.
	\begin{proof}
		A computation similar to that of Lemma \ref{lemmaTMemberPG3} shows that $I + \begin{pmatrix} 0 & 0 \\ 1 & 0 \end{pmatrix} p \in H$. By Lemma \ref{lemmaTMember}, $T \leq H$. 
	\end{proof}
\end{lemma}

\begin{lemma} \label{lemmaSubImSL2}
	Suppose that $p > 3$. If $H \leq \GL_2(\mathbb{Z}/p^2\mathbb{Z})$ with $\varphi(H) = \SL_2(\mathbb{Z}/p\mathbb{Z})$, then $H = \SL_2(\mathbb{Z}/p^2\mathbb{Z})$ or $\varphi^{-1}(\SL_2(\mathbb{Z}/p\mathbb{Z}))$. 
	\begin{proof}
		By Lemma \ref{lemmaKerSL}, $T \leq H$. Therefore, $H$ has an element $\tau$ of the form $\tau = \begin{pmatrix} 1 & 1 \\ 0 & 1 \end{pmatrix} + \begin{pmatrix} a & 0 \\ 0 & a \end{pmatrix} p$. Since $p > 3$, there is some $w \in \mathbb{Z}/p\mathbb{Z}$ such that $w^2 \neq 1$. In particular, $H$ has an element of the form $\begin{pmatrix} w & 0 \\ 0 & \frac{1}{w} \end{pmatrix} + \begin{pmatrix} 0 & b \\ c & 0 \end{pmatrix} p$ by Corollary \ref{corDiagMember}. Since $I + \begin{pmatrix} 0 & 1 \\ 0 & 0 \end{pmatrix}p , I + \begin{pmatrix} 0 & 0 \\ 1 & 0 \end{pmatrix} p \in H$, $H$ has $\begin{pmatrix} w & 0 \\ 0 & \frac{1}{w} \end{pmatrix}$ as well. If $a \neq 0$, then a computation in Lemma \ref{lemmaConvolutedBp} shows that $I + \begin{pmatrix} 1 & 0 \\ 0 & 1 \end{pmatrix} p$, and so $\ker \varphi \leq H$. In this case, $H = \varphi^{-1}(\SL_2(\mathbb{Z}/p\mathbb{Z}))$. Similarly, $H$ has an element of the form $\begin{pmatrix} 1 & 1 \\ 0 & 1 \end{pmatrix} + \begin{pmatrix} a' & 0 \\ 0 & a' \end{pmatrix} p$. If $a' \neq 0$, then $\ker \varphi \leq H$. If $a = a' = 0$, then $\begin{pmatrix} 1 & 1 \\ 0 & 1 \end{pmatrix}, \begin{pmatrix} 1 & 0 \\ 1 & 1 \end{pmatrix} \in H$ and so $\SL_2(\mathbb{Z}/p^2\mathbb{Z}) \leq H$. \par
		It now suffices to show that if $\SL_2(\mathbb{Z}/p^2\mathbb{Z}) \neq H$, then $H = \varphi^{-1}(\SL_2(\mathbb{Z}/p\mathbb{Z}))$. Choose some $h \in H$ that is not in $\SL_2(\mathbb{Z}/p^2\mathbb{Z})$. Since $\varphi(H) = \SL_2(\mathbb{Z}/p\mathbb{Z})$, $\det(h) \equiv 1 \pmod{p}$. On the other hand, $\det(h) \neq 1$ because $h \not\in \SL_2(\mathbb{Z}/p^2\mathbb{Z})$. Choose some $s \in \SL_2(\mathbb{Z}/p^2\mathbb{Z})$ such that $\varphi(s) = \varphi(h)$. Note that $sh^{-1} \in \ker \varphi$, and $\det(sh^{-1}) \neq 1$, and so $sh^{-1}$ is of the form $I + \begin{pmatrix} a & b \\ c & d \end{pmatrix} p$ where $a+d \neq 0$. Hence, $\ker \varphi \leq H$, in which case $H = \varphi^{-1}(\SL_2(\mathbb{Z}/p\mathbb{Z}))$.
	\end{proof}
\end{lemma}

\begin{proposition} \label{propSL}
	Suppose that $p > 3$. Let $H_1, H_2 \leq \GL_2(\mathbb{Z}/p^2\mathbb{Z})$. If $H_1$ and $H_2$ are locally conjugate with $\SL_2(\mathbb{Z}/p\mathbb{Z}) \leq \varphi(H_i)$, then $H_1 = H_2$.
	\begin{proof}
		Since $\varphi(H_1)$ and $\varphi(H_2)$ are locally conjugate by Proposition \ref{propNecessaryTwo}, they are conjugate, and in fact equal, by Theorem \ref{theoremLocConjGL2p}. Clearly, $\varphi(H_i) = \SL_2(\mathbb{Z}/p\mathbb{Z}) \rtimes D$, where $D$ is a subgroup of the group $\left\{ \begin{pmatrix} 1 & 0 \\ 0 & d \end{pmatrix} \in \GL_2(\mathbb{Z}/p\mathbb{Z}) \right\}$. Say that $\begin{pmatrix} 1 & 0 \\ 0 & d_0 \end{pmatrix}$ generates $D$. By Corollary \ref{corDiagMember} and since $I + \begin{pmatrix} 0 & 1 \\ 0 & 0 \end{pmatrix} p, I + \begin{pmatrix} 0 & 0 \\ 1 & 0 \end{pmatrix} p \in H_i$, $H_i$ has $\begin{pmatrix} 1 & 0 \\ 0 & d_0 \end{pmatrix}$. Therefore, $H_i = (H_i \cap \varphi^{-1} (SL_2(p))) \rtimes \left \langle \begin{pmatrix} 1 & 0 \\ 0 & d_0 \end{pmatrix} \right \rangle$. By Lemma \ref{lemmaSubImSL2}, $H_i \cap \varphi^{-1} (\SL_2(\mathbb{Z}/p\mathbb{Z})) = \SL_2(\mathbb{Z}/p\mathbb{Z})$ or $\varphi^{-1}(\SL_2(\mathbb{Z}/p\mathbb{Z}))$. $H_1 \cap \varphi^{-1}(\SL_2(\mathbb{Z}/p\mathbb{Z}))$ and $H_2 \cap \varphi^{-1}(\SL_2(\mathbb{Z}/p\mathbb{Z}))$ are thus equal, and so $H_1 = H_2$. 
	\end{proof}
\end{proposition}

For $p = 3$, the subgroups of $\GL_2(\mathbb{Z}/p^2\mathbb{Z})$ whose images under $\varphi$ are $\SL_2(\mathbb{Z}/p\mathbb{Z})$ are conjugate to one of the following:
\begin{enumerate}
	\item $\left \langle \begin{pmatrix} 7 & 6 \\ 4 & 4 \end{pmatrix}, \begin{pmatrix} 7 & 4 \\ 6 & 4 \end{pmatrix} \right \rangle$
	\item $\left \langle \begin{pmatrix} 1 & 1 \\ 0 & 1 \end{pmatrix}, \begin{pmatrix} 1 & 0 \\ 7 & 1 \end{pmatrix} \right \rangle = \SL_2(\mathbb{Z}/p^2\mathbb{Z})$
	\item $\left \langle \begin{pmatrix} 4 & 1 \\ 0 & 1 \end{pmatrix}, \begin{pmatrix} 7 & 0 \\ 4 & 1 \end{pmatrix} \right \rangle$
	\item $\left \langle \begin{pmatrix} 1 & 1 \\ 0 & 1 \end{pmatrix}, \begin{pmatrix} 4 & 0 \\ 1 & 1 \end{pmatrix} \right \rangle = \varphi^{-1}(\SL_2(\mathbb{Z}/p\mathbb{Z}))$ 
	\item $\left \langle \begin{pmatrix} 7 & 6 \\ 1 & 7 \end{pmatrix}, \begin{pmatrix} 7 & 4 \\ 6 & 4 \end{pmatrix} \right \rangle$
	\item $\left \langle \begin{pmatrix} 1 & 6 \\ 7 & 7 \end{pmatrix}, \begin{pmatrix} 4 & 7 \\ 6 & 4 \end{pmatrix} \right \rangle$.
\end{enumerate}
No two distinct subgroups among these are locally conjugate. Therefore, the same result as Proposition \ref{propSL} holds for $p = 3$. 

\section{The Exceptional cases}
This section determines local conjugacy in $\GL_2(\mathbb{Z}/p^2\mathbb{Z})$ for the remaining cases. 

\begin{lemma} \label{lemmaExceptionalKer}
	Let $H$ be a subgroup of $\GL_2(\mathbb{Z}/p^2\mathbb{Z})$. If the image of $H$ in $\PGL_2(\mathbb{Z}/p^2\mathbb{Z})$ is isomorphic to $A_4, S_4$ or $A_5$, then $H \cap \ker \varphi$ is one of the following:
	\begin{enumerate}
		\item $\left \langle I \right \rangle$
		\item $\left \langle I + \begin{pmatrix} 1 & 0 \\ 0 & 1 \end{pmatrix} p \right \rangle$
		\item $T$
		\item $\ker \varphi$
	\end{enumerate}
	\begin{proof}
		Let $\overline{H}$ denote the image of $H$ in $\PGL_2(\mathbb{Z}/p\mathbb{Z})$. Note that $\overline{H}$ contains a subgroup which is isomorphic to $A_4$. Therefore, $\overline{H}$ has an element $\overline{h_1}$ of order $3$ and an element $\overline{h_2}$ of order $2$ which do not commute. In fact, $\overline{h_1} \overline{h_2}$ has order $3$. Furthermore, $\overline{H}$ is generated by $\overline{h_1}$ and $\overline{h_2}$. Choose $h_1,h_2 \in \varphi(H)$ such that the images of $h_1$ and $h_2$ in $\overline{H}$ are $\overline{h_1}$ and $\overline{h_2}$ respectively. By Proposition \ref{propSubGpGL2p}, $\varphi(H)$ has no element of order $p$. In particular, $\varphi(H)$ has no element which is conjugate to a matrix of the form $\begin{pmatrix} w & 1 \\ 0 & w \end{pmatrix}$ where $w \in (\mathbb{Z}/p\mathbb{Z})^\times$. Thus, $h_1$ is conjugate to an element in $C_s(p)$ or an element in $C_{ns}(p)$. Assume for the rest of the proof that $h_1$ is conjugate to an element in $C_s(p)$. The case where $h_1$ is conjugate to an element of $C_{ns}(p)$ works similarly. Replace $H$ with a conjugate so that $h_1$ is replace with a matrix of the form $\begin{pmatrix} w & 0 \\ 0 & z \end{pmatrix} \in C_s(p)$. Note that $w \neq \pm z$. Express $h_2$ as $h_2 = \begin{pmatrix} \alpha & \beta \\ \gamma & \delta \end{pmatrix}$. Compute
		\begin{align*}
			h_2^2 = \begin{pmatrix} \alpha^2 + \beta \gamma & (\alpha + \delta) \beta \\ (\alpha + \delta) \gamma & \beta \gamma + \delta^2 \end{pmatrix}.
		\end{align*}
		Since $\overline{h_2}$ is an element of $\PGL_2(\mathbb{Z}/p\mathbb{Z})$ of order $2$, $\alpha^2 = \delta^2$ and $(\alpha + \delta) \beta = (\alpha + \delta) \gamma = 0$. If $\alpha = \delta$, then $h_2 = \begin{pmatrix} \alpha & 0 \\ 0 & \alpha \end{pmatrix}$ or $\begin{pmatrix} 0 & \beta \\ \gamma & 0 \end{pmatrix}$. However, $\overline{H}$ would be cyclic or dihedral, which is a contradiction to Proposition \ref{propSubGpGL2p}. Hence, $\alpha = -\delta \neq 0$ and we express $h_2$ as $h_2 = \begin{pmatrix} \alpha & \beta \\ \gamma & -\alpha \end{pmatrix}$, where not both $\beta$ and $\gamma$ are $0$. Suppose that $\gamma = 0$. Compute
		\begin{align*}
			(h_1h_2)^3 &= \begin{pmatrix} \alpha w & \beta w \\ 0 & -\alpha z \end{pmatrix}^3 \\
								&= \begin{pmatrix} \alpha^3 w^3 & * \\ 0 & -\alpha^3 z^3 \end{pmatrix}.
		\end{align*}
		Since $\overline{h_1}\overline{h_2}$ has order $3$, $\alpha^3 w^3 = - \alpha^3 z^3$. However, $w^3 = z^3$ because $\overline{h_1}$ has order $3$, which is a contradiction. Hence, $\gamma \neq 0$ and similarly, $\beta \neq 0$.  
		\par
		An argument similar to that of Lemma \ref{lemmaCartanKer} shows that $I + \begin{pmatrix} a & b \\ c & d \end{pmatrix} \in H$ if and only if $I + \begin{pmatrix} a & 0 \\ 0 & d \end{pmatrix} p, I + \begin{pmatrix} 0 & b \\ c & 0 \end{pmatrix} p \in H$. Moreover, if $I + \begin{pmatrix} 0 & b \\ c & 0 \end{pmatrix} p \in H$ with $b$ and $c$ nonzero, then $I + \begin{pmatrix} 0 & b \\ c & 0 \end{pmatrix} p$ and $h_1 \left( I + \begin{pmatrix} 0 & b \\ c & 0 \end{pmatrix} p \right) h_1^{-1} = I + \begin{pmatrix} 0 & b \frac{w}{z} \\ c \frac{z}{w} & 0 \end{pmatrix} p$ are linearly independent because $w \neq \pm z$, in which case $I + \begin{pmatrix} 0 & 1 \\ 0 & 0 \end{pmatrix} p, I + \begin{pmatrix} 0 & 0 \\ 1 & 0 \end{pmatrix} p \in H$.  \par
		Suppose that $H$ has an element of the form $I + \begin{pmatrix} 0 & b \\ c & 0 \end{pmatrix} p$ where $b$ or $c$ is nonzero. Compute
		\begin{align*}
			h_2 \left( I + \begin{pmatrix} 0 & b \\ c & 0 \end{pmatrix} p \right) h_2^{-1} = I + \frac{\begin{pmatrix} c \alpha \beta + b \alpha \gamma & c \beta^2 - b \alpha^2 \\ -c \alpha^2 + b \gamma^2 & -c \alpha \beta - b \alpha \gamma \end{pmatrix} p}{-\alpha^2 - \beta \gamma}.
		\end{align*}
		$c\beta^2 - b \alpha^2$ and $-c\alpha^2 + b \gamma^2$ are in the same ratio as $b$ and $c$, i.e. $c(c\beta^2 - b \alpha^2) = b(-c\alpha^2 + b \gamma^2)$, when $b^2 \gamma^2 = c^2 \beta^2$. If this holds, then $b$ and $c$ are both nonzero. Otherwise, $I + \begin{pmatrix} 0 &b \\ c & 0 \end{pmatrix} p$ and $I + \frac{\begin{pmatrix} 0 & c \beta^2 - b \alpha^2 \\ -c \alpha^2 + b \gamma^2 & 0 \end{pmatrix} p}{-\alpha^2 - \beta \gamma}$ are linearly independent. In either case, $I + \begin{pmatrix} 0 & 1 \\ 0 & 0 \end{pmatrix} p, I + \begin{pmatrix} 0 & 0 \\ 1 & 0 \end{pmatrix} p \in H$. Moreover, 
		\begin{align*}
			h_2 \left( I + \begin{pmatrix} 0 & 1 \\ 0 & 0 \end{pmatrix} p \right) h_2^{-1} = I + \frac{ \begin{pmatrix} \alpha \gamma & -\alpha^2 \\ \gamma^2 & -\alpha \gamma \end{pmatrix} p}{-\alpha^2 - \beta \gamma}
		\end{align*}\par
		and since $\alpha$ and $\gamma$ are nonzero, $I + \begin{pmatrix} 1 & 0 \\ 0 & -1 \end{pmatrix} p \in H$. \par
		Suppose that there is some $I + \begin{pmatrix} a & 0 \\ 0 & d \end{pmatrix} p \in H$ such that $a \neq d$. Compute
		\begin{align*}
			h_2 \left( I + \begin{pmatrix} a & 0 \\ 0 & d \end{pmatrix} p \right) h_2^{-1} = I + \frac{\begin{pmatrix} a \alpha^2 + d \beta \gamma & a \alpha \beta - d \alpha \beta \\ a \alpha \gamma - d \alpha \gamma & a \beta \gamma - d \alpha^2 \end{pmatrix} p}{-\alpha^2 - \beta \gamma}.
		\end{align*}
		Since $\alpha \neq 0$ and $\beta$ and $\gamma$ are nonzero, $H$ has an element of the form $I + \begin{pmatrix} 0 & b \\ c & 0 \end{pmatrix} p$ where $b$ and $c$ are nonzero.  \par
		From the last two paragraphs, it is not difficult to see that $H \cap \ker \varphi$ is one of the four groups as claimed.
	\end{proof}
\end{lemma}

\begin{proposition}
	Let $H_1, H_2 \leq \GL_2(\mathbb{Z}/p^2\mathbb{Z})$ be locally conjugate and suppose that the image of $H_i$ in $\PGL_2(\mathbb{Z}/p\mathbb{Z})$ is isomorphic to $A_4$, $S_4$ or $A_5$. $H_1$ and $H_2$ are conjugate.
	\begin{proof}
		By Proposition \ref{propNecessaryOne}, $H_1 \cap \ker \varphi$ and $H_2 \cap \ker \varphi$ are locally conjugate. Lemma \ref{lemmaExceptionalKer} shows that $H_1 \cap \ker \varphi = H_2 \cap \ker \varphi$. Moreover, $\varphi(H_1)$ and $\varphi(H_2)$ are locally conjugate by Proposition \ref{propNecessaryTwo}. They must be conjugate by Proposition \ref{propSubGpGL2p} and Theorem \ref{theoremLocConjGL2p}. Replace $H_2$ with a conjugate so that $\varphi(H_1) = \varphi(H_2)$. All of the subgroups of $\ker \varphi$ which are listed in Lemma \ref{lemmaExceptionalKer} are normal in $\GL_2(\mathbb{Z}/p^2\mathbb{Z})$ and so $H_2 \cap \ker \varphi$ is preserved by the conjugation. \par
		Note that $|\varphi(H_i)|$ divides $12(p-1)$, $24(p-1)$ or $60(p-1)$ because the image of $H_i$ in $\PGL_2(\mathbb{Z}/p\mathbb{Z})$ is isomorphic to $A_4$, $S_4$ or $A_5$ and the kernel of the natural homomorphism $\GL_2(\mathbb{Z}/p\mathbb{Z}) \rightarrow \PGL_2(\mathbb{Z}/p\mathbb{Z})$ is $Z(p)$, which has order $p-1$. Therefore, if $p > 5$, then $p$ does not divide $|\varphi(H_i)|$ \cite[Table 2]{Sutherland} shows that the image of $H_i$ in $\PGL_2(\mathbb{Z}/p\mathbb{Z})$ cannot be isomorphic to $A_4$, $S_4$ or $A_5$ if $p = 3$ and cannot be isomorphic to $A_5$ if $p = 5$. In any case, $p$ does not divide $|\varphi(H_i)|$ and so $H_1$ and $H_2$ are conjugate by Proposition \ref{propSchurZassenhaus}.
	\end{proof}
\end{proposition}

\section{Conclusion}

Theorem \ref{theoremConclusion} summarizes the categorization of pairs of nontrivially locally conjugate subgroups of $\GL_2(\mathbb{Z}/p^2\mathbb{Z})$.

\begin{theorem}\label{theoremConclusion}
	The pairs of nontrivially locally conjugate subgroups of $\GL_2(\mathbb{Z}/p^2\mathbb{Z})$ are, up to conjugation, those listed in Propositions \ref{propCartan} and \ref{propBorel} and for $p = 3$, the ones described at the end of Section \ref{sectionBorel}. 
\end{theorem}

\section*{Acknowledgements}
 
I would like to thank my mentor Atticus Christensen for his invaluable insight on approaching the main problem, his daily help on the complicated details towards solving it, and his guidance in constructing this paper. I would also like to thank Professor David Jerison, Professor Ankur Moitra, and Dr. Andrew Sutherland for their encouragement and experienced advice on working with mathematics. In particular, Dr. Sutherland motivated this project and guided me on it as my supervisor for the MIT Undergraduate Research Opportunities Program. I give my final thanks to Dr. Slava Gerovitch and MIT Mathematics for making the Summer Program in Undergraduate Research possible.


\end{document}